\documentclass[12pt,notitlepage]{article}
\title{Limit groups with respect to Thompson's group $F$ and other finitely generated groups}
\author{Roland Zarzycki\footnote{The author was partially supported by the national grant 3707/B/H03/2009/36.}}
\date{27.06.2009}
\usepackage[latin2]{inputenc}
\usepackage{amssymb}
\usepackage{indentfirst}
\usepackage{latexsym}
\begin{document}

\maketitle

\begin{abstract}

 Let $F$ be the (Thompson's) group $\langle x_{0}, x_{1}| [x_{0}x_{1}^{-1}, x_{0}^{-i}x_{1}
x_{0}^{i}], i=1,2 \rangle$. We study the structure of $F$-limit groups. Let $G_{n}= \langle
y_{1},\ldots , y_{m},$ $x_{0},x_{1} \vert [x_{0}x_{1}^{-1},x_{0}^{-1}x_{1}x_{0}],[x_{0}x_{1}
^{-1},x_{0}^{-2}x_{1}x_{0}^{2}], y_{j}^{-1}g_{j,n}(x_{0},x_{1}), 1\leq j\leq m \rangle$,
where $g_{j,n}(x_{0},x_{1})\in F$, $n\in\mathbb{N}$, be a family of groups marked by $m+2$
elements. If the sequence $(G_{n})_{n<\omega}$ is convergent in the space of marked
groups and $G$ is the corresponding limit we say that $G$ is an $F$-limit group. Primarily
the paper is devoted to the study of $F$-limit groups.\parskip0pt

 The results are based on some theorems concerning laws with parameters in $F$. In
particular several constructions of such laws are given. On the other hand we formulate
some very general conditions on words with parameters $w(y,a_{1},\ldots ,a_{n})$ over
$F$ which guarantee that the inequality $w(y,\bar{a})\neq 1$ has a solution in $F$.\parskip0pt

 Some of the results are of a more general nature and can be applied to study limit groups with
respect to other finitely generated groups and classes of finitely generated groups, in particular
to the case of the Grigorchuk group.

\end{abstract}

\newpage

\tableofcontents

\begin{flushleft}
\hspace*{16pt}\textbf{References}\hspace*{294pt} \textbf{49}
\end{flushleft}

\newtheorem{theo}{Theorem}[section]
\newtheorem{lemm}[theo]{Lemma}
\newtheorem{coro}[theo]{Corollary}
\newtheorem{defi}[theo]{Definition}
\newtheorem{prop}[theo]{Proposition}
\newtheorem{fact}[theo]{Fact}
\newtheorem{rema}[theo]{Remark}
\newtheorem{exam}[theo]{Example}
\newtheorem{ques}[theo]{Question}
\newtheorem{clai}[theo]{Claim}

\section{Introduction}

\subsection{Outline}

 The notion of \emph{limit groups} was introduced by Z. Sela in his work on characterization
of elementary equivalence of free groups \cite{Sel}. The idea has been extended in the
paper of C. Champetier and V. Guirardel \cite{GC}, where the authors look at \emph{limit
groups} as limits of convergent sequences in spaces of \emph{marked groups}. They have
given a description of Sela's limit groups in these terms (with respect to the class of free
groups). This approach has been also aplied by L. Guyot and Y. Stalder \cite{SG} to the
class of Baumslag-Solitar groups and recently by L. Guyot \cite{G} to the class of dihedral
groups. \parskip0pt

 Thompson's group $F$ has remained one of the most interesting objects in geometric group
theory. This makes natural to consider limits with respect to $F$. Let $F$ be given by the
standard presentation $\langle x_{0}, x_{1}\ |\ [x_{0}x_{1}^{-1}, x_{0}^{-i}x_{1}x_{0}^{i}], i=1,2
\rangle$. Now let
$$G_{n}= \Big\langle x_{0},x_{1}, y_{1},\ldots ,y_{m}\ \Big|\ [x_{0}x_{1}^{-1},x_{0}^{-i}x_{1}x_{0}^{i}], y
	_{j}^{-1}g_{j,n}(x_{0},x_{1}), i\in\{1, 2\} , j\leq m\Big\rangle ,$$
where $g_{j,n}(x_{0},x_{1})\in F$, $n\in\mathbb{N}$, be a family of groups isomorphic to $F$
and marked by $m+2$ elements. If the sequence $(G_{n})_{n<\omega}$ is convergent in the
space of marked groups and $G$ is the corresponding limit, we say that $G$ is an $F$-limit
group. Our paper is devoted to a description of $F$-limit groups. \parskip0pt

 In our study we employ mainly two technical tools. The first one is Theorem \ref{ab}, which
gives a sufficient condition for the existence of a solution of a given inequality in the (more
general) context of a group having \emph{hereditarily separating action}. We develop this
approach in Theorem \ref{uab}, which deals with a finite set of inequalities and uses weaker
assumptions.\parskip0pt

 The second tool is Proposition \ref{lwc2} (also Proposition \ref{lwc4}), which enables to
construct \emph{laws with parameters} in the Thompson's group $F$. Note that in 1985
Brin and Squier \cite{BS} showed that Thompson's group $F$ does not satisfy any law
(also see Abert's paper \cite{A} for a shorter proof). However, we will show that there are
certain non-trivial words with constants over $F$ (which will be later called \emph{laws with
parameters}), which are equal to identity for each evaluation in $F$.\parskip0pt

 One of the most important consequences of the above theorems is the proof that the
HNN-extension of the form $F\ast _{H_{1}}$ for a certain infinitely generated subgroup $H_{1}
<F$ can be realized as an $F$-limit group. It is the subject of Theorem \ref{h1}. Another
application gives Theorem \ref{free}, which shows that no free product of the form $F\ast G$
occurs as a limit of a sequence $(G_{n})_{n<\omega}$ as above. Similarly by Theorem \ref{mt},
in the case of convergent sequence $(G_{n})_{n<\omega}$ marked by $3$ elements we
cannot get an HNN-extension over a finitely generated subgroup. Theorem \ref{dmt} gives an
anologous statement in the case of centralized HNN-extensions over infinitely generated
subgroups satisfying certain technical condition.\parskip0pt

 This paper have been submitted as the doctoral dissertation at the Univeristy of Wroc{\l}aw. It
was arranged as an integral whole and some of the results draw their meaning from the context
of the others. This makes the text difficult to split. Therefore although it may seem to be too long,
we decided to leave it in its integral form.\parskip0pt

 Some results concerning the existence of laws with parameters have been already given in our
previous paper (\cite{Zar}). However, in the present paper we give new construction of such laws,
simplify proofs and generalize applications of the previously obtained results.\parskip18pt

 On a personal note, I would like to express my deepest gratitude to my supervisor, professor
Alexander Iwanow, for his guidance, encouragement and priceless remarks.\parskip0pt

\subsection{Preliminaries}

 A \emph{marked group} $(G,S)$ is a group G with a distinguished set of generators $S =
(s_{1}, s_{2}, \ldots , s_{n})$. For fixed $n$, let $\mathcal{G} _{n}$ be the set of all
$n$-generated groups marked by $n$ generators (up to isomorphism of marked groups).
Following \cite{GC} we put certain metric on $\mathcal{G} _{n}$. We will say that two marked
groups $(G,S), (G',S')\in\mathcal{G} _{n}$ are at distance less or equal to $e^{-R}$ if they
have exactly the same relations of length at most $R$. The set $\mathcal{G} _{n}$ equiped
with this metric is a compact space \cite{GC}. \emph{Limit groups} are simply limits of
convergent sequences in this metric space.

\begin{defi} Let $G$ be an $n$-generated group. A marked group in $\mathcal{G} _{n}$ is
	a \emph{limit group with respect to $G$} if it is a limit of marked groups each
	isomorphic to $G$.
\end{defi}

 To introduce the Thompson's group $F$ we will follow \cite{CFP}.

\begin{defi} Thompson's group $F$ is the group given by the following infinite group
presentation:
	$$\Big\langle x_{0}, x_{1}, x_{2}, \ldots \Big| x_{j}x_{i}=x_{i}x_{j+1} (i<j) \Big\rangle .$$
\end{defi}

 In fact $F$ is finitely presented:
	$$ F = \Big\langle x_{0}, x_{1}\ \Big|\ [x_{0}x_{1}^{-1}, x_{0}^{-i}x_{1}x_{0}^{i}], i=1,2
		\Big\rangle . $$
 Every non-trivial element of $F$ can be uniquely expressed in the normal form:
	$$ x_{0}^{b_{0}}x_{1}^{b_{1}}x_{2}^{b_{2}}\ldots x_{n}^{b_{n}}x_{n}^{-a_{n}}\ldots
	x_{2}^{-a_{2}}x_{1}^{-a_{1}}x_{0}^{-a_{0}},$$
where $n$, $a_{0}, \ldots , a_{n}$, $b_{0}, \ldots , b_{n}$ are non-negative integers such
that: \\
	i) exactly one of $a_{n}$ and $b_{n}$ is nonzero; \\
	ii) if $a_{k}>0$ and $b_{k}>0$ for some integer $k$ with $0\leq k < n$, then
	$a_{k+1}>0$ or $b_{k+1}>0$. \parskip0pt

 We study properties of limit groups with respect to $F$. For this purpose let us consider
a sequence, $(g_{i,n})_{n<\omega}$, $1\leq i\leq m$, of elements taken from the group
$F$ and the corresponding sequence of limit groups marked by $m+2$ elements, $G_{n}
=(F,(x_{0},x_{1},g_{1,n},\ldots ,g_{m,n}))$, $n\in\mathbb{N}$, where $x_{0}$ and $x_{1}$
are the standard generators of $F$. Assuming that such a sequence is convergent in the
space of groups marked by $m+2$ elements, denote by
$$G = \Big(\Big\langle x_{0}, x_{1}, g_{1}\ldots , g_{m}\ \Big|\ R_{F} \cup R_{G} \Big\rangle,
	\Big( x_{0}, x_{1}, g_{1},\ldots , g_{m}\Big)\Big)$$
the limit group formed in that manner; here $x_{0}$, $x_{1}$ are "limits" of constant
sequences $(x_{0})_{n<\omega}$ and $(x_{0})_{n<\omega}$, $g_{i}$ is the "limit" of $(g_{i,
n})_{n<\omega}$ for $1\leq i\leq m$, $R_{F}$ and $R_{G}$ refer respectively to the set of
standard relations taken from $F$ and the set (possibly infinite) of new relations. \parskip0pt

 It has been shown in \cite{GC} that in the case of free groups some standard constructions
can be obtained as limits of free groups. For example, it is possible to get $\mathbb{Z}
^{k}$ as a limit of $\mathbb{Z}$ and $\mathbb{F}_{k}$ as a limit of $\mathbb{F} _{2}$.  On
the other hand, the direct product of $\mathbb{F} _{2}$ and $\mathbb{Z}$ cannot be
obtained as a limit group with respect to $\mathbb{F} _{2}$. HNN-extensions often occur in
the class of limit groups (with respect to free groups). For example, the following groups
are the limits of convergent sequences in the space of free groups marked by three
elements: the free group of rank $3$, the free abelian group of rank $3$ or a HNN-extension
over a cyclic subgroup of the free group of rank $2$ (\cite{FGM}). All non-exceptional
surface groups form another broad class of interesting examples (\cite{BaB},
\cite{BaG}).\parskip0pt

 In the case of Thompson's group the situation is not so clear. Since the centrum of $F$ is
trivial, we can only say that any sequence of groups of the form given above cannot
converge to any direct product with the whole group as an $F$-limit group. It is worth noting
that the case when we do not fix the generators $x_{0}$, $x_{1}$ has been studied before.
Recently M. Brin have proved that it is possible to obtain the non-abelian free group of rank
$2$ as the limit of marked sequences of the form $((F, (g_{0,n}, g_{1,n})))_{n<\omega}$
for certain $(g_{0,n})_{n<\omega}$ and $(g_{1,n})_{n<\omega}$ from $F$ (see \cite{Bri}
for details). A similar result have been announced by J. Taback during the conference on
Geometric and Asymptotic Group Theory with Applications, which took place in Hoboken
NJ (\cite{Ast}).\parskip0pt

 There are many geometric interpretations of $F$, but here we use the following one.
Consider the set of all strictly increasing continuous piecewise-linear functions from the
closed unit interval onto itself. Then the group $F$ is realized by the set of all such
functions, which are differentiable except at finitely many dyadic rational numbers and such
that all slopes (deriviatives) are integer powers of 2. The corresponding group operation is
just the composition. For the further reference it will be usefull to give an explicit form of
the generators $x_{n}$, for $n\geq 0$, in terms of piecewise-linear functions:

$$ x_{n}(t) = \left\{ \begin{array}{ll} t & \textrm{, $t\in [0,\frac{2^{n}-1}{2^{n}} ]$} \\
	\frac{1}{2}t + \frac{2^{n}-1}{2^{n+1}} & \textrm{, $t\in [\frac{2^{n}-1}{2^{n}}, \frac{2^{n+1}-1}
	{2^{n+1}} ]$} \\ t - \frac{1}{2^{n+2}} & \textrm{, $t\in [\frac{2^{n+1} -1}{2^{n+1}},
	\frac{2^{n+2}-1}{2^{n+2}}]$} \\	2t-1 & \textrm{, $t\in [\frac{2^{n+2}-1}{2^{n+2}},
	1]$} \end{array}\right.$$ \parskip0pt

 For any dyadic subinterval $[a,b]\subset [0,1]$, let us consider the set, of elements in $F$,
which are trivial on its complement, and denote it by $F _{[a,b]}$. We know that it forms a
subgroup of $F$, which is isomorphic to the whole group. Let us denote its standard infinite
set of generators by $x_{[a,b],0}, x_{[a,b],1}, x_{[a,b],2}, \ldots$, where for $n\geq 0$ we have:
$$x_{[a,b],n}(t) = \left\{ \begin{array}{ll} t & \textrm{, $t\in [0,a+\frac{(2^{n}-1)(b-a)}{2^{n}}]$} \\
	\frac{1}{2}t + \frac{1}{2}(a+\frac{2^{n}-1}{2^{n}}) & \textrm{, $t\in [a+\frac{(2^{n}-1)(b-a)}{2^{n}},
			a+\frac{(2^{n+1}-1)(b-a)}{2^{n+1}}]$} \\
	 t - \frac{b-a}{2^{n+2}} & \textrm{, $t\in [a+\frac{(2^{n+1}-1)(b-a)}{2^{n+1}},
			a+\frac{(2^{n+2}-1)(b-a)}{2^{n+2}}]$} \\
	2t - b & \textrm{, $t\in [a+\frac{(2^{n+2}-1)(b-a)}{2^{n+2}},b]$} \\
	t & \textrm{, $t\in [b,1]$} \end{array}\right.$$
 Moreover, if $\iota _{[a,b]}$ denotes the natural isomorphism between $F$ and $F_{[a,b]}$
sending $x_{n}$ to $x_{[a,b],n}$ for all $n\geq 0$, then for any $f\in F$ by $f_{[a,b]}$ we denote
the element $\iota _{[a,b]} (f)\in F_{[a,b]}<F$. \parskip0pt

 Let us consider an arbitrary element $g$ in $F$ and treat it as a piecewise-linear
homeomorphism of the interval $[0,1]$. Let $supp(g)$ be the set $\{ x\in [0,1] : g(x) \neq x \}$
and $\overline{supp}(g)$ the topological closure of $supp(g)$. We will call each point from the
set $P_{g} = (\overline{supp}(g)\setminus supp(g))\cap\mathbb{Z} [\frac{1}{2}]$ a \emph{dividing
point} of $g$. This set is obviously finite and thus we get a finite subdivision of $[0,1]$ of the
form $[0=p_{0}, p_{1}], [p_{1}, p_{2}], \ldots , [p_{n-1}, p_{n}=1]$ for some natural $n$. It is
easy to see that $g$ can be presented as $g = g_{1}g_{2}\ldots g_{n}$, where $g_{i}\in F
_{[p_{i-1}, p_{i}]}$ for each $i$. Since $g$ can act trivially on some of these subintervals, some
of the elements $g_{1},\ldots ,g_{n}$ may be trivial. We call the set of all non-trivial elements
from $\{ g_{1},\ldots ,g_{n}\}$ the \emph{defragmentation} of $g$.\parskip0pt

\begin{fact}\emph{(Corollary 15.36 in \cite{GS}, Proposition 3.2 in \cite{KM})}\label{GSf}
	The centralizer of any element $g\in F$ is the direct product of finitely many cyclic
	groups and finitely many groups isomorphic to $F$.
\end{fact}

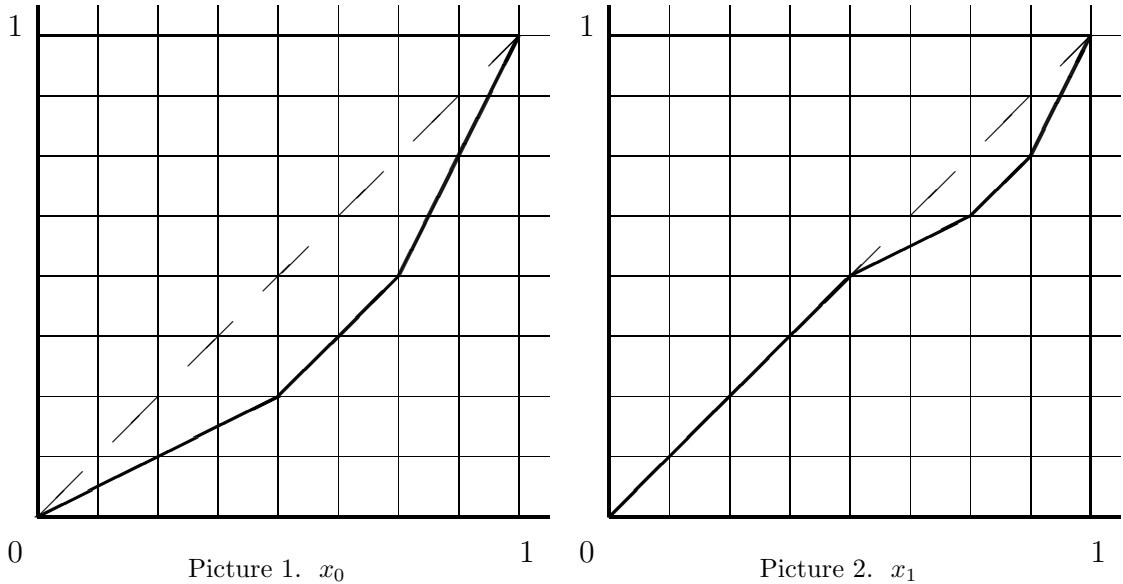
\setlength{\unitlength}{2mm}
\begin{picture}(45,45)
	\linethickness{0.05mm}
		\multiput(-4,8)(4,0){9}%
			{\line(0,1){34}}
		\multiput(-4,8)(0,4){9}%
			{\line(1,0){34}}

	\linethickness{0.375mm}
		\put(-4,8){\line(0,1){34}}
		\put(-4,8){\line(1,0){34}}

	\linethickness{0.25mm}
		\put(-4,40){\line(1,0){32}}
		\put(28,8){\line(0,1){32}}

	\linethickness{0.125mm}

		\put(-4,8){\line(1,1){3}}
		\put(1,13){\line(1,1){3}}
		\put(6,18){\line(1,1){3}}
		\put(11,23){\line(1,1){3}}
		\put(16,28){\line(1,1){3}}
		\put(21,33){\line(1,1){3}}
		\put(26,38){\line(1,1){2}}

		\put(-4.15,8){\line(2,1){16}}
		\put(-4.075,8){\line(2,1){16}}
		\put(-4,8){\line(2,1){16}}
		\put(-3.925,8){\line(2,1){16}}
		\put(-3.85,8){\line(2,1){16}}

		\put(11.9,16){\line(1,1){8}}
		\put(11.95,16){\line(1,1){8}}
		\put(12,16){\line(1,1){8}}
		\put(12.05,16){\line(1,1){8}}
		\put(12.1,16){\line(1,1){8}}

		\put(19.9,24){\line(1,2){8}}
		\put(19.95,24){\line(1,2){8}}
		\put(20,24){\line(1,2){8}}
		\put(20.05,24){\line(1,2){8}}
		\put(20.1,24){\line(1,2){8}}

	\put(6,4){{\footnotesize Picture 1.\ \ $x_{0}$}}
	\put(-6,5){$0$}
	\put(28,5){$1$}
	\put(-6,40){$1$}

	\linethickness{0.05mm}
		\multiput(34,8)(4,0){9}%
			{\line(0,1){34}}
		\multiput(34,8)(0,4){9}%
			{\line(1,0){34}}

	\linethickness{0.375mm}
		\put(34,8){\line(0,1){34}}
		\put(34,8){\line(1,0){34}}

	\linethickness{0.25mm}
		\put(34,40){\line(1,0){32}}
		\put(66,8){\line(0,1){32}}

	\linethickness{0.125mm}

		\put(34,8){\line(1,1){3}}
		\put(39,13){\line(1,1){3}}
		\put(44,18){\line(1,1){3}}
		\put(49,23){\line(1,1){3}}
		\put(54,28){\line(1,1){3}}
		\put(59,33){\line(1,1){3}}
		\put(64,38){\line(1,1){2}}

		\put(33.9,8){\line(1,1){16}}
		\put(33.95,8){\line(1,1){16}}
		\put(34,8){\line(1,1){16}}
		\put(34.05,8){\line(1,1){16}}
		\put(34.1,8){\line(1,1){16}}

		\put(49.85,24){\line(2,1){8}}
		\put(49.925,24){\line(2,1){8}}
		\put(50,24){\line(2,1){8}}
		\put(50.075,24){\line(2,1){8}}
		\put(50.15,24){\line(2,1){8}}

		\put(57.9,28){\line(1,1){4}}
		\put(57.95,28){\line(1,1){4}}
		\put(58,28){\line(1,1){4}}
		\put(58.05,28){\line(1,1){4}}
		\put(58.1,28){\line(1,1){4}}

		\put(61.9,32){\line(1,2){4}}
		\put(61.95,32){\line(1,2){4}}
		\put(62,32){\line(1,2){4}}
		\put(62.05,32){\line(1,2){4}}
		\put(62.1,32){\line(1,2){4}}

	\put(44,4){{\footnotesize Picture 2.\ \ $x_{1}$}}
	\put(32,5){$0$}
	\put(66,5){$1$}
	\put(32,40){$1$}

\end{picture}

\begin{rema}\label{GS}
	\emph{Moreover if the element $g\in F$ has the defragmentation $g=g_1 \ldots
	g_n$, then some roots of the elements $g_1, \ldots , g_n$ are the generators of
	cyclic components of the decomposition of the centralizer above. The components
	of this decomposition, which are isomorphic to $F$, are just the groups of the form
	$F_{[a,b]}$, where $[a,b]$ is one of the subintervals $[p_{i-1},p_{i}]\subseteq [0,1]$,
	which are stabilized pointwise by $g$.}
\end{rema}

\begin{fact}\label{GSc}
	The center of $F$ is trivial.
\end{fact}

\begin{rema}
	\emph{Note that it follows from Fact \ref{GSc} that the center of any $F_{[a,b]}$ is trivial.}
\end{rema}

Generally, if we interpret the elements of $F$ as functions, the relations occuring in the
presentation of $F$, $[x_{0}x_{1}^{-1}, x_{0}^{-i}x_{1}x_{0}^{i}]$ for $i = 1,2$, have to assure, that
two functions, which have mutually disjoint supports except of finitely many points, commute. In
particular, these relations imply analogous relations for different $i>2$. According to the fact that
$x_{0}^{-i}x_{1}x_{0}^{i} = x_{i+1}$, we conclude that all the relations of the form $[x_{0}x_{1}^{-1},
x_{M}]$, $M>1$, hold in Thompson's group $F$. We often refer to these geometrical
observations. \parskip0pt

 The following remark and fact are another useful geometrical tools.

\setlength{\unitlength}{2mm}
\begin{picture}(45,45)
	\linethickness{0.05mm}
		\multiput(-4,8)(4,0){5}%
			{\line(0,1){34}}
		\multiput(16,8)(4,0){3}%
			{\line(0,1){2}}
		\multiput(16,18)(4,0){4}%
			{\line(0,1){24}}

		\multiput(-4,12)(0,4){2}%
			{\line(1,0){18}}
		\multiput(26,12)(0,4){2}%
			{\line(1,0){4}}
		\multiput(-4,20)(0,4){6}%
			{\line(1,0){34}}

	\linethickness{0.375mm}
		\put(-4,8){\line(0,1){34}}
		\put(-4,8){\line(1,0){34}}

	\linethickness{0.25mm}
		\put(-4,40){\line(1,0){32}}
		\put(28,8){\line(0,1){32}}

	\linethickness{0.125mm}

		\put(-4,8){\line(1,1){3}}
		\put(1,13){\line(1,1){3}}
		\put(6,18){\line(1,1){3}}
		\put(11,23){\line(1,1){3}}
		\put(16,28){\line(1,1){3}}
		\put(21,33){\line(1,1){3}}
		\put(26,38){\line(1,1){2}}

		\put(-4.05,8){\line(1,1){4}}
		\put(-4,8){\line(1,1){4}}
		\put(-3.95,8){\line(1,1){4}}

		\put(-0.075,12){\line(2,1){4}}
		\put(0,12){\line(2,1){4}}
		\put(0.075,12){\line(2,1){4}}

		\put(3.95,14){\line(1,1){2}}
		\put(4,14){\line(1,1){2}}
		\put(4.05,14){\line(1,1){2}}

		\put(5.95,16){\line(1,4){2}}
		\put(6,16){\line(1,4){2}}
		\put(6.05,16){\line(1,4){2}}

		\put(7.95,24){\line(1,1){2}}
		\put(8,24){\line(1,1){2}}
		\put(8.05,24){\line(1,1){2}}

		\put(9.925,26){\line(4,1){4}}
		\put(10,26){\line(4,1){4}}
		\put(10.075,26){\line(4,1){4}}

		\put(13.925,27){\line(2,1){2}}
		\put(14,27){\line(2,1){2}}
		\put(14.075,27){\line(2,1){2}}

		\put(15.95,28){\line(1,1){4}}
		\put(16,28){\line(1,1){4}}
		\put(16.05,28){\line(1,1){4}}

		\put(19.925,32){\line(2,1){4}}
		\put(20,32){\line(2,1){4}}
		\put(20.075,32){\line(2,1){4}}

		\put(23.95,34){\line(1,1){2}}
		\put(24,34){\line(1,1){2}}
		\put(24.05,34){\line(1,1){2}}

		\put(25.95,36){\line(1,2){2}}
		\put(26,36){\line(1,2){2}}
		\put(26.05,36){\line(1,2){2}}

		\put(0,12){\circle*{1}}
		\put(6.7,18.7){\circle{1}}
		\put(16,28){\circle*{1}}
		\put(20,32){\circle*{1}}
		\put(28,40){\circle*{1}}
		
	\put(-1,5){{\footnotesize Picture 3.\ \ \emph{Defragmentation} of}}
	\put(0,2.5){{\footnotesize an exemplary element $g\in F$}}
	\put(-6,5){$0$}
	\put(28,5){$1$}
	\put(-6,40){$1$}

	\put(15,14){\circle*{1}}
	\put(15,11.5){\circle{1}}
	\put(15,16){{\tiny Dividing points:}}
	\put(16,13.5){{\tiny dyadic}}
	\put(16,11){{\tiny non-dyadic}}
	\put(5,21){{\small $g_{1}$}}
	\put(21,37){{\small $g_{2}$}}

	\linethickness{0.05mm}
		\multiput(34,8)(4,0){2}%
			{\line(0,1){34}}
		\multiput(46,8)(8,0){3}%
			{\line(0,1){34}}
		\multiput(34,8)(0,4){9}%
			{\line(1,0){34}}

		\put(42,8){\line(0,1){0.5}}
		\put(42,11.5){\line(0,1){30.5}}
		\put(50,8){\line(0,1){8.5}}
		\put(50,19.5){\line(0,1){22.5}}
		\put(58,8){\line(0,1){16.5}}
		\put(58,27.5){\line(0,1){14.5}}
		\put(66,8){\line(0,1){24.5}}
		\put(66,35.5){\line(0,1){6.5}}

	\linethickness{0.375mm}
		\put(34,8){\line(0,1){34}}
		\put(34,8){\line(1,0){34}}

		\put(38,12){\line(0,1){16}}
		\put(38,12){\line(1,0){16}}
		\put(54,12){\line(0,1){16}}
		\put(38,28){\line(1,0){16}}

		\put(58,32){\line(0,1){8}}
		\put(58,32){\line(1,0){8}}
		\put(66,32){\line(0,1){0.5}}
		\put(66,35.5){\line(0,1){4.5}}
		\put(58,40){\line(1,0){8}}

	\linethickness{0.25mm}
		\put(34,40){\line(1,0){32}}
		\put(66,8){\line(0,1){24.5}}
		\put(66,35.5){\line(0,1){4.5}}

	\linethickness{2mm}
		\put(34.5,8){\line(0,1){4}}
		\put(35.5,8){\line(0,1){4}}
		\put(36.5,8){\line(0,1){4}}
		\put(37.5,8){\line(0,1){4}}

		\put(54.5,28){\line(0,1){4}}
		\put(55.5,28){\line(0,1){4}}
		\put(56.5,28){\line(0,1){4}}
		\put(57.5,28){\line(0,1){4}}

	\linethickness{0.125mm}

		\put(34,8){\line(1,1){32}}

		\put(34,8){\circle*{1}}
		\put(38,12){\circle*{1}}
		\put(54,28){\circle*{1}}
		\put(58,32){\circle*{1}}
		\put(66,40){\circle*{1}}

		\put(38,12){\line(2,1){4}}
		\put(38,12){\line(1,2){2}}

		\put(42,14){\line(1,1){2}}
		\put(40,16){\line(1,1){2}}

		\put(44,16){\line(1,4){2}}
		\put(42,18){\line(4,1){8}}

		\put(46,24){\line(1,1){2}}
		\put(50,20){\line(1,1){2}}

		\put(48,26){\line(4,1){4}}
		\put(52,22){\line(1,4){1}}

		\put(52,27){\line(2,1){2}}
		\put(53,26){\line(1,2){1}}

		\put(58,32){\line(2,1){4}}
		\put(58,32){\line(1,2){2}}

		\put(62,34){\line(1,1){2}}
		\put(60,36){\line(1,1){2}}

		\put(64,36){\line(1,2){2}}
		\put(62,38){\line(2,1){4}}

	\put(63,33.5){{\scriptsize $\langle\sqrt{g_{2}}\rangle$}}

	\put(55,26){{\scriptsize $F_{[\frac{5}{8},\frac{3}{4}]}$}}

	\put(47,17.5){{\scriptsize $\langle\sqrt{g_{1}}\rangle$}}

	\put(39,10){{\scriptsize $F_{[0,\frac{1}{4}]}$}}

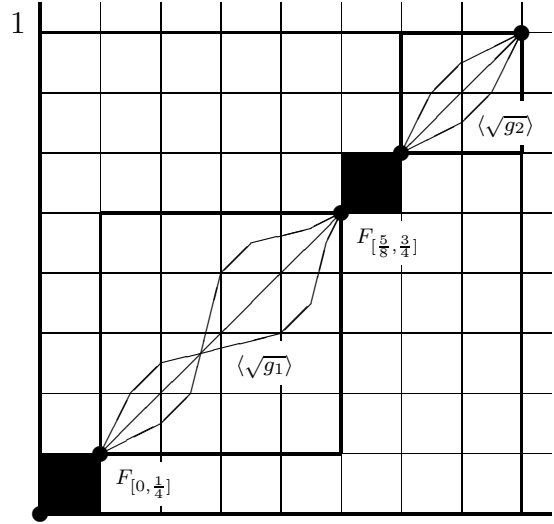
	\put(39,4){{\footnotesize Picture 4.\ \ $C_{F}(g)\cong F^{2}\oplus\mathbb{Z}^{2}$}}
	\put(32,5){$0$}
	\put(66,5){$1$}
	\put(32,40){$1$}

\end{picture}

\begin{rema}\label{rs}
 For any $f,f'\in F$ we have $supp(f^{f'})=(f')^{-1}(supp(f))$.
\end{rema}

\begin{fact}\emph{(Lemma 4.2 in \cite{CFP}, Lemma 2.4 in \cite{KM})}\label{CFP}
	If $0=x_{0}<x_{1}<x_{2}<\ldots <x_{n}=1$ and $0=y_{0}<y_{1}<y_{2}<\ldots <y_{n}=
	1$ are partitions of $[0,1]$ consisting of dyadic rational numbers, then there exists
	$f\in F$ such that $f(x_{i})=y_{i}$ for $i=0,\ldots , n$. Furthermore, if $x_{i-1}=y_{i-1}$
	and $x_{i}=y_{i}$ for some $i$ with $1\leq i\leq n$, then $f$ can be taken to be
	trivial on the interval $[x_{i-1},x_{i}]$.
\end{fact}

 Many examples, which occur in this paper, can be easily exposed using the \emph{rectangle
diagrams} introduced by W. Thurston (\cite{CFP}). Originally rectangle diagrams are drawn
horizontally, but for our purposes it will be more convenient to draw them vertically.\parskip0pt

 The information given by a function $f\in F$ is encoded in the partitions of the domain and
the range of $f$ determined by the breakpoints of $f$ and their images, respectively.
Therefore for a given element $f\in F$ we construct a rectangle with a left side, which is
identified with the domain of $f$, and a right side, which is identified with the range of $f$.
Next for every point $t$ on the left side where $f$ is not differentiable, we draw a line
segment from $t$ to $f(t)$ on the right side. We call the described picture the
\emph{rectangle diagram} of $f$. A particularly useful property of these diagrams is that
given two elements in $F$ we can easily draw the rectangle diagram of their composition by
juxtaposing their individual rectangle diagrams.\parskip0pt

\setlength{\unitlength}{2mm}
\begin{picture}(40,40)
	\linethickness{0.25mm}
		\put(2,6){\line(0,1){32}}
		\put(2,6){\line(1,0){8}}
		\put(10,6){\line(0,1){32}}
		\put(2,38){\line(1,0){8}}

			\put(2,22){\line(1,-1){8}}
			\put(2,30){\line(1,-1){8}}
			\put(5,4){\footnotesize $x_{0}$}

		\put(14,6){\line(0,1){32}}
		\put(14,6){\line(1,0){8}}
		\put(22,6){\line(0,1){32}}
		\put(14,38){\line(1,0){8}}

			\put(14,22){\line(1,0){8}}
			\put(14,30){\line(2,-1){8}}
			\put(14,34){\line(2,-1){8}}
			\put(17,4){\footnotesize $x_{1}$}

		\put(26,6){\line(0,1){32}}
		\put(26,6){\line(1,0){8}}
		\put(34,6){\line(0,1){32}}
		\put(26,38){\line(1,0){8}}

			\put(26,22){\line(1,-1){8}}
			\put(26,30){\line(1,-1){8}}
				\put(26,34){\line(2,-1){2}}
				\put(29,32.5){\line(2,-1){2}}
				\put(32,31){\line(2,-1){2}}
				\put(26,36){\line(4,-1){2}}
				\put(29,35.25){\line(4,-1){2}}
				\put(32,34.5){\line(4,-1){2}}
			\put(29,4){\footnotesize $x_{0}$}

		\put(34,6){\line(0,1){32}}
		\put(34,6){\line(1,0){8}}
		\put(42,6){\line(0,1){32}}
		\put(34,38){\line(1,0){8}}

			\put(34,22){\line(1,0){8}}
			\put(34,30){\line(2,-1){8}}
			\put(34,34){\line(2,-1){8}}
				\put(34,14){\line(1,0){2}}
				\put(37,14){\line(1,0){2}}
				\put(40,14){\line(1,0){2}}
			\put(37,4){\footnotesize $x_{1}$}

		\put(42,6){\line(0,1){32}}
		\put(42,6){\line(1,0){8}}
		\put(50,6){\line(0,1){32}}
		\put(42,38){\line(1,0){8}}

			\put(42,14){\line(1,1){8}}
			\put(42,22){\line(1,1){8}}
				\put(42,26){\line(4,3){2}}
				\put(45,28.25){\line(4,3){2}}
				\put(48,30.5){\line(4,3){2}}
				\put(42,30){\line(2,1){2}}
				\put(45,31.5){\line(2,1){2}}
				\put(48,33){\line(2,1){2}}
			\put(45,4){\footnotesize $x_{0}^{-1}$}

		\put(54,6){\line(0,1){32}}
		\put(54,6){\line(1,0){8}}
		\put(62,6){\line(0,1){32}}
		\put(54,38){\line(1,0){8}}

			\put(54,30){\line(1,0){8}}
			\put(54,34){\line(4,-1){8}}
			\put(54,36){\line(4,-1){8}}
				\put(54,22){\line(1,0){2}}
				\put(57,22){\line(1,0){2}}
				\put(60,22){\line(1,0){2}}
			\put(55,4){\footnotesize $x_{0}^{-1}x_{1}x_{0}$}

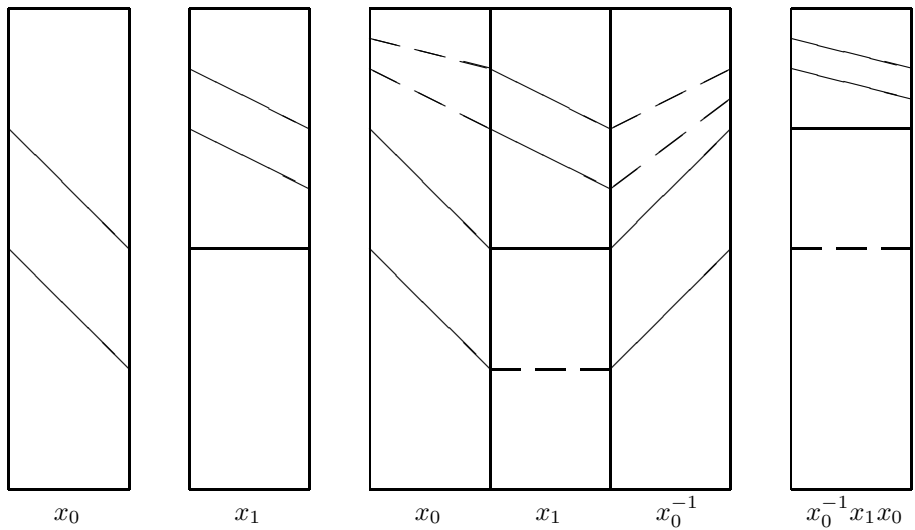
	\put(8,1){{\footnotesize Picture 5.\ \ Computing the \emph{rectangle diagram} of
		 $x_{0}^{-1}x_{1}x_{0}=x_{2}$}}

\end{picture}

 Since $F$ is a permutation group (on $[0,1]$) we sometimes apply standard terminology
and notation of the area of permutation groups (see \cite{H}).\parskip0pt

 Some of the methods, which we use in the case of Thompson's group, have more general
nature. In particular we apply them to the class of \emph{weakly branch groups}. To define
this class we will follow \cite{A}. We consider a finitely generated group $G$, which acts on
some rooted tree $T$. The vertices in $T$, which are at the same distance from the root are
said to be at the same \emph{level}. For any vertex $t\in T$ at the distance $k$ from the root,
all the vertices, which are at distance $1$ from $t$ and at distance $k+1$ from the root, are
called \emph{descendants} of $t$. We say that the action of the group $G$ on $T$ is
\emph{spherically transitive} if $G$ acts transitively on each level of $T$. For any vertex $t\in
T$ we define its \emph{rigid stabilizer} to be the set of all elements from $G$, which move
only descendants of $v$.\parskip0pt

\begin{defi}
	 A group $G$ is called \emph{weakly branch group} if it acts spherically transitively
	on some rooted tree $T$ so that the rigid stabilizer of every vertex is non-trivial.
\end{defi}

 The \emph{boundary} of a tree $T$, denoted by $\partial T$, consists of the infinte branches
starting at the root. A weakly branch group acts on the boundary $\partial T$ as well. The
class of weakly branch groups includes many groups with interesting properties. In particular,
it contains the first Grigorchuk group (see \cite{Gr} and \cite{Gri}), which is given by the following
presentation:
$$\Big\langle a,b,c,d\ \Big|\ 1=a^{2}=b^{2}=c^{2}=bcd=\sigma ^{k}((ad)^{4})=\sigma ^{k}((adacac)
	^{4}), k=0, 1,\ldots\Big\rangle ,$$
where the substitution $\sigma$ is defined by
$$\sigma :=\left\{ \begin{array}{ll} a\to aca\\
	b\to d\\
	c\to b\\
	d\to c.\end{array}\right.$$

 For any group $G$ acting on some set $X$ and any subset $A\subseteq X$, by $stab_{G}(A)$
we denote the set of all elements from $G$, which stabilize $A$ pointwise.\parskip18pt

 The structure of the paper is as follows. In Section 2 we consider the sets of inequalities in
groups having \emph{hereditarily separating} action. In Section 3 we apply the obtained
results to investigate the partial (relative) convergence (with respect to a fixed class of words)
of sequences of marked groups, in particular with respect to Thompson's group $F$ and
Grigorchuk group. In Section 4 we study the inequalities and laws with constants in the case
of Thompson's group $F$. Next we apply obtained results to show that among $F$-limit groups
marked by $x_{0}$, $x_{1}$ and some additional markers there are no free products of $F$
with any non-trivial group. We also prove that no non-trivial HNN-extensions over a finitely
generated subgroup occur as such a limit of $F$ marked by three elements.\parskip0pt

\section{Inequalities in groups with hereditarily separating action}

 In this section we prove some general results concerning inequalities over groups having
hereditarily separating actions in the sense of the paper of M. Abert \cite{A}. Then we apply
these results to limits of Thompson's group.

\subsection{Solving an inequality in the case of one oscillating word}

\begin{defi}
	Let $G$ be a permutation group acting on an infinite set $X$. We say that $G$
	\emph{separates} $X$ if, for any finite subset $Y\subset X$, the pointwise stabilizer
	$stab _{G}(Y)$ does not stabilize any point outside $Y$. \parskip0pt

	Assume that $X$ is a metric space, $G$ consists of homeomorphisms of $X$ and
	the set of fixed points, $Fix(G)$, is finite. We say that $G$ \emph{hereditarily separates}
	$X$ if, for any open and infinite subset $Z\subseteq X$ and for any finite subset $Y\subset
	Z$, the subgroup $stab _{G}((X\setminus	Z)\cup Y)$ does not stabilize any point from
	$Z\setminus (Y\cup Fix(G))$.
\end{defi}

\begin{rema} \label{rema}
	 \emph{ (taken from \cite{A})\ For any separating action of a group $G$ on $X$ and any
	finite $Y\subset X$ the orbits of the action of the pointwise stabilizer of $Y$ on $X
	\setminus Y$ are infinite. Indeed, let $X'$ be a non-trivial orbit of the action of $stab _{G}
	(Y)$. If $X'$ is finite consider $stab _{G}(Y')$, where $Y':=Y\cup(X'\setminus\{ x\})$ for
	some $x\in X'$. Since the set $X'\setminus\{ x\}$ is fixed by the action of $stab _{G}(Y')$,
	$x'$ has to be mapped onto itself. This contradicts separability of the action of
	$G$.} \parskip0pt

	 \emph{ As a corollary we see that for a hereditarily separating action of $G$ on $X$ and
	an open and infinite subset $Z\subseteq X$, the action of the stabilizer $stab _{G}(X\setminus
	Z)$ on $Z\setminus Fix(G)$ has only infinite orbits.}
\end{rema}

\begin{exam} \label{f}
	\emph{ Observe that Thompson's group $F$ is hereditarily separating with respect to its
	standard action on $[0,1]$. For this purpose suppose that $Z$ is and open subset of $[0,
	1]$, $\{ y_{1}, y_{2}, \ldots y_{s}\} :=Y\subset Z$ and $t\in Z\setminus (Y\cup \{ 0, 1\})$.
	There is some non-trivial dyadic segment $[p,q]\subseteq Z$ containing $t$ such that $Y
	\cap [p,q]=\emptyset$. Obviously $x_{[p,q],0}\in stab_{F} (([0,1]\setminus Z)\cup Y)
	\setminus stab _{F} (\{ t\} )$.}
\end{exam}

\begin{exam} \label{grig}
	\emph{ (related to an argument from \cite{A})\ The action of any finitely generated weakly
	branch group on the boundary space of the corresponding infinite rooted tree is also
	hereditarily separating. To see this fix any such a group $G$ and the corresponding tree
	$T$. Now let $X:=\partial T$.} \parskip0pt

	\emph{ To see that $G$ hereditarily separates $X$, suppose that $Z$ is an open
	subset of	$X$, $\{ y_{1}, y_{2}, \ldots y_{t}\} :=Y\subset Z$ and $x\in Z\setminus Y$.
	Wlog assume that $t$ is the vertex of the infinite ray $x$ such that if any $y\in X$
	contains $t$ then $y\in Z$. Let $k$ be the level of $t$. Now choose a level $k'\geq k$
	such that the vertices in the infinite rays $y_{1},\ldots , y_{t}, x$ at the $k'$-th level are
	all distinct. Let $t_{0}$ be the vertex of $x$ at level $k'$. Let $S$ be the stabilizer of
	$t_{0}$ in $G$ and let $R$ be the rigid vertex stabilizer of $t_{0}$ in $G$. Then $S$
	acts spherically transitively on the infinite subtree $T_{t_{0}}$ rooted at $t_{0}$.
	Indeed, if $a$ and $b$ are both descendants of $t_{0}$ at the same level, then there
	is $g\in G$ such that $g(a)=b$, and clearly $g$ must stabilize $t_{0}$.}\parskip0pt

	\emph{Now suppose that there is some infinte ray of the form $t_{0}t_{1}t_{2}\ldots$
	such that for any $i\in\mathbb{N}$ and any $r\in R$, $r(t_{i})=t_{i}$. Fix any $t\in T_{t
	_{0}}$. Since $S$ acts spherically transitively on $T_{t_{0}}$, there is some $s_{0}\in S$
	such that $s_{0}(t)=t_{i}$ for some $i\in\mathbb{N}$. Clearly for any $r\in R$ we have
	$s_{0}^{-1}rs_{0}(t)=t$. Since $R$ is normal in $S$ we see that $R=R^{s_{0}^{-1}}$
	stabilizes every $t\in T_{t_{0}}$, a contradiction. Thus we assume that $R$ cannot
	stabilize any infinite ray going trough $t_{0}$. In particular, there exists $r\in R$ such
	that $r(x)\neq x$. On	the other hand, $r$ stabilizes every ray not going through $t_{0}$.
	It follows that $r\in stab _{G}((X\setminus Z)\cup Y)$.}
\end{exam}

 Now let $G$ be a permutation group on $X$. We distinguish a specific type of words over $G$ with
respect to the action on $X$. Let $w$ be a word over $G$ on $t$ variables $y_{1},\ldots , y_{t}$. It can
be considered as an element of $\mathbb{F} _{t}\ast G$. We assume that $w$ is reduced in
$\mathbb{F} _{t}\ast G$. If $w\notin\mathbb{F} _{t}$, we usually assume that $w$ is in the form $w=u
_{k}v_{k}u_{k-1}v_{k-1}\ldots u_{1}v_{1}$, where $k\in\mathbb{N}$, $u_{i}$ depends only on variables
and $v_{i}\in G\setminus\{ 1\}$ for each $i\leq k$. In this case define:
$$O_{w}:=\bigcap _{i=0} ^{k-1} v_{0}^{-1}v_{1}^{-1}\ldots v_{i}^{-1}\Big( supp(v_{i+1})\Big) ,$$
where $v_{0}=1$. If $w\in\mathbb{F} _{t}$ then let $O_{w}:=X\setminus Fix(G)$.

\begin{defi} \label{osc}
	We say that a non-trivial word $w\in\mathbb{F} _{t}\ast G$ is \emph{oscillating} if $w\in
	\mathbb{F} _{t}$ or $w$ is in the form above and $O_{w}\neq\emptyset$.
\end{defi}

\begin{rema}
	 \emph{In the situation where $G$ acts on $X$ by homeomorphisms the set $O_{w}$ is
	open.}
\end{rema}

 Note that if all $v_{i}$ are taken from the same cyclic subgroup of $G$, then $w$ is oscillating. The
next example is different.

\begin{exam} \label{1}
	\emph{ As an example once again consider the Thompson's group $F$ with its standard
	action on $[0,1]$. Let $w_{1}=yx_{1}y^{-1}x_{2}y^{2}x_{1}^{-1}$, where $y$ is a variable.
	In the notation from the definitions we have $v_{1}=x_{1}^{-1}$, $v_{2}=x_{2}$, $v_{3}=x_{1}$
	and hence:
	$$O_{w_{1}}=x_{1}x_{2}^{-1}\bigg(\Big(\frac{1}{2},1\Big)\bigg)\cap x_{1}\bigg(\Big(\frac{3}{4},1
		\Big)\bigg)\cap\Big(\frac{1}{2},1\Big) = \Big(\frac{5}{8},1\Big) .$$
	Thus we see that $w_{1}$ is oscillating (in $(\frac{5}{8},1)$).}
\end{exam}

\setlength{\unitlength}{2mm}
\begin{picture}(45,45)
	\linethickness{0.05mm}
		\put(-4,9){\line(0,1){32}}
		\put(-4,9){\line(1,0){8}}
		\put(4,9){\line(0,1){32}}
		\put(-4,41){\line(1,0){8}}

			\put(-4,25){\line(1,0){8}}
			\put(-4,33){\line(2,-1){8}}
			\put(-4,35){\line(4,-1){8}}
			\put(-4,37){\line(1,0){8}}
			\put(-5,7){\tiny $x_{1}x_{2}^{-1}(supp(v_{3}))$}

		\put(7,9){\line(0,1){32}}
		\put(7,9){\line(1,0){8}}
		\put(15,9){\line(0,1){32}}
		\put(7,41){\line(1,0){8}}

			\put(18,25){\line(1,0){8}}
			\put(7,33){\line(2,-1){8}}
			\put(7,37){\line(2,-1){8}}
			\put(7.1,7){\tiny $x_{1}(supp(v_{2}))$}

		\put(18,9){\line(0,1){32}}
		\put(18,9){\line(1,0){8}}
		\put(26,9){\line(0,1){32}}
		\put(18,41){\line(1,0){8}}

			\put(7,25){\line(1,0){2}}
			\put(10,25){\line(1,0){2}}
			\put(13,25){\line(1,0){2}}

			\put(19.5,7){\tiny $supp(v_{1})$}

	\linethickness{0.375mm}

		\put(-4,25){\line(0,1){16}}
		\put(4,25){\line(0,1){16}}
		\put(7,33){\line(0,1){8}}
		\put(15,29){\line(0,1){12}}
		\put(18,25){\line(0,1){16}}
		\put(26,25){\line(0,1){16}}

	\linethickness{0.125mm}
		\put(4,29){\line(1,0){0.5}}
		\put(5.5,29){\line(1,0){0.5}}
		\put(7,29){\line(1,0){0.5}}
		\put(8.5,29){\line(1,0){0.5}}
		\put(10,29){\line(1,0){0.5}}
		\put(11.5,29){\line(1,0){0.5}}
		\put(13,29){\line(1,0){0.5}}
		\put(14.5,29){\line(1,0){0.5}}
		\put(16,29){\line(1,0){0.5}}
		\put(17.5,29){\line(1,0){0.5}}
		\put(19,29){\line(1,0){0.5}}
		\put(20.5,29){\line(1,0){0.5}}
		\put(22,29){\line(1,0){0.5}}
		\put(23.5,29){\line(1,0){0.5}}
		\put(25,29){\line(1,0){0.5}}

		\put(4,41){\line(1,0){0.5}}
		\put(5.5,41){\line(1,0){0.5}}
		\put(7,41){\line(1,0){0.5}}
		\put(8.5,41){\line(1,0){0.5}}
		\put(10,41){\line(1,0){0.5}}
		\put(11.5,41){\line(1,0){0.5}}
		\put(13,41){\line(1,0){0.5}}
		\put(14.5,41){\line(1,0){0.5}}
		\put(16,41){\line(1,0){0.5}}
		\put(17.5,41){\line(1,0){0.5}}
		\put(19,41){\line(1,0){0.5}}
		\put(20.5,41){\line(1,0){0.5}}
		\put(22,41){\line(1,0){0.5}}
		\put(23.5,41){\line(1,0){0.5}}
		\put(25,41){\line(1,0){0.5}}

		\put(26,9){\line(1,0){0.5}}
		\put(27,9){\line(1,0){0.5}}
		\put(29.5,9){\line(1,0){0.5}}
		\put(30.5,9){\line(1,0){0.5}}
		\put(26,25){\line(1,0){0.5}}
		\put(27,25){\line(1,0){0.5}}
		\put(29.5,25){\line(1,0){0.5}}
		\put(30.5,25){\line(1,0){0.5}}
		\put(26,29){\line(1,0){0.5}}
		\put(27,29){\line(1,0){0.5}}
		\put(29.5,29){\line(1,0){0.5}}
		\put(30.5,29){\line(1,0){0.5}}
		\put(26,33){\line(1,0){0.5}}
		\put(27,33){\line(1,0){0.5}}
		\put(29.5,33){\line(1,0){0.5}}
		\put(30.5,33){\line(1,0){0.5}}
		\put(26,37){\line(1,0){0.5}}
		\put(27,37){\line(1,0){0.5}}
		\put(29.5,37){\line(1,0){0.5}}
		\put(30.5,37){\line(1,0){0.5}}
		\put(26,41){\line(1,0){0.5}}
		\put(27,41){\line(1,0){0.5}}
		\put(29.5,41){\line(1,0){0.5}}
		\put(30.5,41){\line(1,0){0.5}}

		\put(28.1,8.7){{\footnotesize $0$}}
		\put(28.1,24.9){{\footnotesize $\frac{1}{2}$}}
		\put(28.1,28.9){{\footnotesize $\frac{5}{8}$}}
		\put(28.1,32.9){{\footnotesize $\frac{3}{4}$}}
		\put(28.1,36.9){{\footnotesize $\frac{7}{8}$}}
		\put(28.1,40.7){{\footnotesize $1$}}

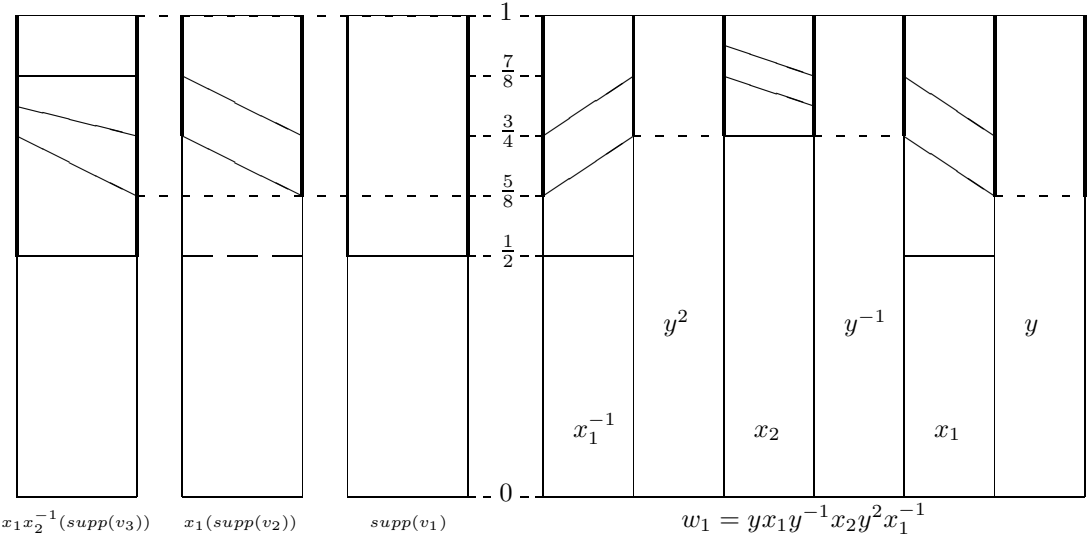
	\put(-5,4){{\footnotesize Picture 6.\ \ Supports of constant segments}}
 	\put(1,1.5){{\footnotesize in an \emph{oscillating} word $w_{1}$}}

	\linethickness{0.05mm}
		\put(31,9){\line(1,0){36}}
		\put(31,41){\line(1,0){36}}
		\multiput(31,9)(6,0){7}%
			{\line(0,1){32}}

			\put(31,25){\line(1,0){6}}
			\put(31,29){\line(3,2){6}}
			\put(31,33){\line(3,2){6}}
			\put(33,13){\footnotesize $x_{1}^{-1}$}

			\put(39,20){\footnotesize $y^{2}$}

			\put(43,33){\line(1,0){6}}
			\put(43,37){\line(3,-1){6}}
			\put(43,39){\line(3,-1){6}}
			\put(45,13){\footnotesize $x_{2}$}

			\put(51,20){\footnotesize $y^{-1}$}

			\put(55,25){\line(1,0){6}}
			\put(55,33){\line(3,-2){6}}
			\put(55,37){\line(3,-2){6}}
			\put(57,13){\footnotesize $x_{1}$}

			\put(63,20){\footnotesize $y$}

			\put(40.2,7){\footnotesize $w_{1}=yx_{1}y^{-1}x_{2}y^{2}x_{1}^{-1}$}

	\linethickness{0.375mm}

			\put(31,29){\line(0,1){12}}
			\put(37,33){\line(0,1){8}}
			\put(43,33){\line(0,1){8}}
			\put(49,33){\line(0,1){8}}
			\put(55,33){\line(0,1){8}}
			\put(61,29){\line(0,1){12}}
			\put(67,29){\line(0,1){12}}

	\linethickness{0.125mm}

		\put(37,33){\line(1,0){0.5}}
		\put(38.5,33){\line(1,0){0.5}}
		\put(40,33){\line(1,0){0.5}}
		\put(41.5,33){\line(1,0){0.5}}

		\put(49,33){\line(1,0){0.5}}
		\put(50.5,33){\line(1,0){0.5}}
		\put(52,33){\line(1,0){0.5}}
		\put(53.5,33){\line(1,0){0.5}}

		\put(61,29){\line(1,0){0.5}}
		\put(62.5,29){\line(1,0){0.5}}
		\put(64,29){\line(1,0){0.5}}
		\put(65.5,29){\line(1,0){0.5}}

	\put(33.5,4){{\footnotesize Picture 7.\ \ The images of points from $O_{w_{1}}$}}
 	\put(34.5,1.5){{\footnotesize always meet the supports of constants}}

\end{picture}

 We introduce the following notation: for any $w$ given in the form $w=u_{k}v_{k}u_{k-1}v_{k-1}\ldots u
_{1}v_{1}$ as above and for any set $A\subseteq X$ by $\mathcal{V}_{w}(A)$ we will denote the set
$\bigcup _{j=1}^{k} v_{j}\ldots v_{1}(A)$ and by $\mathcal{V}_{w}^{-1}(A)$ the set $\bigcup _{j=1} ^{k}v
_{1}^{-1}\ldots v_{j}^{-1}(A)$. \parskip0pt

 The following statement is related to Theorem 1.1 from \cite{A}.

\begin{theo} \label{ab}
	 Let $G$ be a group acting on a perfect Polish space $X$ by homeomorphisms. Let $w$
	be a word over $G$ on $t$ variables, $y_{1},\ldots , y_{t}$, which is reduced and
	non-constant (i. e. $w\notin G$) in $\mathbb{F}_{t}\ast G$. If $G$ hereditarily separates
	$X$ and $w$ has a conjugate in $\mathbb{F}_{t}\ast G$, which is oscillating, then the
	inequality $w\neq 1$ has a solution in $G$. \parskip0pt

	 Moreover, assume that $w\in\mathbb{F} _{t}$ or $w=u_{k}v_{k}\ldots u_{1}v_{1}$,
	where for each $s$, $1\leq s\leq k$, $v_{s}\in G\setminus\{ 1\}$ and $u_{s}$ contains
	only variables. Then for any open set $O'\subseteq O_{w}$, there is a tuple $\bar{g}=(g_{1},
	\ldots , g_{t})$ such that $w(\bar{g})\neq 1$, and for all $i$, $1\leq i\leq t$, the following
	conditions hold: \parskip2pt

	\ \ $\bullet$\ $supp(g_{i})\subseteq\mathcal{V} _{w}(O')$, \parskip2pt

	\ \ $\bullet$\ $g_{i}(v_{s}\ldots v_{1}(O'))=v_{s}\ldots v_{1}(O')$ for any $s$, $1\leq s
		\leq k$.
\end{theo}

\emph{Proof.}\ Using conjugation (if necessary) we assume that if $w\notin\mathbb{F} _{t}$ then $w$
is written in the form:
$$w=u_{k,l_{k}}\ldots u_{k,1}v_{k}\ldots u_{2,l_{2}}\ldots u_{2,1}v_{2}u_{1,l_{1}}\ldots  u_{1,1}v_{1},$$
where $u_{j,i_{j}}\in\{ y_{1}^{\pm 1},\ldots , y_{t}^{\pm 1}\}$, $1\leq i_{j}\leq l_{j}$, and $v_{j}\in G
\setminus\{ 1\}$, $1\leq j\leq k$. Let $L_{j}:=\sum_{i=1}^{j} l_{i}$, $1\leq j\leq k$. For any $1\leq r
\leq L_{k}$ we denote the final segment of $w$ determined by final $r$ occurances of letters from $\{
y_{1}^{\pm 1},\ldots , y_{t}^{\pm 1}\}$ by:
$$(w)_{r}=u_{d,s}\ldots u_{d,1}v_{d}\ldots u_{2,l_{2}}\ldots u_{2,1}v_{2}u_{1,l_{1}}\ldots  u_{1,1}v
	_{1},$$
where $r=L_{d-1}+s$, $1\leq s\leq l_{d}$. By $(w)_{r}(g_{1},\ldots , g_{t})$ we denote the value of
$(w)_{r}$ in $G$ via the substitution $y_{i}=g_{i}$, $1\leq i\leq t$, for a tuple of elements $\bar{g}=(g
_{1},\ldots , g_{t})\in G^{t}$. To simplify notation let also $(w)_{0}:=1$. \parskip0pt

 Let $\bar{g}\in G$ and $p\in X$. Define $p_{r,\bar{g}}:=(w)_{r}(\bar{g})(p)$ for all $r$, $1\leq r\leq L
_{k}$. We say that $\bar{g}$ is \emph{distinctive} for the word $w$ and the point $p$, if all the points
$$p=p_{0,\bar{g}}, v_{1}(p_{0,\bar{g}}),\ldots , p_{l_{1},\bar{g}}, v_{2}(p_{l_{1},\bar{g}}),\ldots , p_{n,
\bar{g}}$$
are pairwise distinct. Observe that to prove the proposition we need a weaker condition: find $p\in X$
and a tuple of elements $\bar{g}\in G$ such that $w(\bar{g})(p)\neq p$. \parskip0pt

 Fix $O'$. We will prove the theorem by induction. At $n$-th step we will show that: \parskip2pt

\ \ $\bullet$ \ There is $p\in O'$ and a tuple $\bar{g}=(g_{1},\ldots , g_{t})\in G$ such that $\bar{g}$ is
	distinctive	for $p$ and $(w)_{n}$. \parskip2pt

\ \ $\bullet$ \ In the condition above we can choose $\bar{g}$ so that for all $i$, $1\leq i\leq t$, $supp(g
	_{i})\subseteq\mathcal{V}_{w}(O')$ and $g_{i}(v_{r}\ldots v_{1}(O'))=v_{r}\ldots v_{1}(O')$ for
	any $r\leq k$. \parskip2pt

We make the following observation for further reference.\\

 \textbf{Claim.} For any $r\leq k$ and any $q\in v_{r}\ldots v_{1}(O')\setminus\bigcup _{i=1} ^{k} (v_{i}
\ldots v_{1}(\bar{O'}\setminus O'))$ we can find a neighbourhood $O\subseteq v_{r}\ldots v_{1}(O')$ of
$q$ such that the following condition holds:
$$(\dag )\ \forall s\leq k\ \bigg( O\cap v_{s}\ldots v_{1}(O')\neq\emptyset\ \Rightarrow\ O\subseteq v_{s}
	\ldots v_{1}(O')\bigg) .$$

 To see this we construct inductively the set $O$ for a given
$$q\in v_{r}\ldots v_{1}(O')\setminus\Big(\bigcup _{i=1} ^{k}\Big( v_{i}\ldots v_{1}(\bar{O'}\setminus O')\Big)
	\Big) .$$
 Let $O_{0}:=v_{r}\ldots v_{1}(O')$. At the $n$-th, $1\leq n\leq k$, step of the construction we consider $O
_{n-1}\cap v_{n}\ldots v_{1}(O')$, where $O_{n-1}$ is defined at Step $n-1$. If $q\in O_{n-1}\cap v_{n}
\ldots v_{1}(O')$ then we define $O_{n}:=O_{n-1}\cap v_{n}\ldots v_{1}(O')$. Now suppose that $q\notin O
_{n-1}\cap v_{n}\ldots v_{1}(O')$. Since $q\notin\bigcup _{i=1} ^{k} (v_{i}\ldots v_{1}(\bar{O'}\setminus O'))$,
we find some neighbourhood $O_{n-1}'\subseteq O_{n-1}\setminus v_{n}\ldots v_{1}(O')$ of $q$ and
define $O_{n}:=O_{n-1}'$. After $k$ steps we define $O:=O_{k}$. It is clear that the set $O$ is non-empty
and open. If $w\in\mathbb{F}_{t}$ then we define $O:=O'$, what finishes the proof of the claim.\\

 Fix $p\in O'$. For $n=1$, $(w)_{1}$ is of the form $w=y_{i}^{\pm 1}v_{1}$ for some $1\leq i\leq t$
(when $w\in\mathbb{F} _{t}$ we replace $v_{1}$ by $1$ and follow the argument below). According to
the assumptions, $p\neq v_{1}(p)$ for all $p\in O'$ (except the case, when $w\in\mathbb{F} _{t}$ and
this condition is redundant). Since the set $\bigcup _{i=1} ^{k} (v_{i}\ldots v_{1}(\bar{O'}\setminus O'))$
is nowhere dense in $X$ and the action of $G$ is continuous, we can slightly modify our choice of $p$
(if necessary) so that the inequality $p\neq v_{1}(p)$ is still satisfied, but $v_{1}(p)\notin\bigcup _{i=1}
^{k}(v_{i}\ldots v_{1}(\bar{O'}\setminus O'))$. Now it follows from the claim above that we can find a
neighbourhood $O\subseteq v_{1}(O')$ of $v_{1}(p)$ such that $O$ satisfies $(\dag )$. Wlog suppose
$w=y_{i}v_{1}$. Since $G$ hereditarily separates $X$, the $stab_{G}(X\setminus O)$-orbit of $v_{1}
(p)$ is infinite. Thus we can choose $f\in stab_{G}(X\setminus O)$ such that $f(v_{1}(p))\notin\{ p, v
_{1}(p)\}$. Defining $g_{i}:=f$ and choosing any set of $t-1$ elements from $stab_{G}(X\setminus
O)$ we obtain a tuple $\bar{g}$ distinctive for $p$ and $(w)_{1}$. \parskip0pt

 Since for all $s$, $g_{s}\in stab_{G}(X\setminus O)$, the condition $supp(g_{s})\subseteq\mathcal{V}
_{w}(O')$ is satisfied for all $s\leq t$ . \parskip0pt

 Now fix $i\leq t$ and $r\leq k$. If $O\cap v_{r}\ldots v_{1}(O')=\emptyset$, then obviously $g_{i}(v_{r}
\ldots v_{1}(O'))=v_{r}\ldots v_{1}(O')$. On the other hand, if $O\cap v_{r}\ldots v_{1}(O')\neq
\emptyset$, then it follows from the construction of $O$ that $O\subseteq v_{r}\ldots v_{1}(O')$. Thus
$supp(g_{i})$ is also a subset of $v_{r}\ldots v_{1}(O')$ and the condition $g_{i}(v_{r}\ldots v_{1}(O'))=
v_{r}\ldots v_{1}(O')$ is satisfied. \parskip0pt

 Assume that for
$$(w)_{n-1}=u_{d,s}\ldots u_{d,1}v_{d}\ldots u_{2,1}v_{2}u_{1,l_{1}}\ldots u_{1,1}v_{1} ,$$
where $n-1=L_{d-1}+s$, we can find $p\in O'$ and a tuple $\bar{g}\in G$ such that the induction
hypothesis is satisfied. According to the form of $(w)_{n}$ we consider two cases.\\

\textbf{Case 1.} $(w)_{n}=u_{d,s+1}u_{d,s}\ldots u_{d,1}v_{d}\ldots u_{2,l_{2}}v_{2}u_{1,l_{1}}\ldots
	u_{1,1}v_{1}$, where $n-1=L_{d-1}+s$, $s\geq 1$.\parskip0pt

 If
$$p_{n,\bar{g}}\notin\Big\{ p_{i,\bar{g}}\ \Big|\ 0\leq i\leq n-1\Big\}\cup\Big\{ v_{1}(p_{0,\bar{g}}),\ldots
	, v_{d}(p_{L_{d-1},\bar{g}})\Big\} ,$$
then we have found a right tuple $\bar{g}$. So we assume that $p_{n,\bar{g}}=p_{m,\bar{g}}$ for some
$0\leq m<n$ or $p_{n,\bar{g}}=v_{m+1}(p_{L_{m},\bar{g}})$ for some $0\leq m<d-1$. \parskip0pt

 Let $y_{j}^{\pm 1}$ be the first letter of $(w)_{n}$. Replacing $y_{j}$ by $y_{j}^{-1}$ and $g_{j}$ by
$g_{j}^{-1}$ if necessary, we may assume that $u_{n}=y_{j}$. Put
$$Y:=\Big\{ p_{i,\bar{g}}\ \Big|\ 0\leq i\leq n-2\Big\}\cup\Big\{ v_{1}(p_{0,\bar{g}}),\ldots , v_{d}(p_{L_{d-1},
	\bar{g}})\Big\}$$
for $p$ and $\bar{g}$ chosen at the $n-1$-th step of induction. Since $p\in O'$ and $g_{s}^{\pm 1}(v
_{r}\ldots v_{1}(O'))=v_{r}\ldots v_{1}(O')$ for all $r\leq k$ and all $ s\leq t$ (by the induction
hypothesis), $p_{n-1,\bar{g}}\in v_{d}\ldots v_{1}(O')$. Since the set $\bigcup _{i=1} ^{k} (v_{i}\ldots v
_{1}(\bar{O'}\setminus O'))$ is nowhere dense in $X$ and the action of $G$ is continuous, as in the
case $n=1$, we can slightly modify our choice of $p$ (if necessary) so that all the previous
assumptions and inequalities remain true and additionally $p_{n-1,\bar{g}}\notin\bigcup _{i=1} ^{k}(v
_{i}\ldots v_{1}(\bar{O'}\setminus O'))$. As above we choose a neighbourhood $O\subseteq v_{d}\ldots
v_{1}(O')$ of the point $p_{n-1,\bar{g}}\in v_{d}\ldots v_{1}(O')$ satisfying the condition formulated in
Claim. Since the action of $G$ is hereditarily separating, the $stab_{G}((X\setminus O)\cup Y)$-orbit
of $p_{n-1,\bar{g}}$ is infinite. We now introduce another set of points:
$$Z:=\Big\{ g_{j}^{-1}(p_{i,\bar{g}})\ \Big| \ 0\leq i\leq n-1\Big\}\cup\Big\{ g_{j}^{-1}(v_{i+1}(p_{L_{i},\bar{g}})\
	\Big|\ 0\leq i\leq d-1\Big\} .$$
Since $Z$ is finite, there exists $f\in stab_{G}((X\setminus O)\cup Y)$ taking $p_{n-1,\bar{g}}$
outside $Z$. Replacing $g_{j}$ by $g_{j}f$ we obtain a corrected tuple $\bar{g}$. Since the
element $f$ have been chosen from the stabilizer $stab_{G}(Y)$, the points
$$p_{0,\bar{g}},v_{1}(p_{0,\bar{g}}), p_{1,\bar{g}},\ldots , v_{d}(p_{L_{d-1},\bar{g}}), p_{L_{d-1},
	\bar{g}},\ldots , p_{n-1,\bar{g}}$$
are the same as at the $(n-1)$-th step of induction. On the other hand $p_{n,\bar{g}}$ is distinct from
all elements
$$p_{0,\bar{g}},v_{1}(p_{0,\bar{g}}), p_{1,\bar{g}},\ldots , v_{d}(p_{L_{d-1},\bar{g}}), p_{L_{d-1},
	\bar{g}},\ldots, p_{n-1,\bar{g}} .$$

 Since $f\in stab_{G}((X\setminus O)\cup Y)$, $supp(f)\subseteq\mathcal{V}_{w}(O')$. This
together with the induction hypothesis implies that $supp(g_{j})\subseteq\mathcal{V}_{w}(O')$ for
the corrected $g_{j}$. \parskip0pt

 Let $r\leq k$. If $O\cap v_{r}\ldots v_{1}(O')=\emptyset$, then by the choice of $f$ and induction
we have $g_{j}(v_{r}\ldots v_{1}(O'))=v_{r}\ldots v_{1}(O')$. If $O\cap v_{r}\ldots v_{1}(O')\neq
\emptyset$, then $O\subseteq v_{r}\ldots v_{1}(O')$. Thus $f$ and the original $g_{j}$ (defined at
Step $n-1$) stabilize $v_{r}\ldots v_{1}(O')$ setwise (by the choice of $f$ and induction). Thus $g
_{j}f(v_{r}\ldots v_{1}(O'))=v_{r}\ldots v_{1}(O')$. On the other hand all elements $g_{s}$ for $s
\neq j$ have not been changed and thus automatically satisfy required conditions. This finishes
the proof of Case 1. \\

\textbf{Case 2.} $(w)_{n}=u_{d+1,1}v_{d+1}u_{d,s}\ldots u_{d,1}v_{d}\ldots u_{2,l_{2}}v_{2}u_{1,
	l_{1}}\ldots u_{1,1}v_{1}$, where $n-1=L_{d}$. \parskip0pt

 If
$$p_{n,\bar{g}}\notin\Big\{ p_{i,\bar{g}}\ \Big|\ 0\leq i\leq n-1\Big\}\cup\Big\{ v_{1}(p_{0,\bar{g}}),
	\ldots , v_{d+1}(p_{L_{d},\bar{g}})\Big\}$$
and
$$v_{d+1}(p_{n-1,\bar{g}})\notin\Big\{ p_{i,\bar{g}}\ \Big|\ 0\leq i\leq n-1\Big\}\cup\Big\{ v_{1}(p_{0,
	\bar{g}}),\ldots , v_{d}(p_{L_{d-1},\bar{g}})\Big\} ,$$
then we have found a right tuple $\bar{g}$.\parskip0pt

 Assume the contrary. Once again suppose that $y_{j}^{\pm 1}$ is the first letter of $(w)_{n}$.
Replacing $y_{j}$ by $y_{j}^{-1}$ and $g_{j}$ by $g_{j}^{-1}$ if necessary, we can assume that
$u_{n}=y_{j}$. Let $u_{n-1}=y_{j'}^{\pm 1}$. Then let
$$Y':=\Big\{ p_{i,\bar{g}}\ \Big|\ 0\leq i\leq n-1\Big\}\cup\Big\{ v_{1}(p_{0,\bar{g}}),\ldots , v_{d}(p
	_{L_{d-1},\bar{g}})\Big\}$$
 Assume that $v_{d+1}(p_{n-1,\bar{g}})\in Y'$. By the induction hypothesis, $p_{n-1,\bar{g}}\in v
_{d}\ldots v_{1}(O')$. After a slight modyfication of our choice of $p$ (if necessary), once again
wlog we assume that $p_{n-1,\bar{g}}\notin\bigcup _{i=1} ^{k}(v_{i}\ldots v_{1}(\bar{O'}\setminus
O'))$. Thus we can find some neighbourhood $O\subseteq v_{d}\ldots v_{1}(O')$ of the point $p
_{n-1,\bar{g}}\in v_{d}\ldots v_{1}(O')$ satisfying $(\dag )$. Since the action of $G$ is hereditarily
separating, the $stab_{G}((X\setminus O)\cup (Y'\setminus\{ p_{n-1,\bar{g}}\} ))$-orbit of $p_{n-1,
\bar{g}}$ is infinite. We replace $g_{j'}$ by some $f'g_{j'}$ (or $g_{j'}f'$ in the case $u_{n-1}=y_{j'}
^{-1}$), where
$$f'\in stab_{G}\Big( (X\setminus O)\cup (Y'\setminus\{ p_{n-1,\bar{g}}\})\Big)$$
takes $p_{n-1,\bar{g}}$ outside the finite set $Y'\cup v_{d+1}^{-1}(Y')$. Since
$$O\subseteq v_{d}\ldots v_{1}(O')\subseteq supp(v_{d+1}) ,$$
the corrected $p_{n-1,\bar{g}}$ is not fixed by $v_{d+1}$. Thus the corrected $v_{d+1}(p_{n-1,
\bar{g}})$ surely omits the corrected $Y'$. \parskip0pt

 Moreover, since $f'\in stab_{G}((X\setminus O)\cup (Y'\setminus\{ p_{n-1,\bar{g}}\}))$, its
support is contained in $O$, hence we still have $supp(g_{j'})\subseteq\mathcal{V}_{w}(O')$.
Similarly as in Case 1 we see that the condition $supp(f')\subseteq O$ implies that $g_{j'}^{\pm 1}(v
_{r}\ldots v_{1}(O'))=(v_{r}\ldots v_{1}(O'))$ for any $r\leq k$. \parskip0pt

 So we only have to consider the case when $v_{d+1}(p_{n-1,\bar{g}})\notin Y'$, but either $p_{n,
\bar{g}}=p_{j,\bar{g}}$ for some $0\leq j<n$ or $p_{n,\bar{g}}=v_{j+1}(p_{L_{j}\bar{g}})$ for some
$0\leq j\leq d$. Let:
$$Y:=\Big\{ p_{i,\bar{g}}\ \Big|\ 0\leq i\leq n-1\Big\}\cup\Big\{ v_{1}(p_{0,\bar{g}}),\ldots , v_{d}(p_{L
	_{d-1},\bar{g}})\Big\}$$
and
$$Z:=\Big\{g_{j}^{-1}(p_{i,\bar{g}})\ \Big| \ 0\leq i\leq n-1\Big\}\cup\Big\{ g_{j}^{-1}(v_{i+1}(p_{L_{i},
	\bar{g}}))\ \Big|\ 0\leq	 i\leq d\Big\} .$$

 Now observe that $v_{d+1}(p_{n-1,\bar{g}})\in v_{d+1}\ldots v_{1}(O')\setminus Y$ and choose
the neighbourhood $O\subseteq v_{d+1}\ldots v_{1}(O')$ of the point $v_{d+1}(p_{n-1,\bar{g}})$
satisfying $(\dag )$. Then there exists $f\in stab_{G}((X\setminus O)\cup Y)$ taking $v_{d+1}(p
_{n-1,\bar{g}})$ outside $Z$. Replacing $g_{j}$ by $g_{j}f$ we finish the proof exactly as in Case 1.\\

\ \ \ \ \ \ \ \ \ \ \ \ \ \ \ \ \ \ \ \ \ \ \ \ \ \ \ \ \ \ \ \ \ \ \ \ \ \ \ \ \ \ \ \ \ \ \ \ \ \ \ \ \ \ \ \ \ \ \ \ \ \ \ \ \ \ \ \ \ \ \ \ \ \ \ \ \ $\square$\\

 In fact we also have a topology-free version of Theorem \ref{ab}, which generalizes Theorem 1.1
from \cite{A}. For any $A\subseteq X$ by $A^{0}$ we will denote $X\setminus A$ and by $A^{1}$ the set
$A$. Fix any $w$ such that either $w\in\mathbb{F} _{t}$ or $w=u_{k}v_{k}\ldots u_{1}v_{1}$, where for
each $s\leq k$, $v_{s}\in G\setminus\{ 1\}$ and $u_{s}$ contains only variables. Then for
any $\bar{\varepsilon} =(\varepsilon _{1},\ldots ,\varepsilon _{k})\in\{ 0, 1\} ^{k}$, we denote by $O _{w}
^{\bar{\varepsilon}}$ the set $\bigcap _{s=1} ^{k}(v_{s}\ldots v_{1}(O_{w}))^{\varepsilon _{s}}$.

\begin{theo} \label{gab}
	 Suppose $G$ acts by permutations on some set $X$. Let $w$ be a word over $G$ on
	$t$ variables, $y_{1},\ldots , y_{t}$, which is reduced and non-constant (i. e. $w\notin G$)
	in $\mathbb{F}_{t}\ast G$. Assume also that $w\in\mathbb{F} _{t}$ or $w=u_{k}v_{k}\ldots
	u_{1}v_{1}$, where for each $s\leq k$, $v_{s}\in G\setminus\{ 1\}$ and $u_{s}$ contains
	only variables. If $O_{w}\neq\emptyset$ and for any $\bar{\varepsilon}\in\{ 0, 1\} ^{k}$ such
	that $O _{w}^{\bar{\varepsilon}}\neq\emptyset$, the condition
	$$(\diamondsuit )\ stab _{G}(X\setminus O _{w}^{\bar{\varepsilon}})\ \emph{separates}\
		O _{w}^{\bar{\varepsilon}}$$
	is satisfied, then the inequality $w\neq 1$ has a solution in $G$.
\end{theo}

\emph{Proof.}\ If $w\in\mathbb{F}_{t}$, then $O_{w}=X$ and we simpy apply Theorem 1.1 from
\cite{A}. \parskip0pt

 If $w\notin\mathbb{F}_{t}$, then keeping the notation used previously, we follow the proof of
Theorem \ref{ab}. Note that $O_{w}^{\bar{\varepsilon}}$ is infinite for $O_{w}^{\bar{\varepsilon}}\neq
\emptyset$.\parskip0pt

 For $O':=O_{w}$ we refolmulate the claim from the of Theorem \ref{ab} in the following form:\\

 \textbf{Claim $\#$.} For any $r\leq k$ and any $q\in v_{r}\ldots v_{1}(O_{w})$ there is a unique
tuple $\bar{\varepsilon}\in\{ 0, 1\} ^{k}$ such that $q\in O _{w}^{\bar{\varepsilon}}\subseteq v_{r}\ldots v
_{1}(O_{w})$. The corresponding $O_{w}^{\bar{\varepsilon}}$ satisfies:
$$(\ddag )\ \forall s\leq k\ O _{w}^{\bar{\varepsilon}}\cap v_{s}\ldots v_{1}(O_{w})\neq\emptyset\
	\Rightarrow\ O _{w}^{\bar{\varepsilon}}\subseteq v_{s}\ldots v_{1}(O_{w}).$$
We prove Claim $\#$ as follows. First observe that $\{ O_{w}^{\bar{\varepsilon}}\ |\ \bar{\varepsilon}
\in\{ 0, 1\} ^{k}\}$ is a partition of $X$ and hence for any $q\in X$ there is a unique tuple
$\bar{\varepsilon}\in\{ 0, 1\} ^{k}$ such that $q\in O _{w}^{\bar{\varepsilon}}\subseteq v_{r}\ldots v_{1}
(O_{w})$. Now fix this tuple $\bar{\varepsilon}$ and suppose that $O _{w}^{\bar{\varepsilon}}\cap v
_{s}\ldots v_{1}(O_{w})\neq\emptyset$. It follows that $\varepsilon _{s} =1$ in $\bar{\varepsilon}$.
Thus $O _{w}^{\bar{\varepsilon}}\subseteq v_{s}\ldots v_{1}(O_{w})$.\parskip0pt

 We now apply the proof of Theorem $\ref{ab}$. Let $O'=O_{w}$. At the $n$-th step of induction
we show that: \parskip2pt

\ \ $\bullet$ \ There is $p\in O'$ and a tuple $\bar{g}=(g_{1},\ldots , g_{t})\in G$ such that $\bar{g}$ is
	distinctive	for $p$ and $(w)_{n}$. \parskip2pt

\ \ $\bullet$ \ In the condition above we can choose $\bar{g}$ so that for all $i$, $1\leq i\leq t$,
	$supp(g_{i})\subseteq\mathcal{V}_{w}(O')$ and $g_{i}(v_{r}\ldots v_{1}(O'))=v_{r}\ldots v
	_{1}(O')$ for any $r\leq k$. \parskip2pt

 Fix any $p\in O_{w}$ and any $n\in\mathbb{N}$. By Claim $\#$ we obtain $\bar{\varepsilon}\in\{ 0, 1\}
^{k}$ such that $O _{w}^{\bar{\varepsilon}}\subseteq v_{1}(O_{w})$ contains $p_{n-1,\bar{g}}$ and
satisfies $(\ddag )$. Hence we may now apply $(\diamondsuit )$ assumption and Remark \ref{rema}
to see that for any finite set $Y$, the $stab _{G}((X\setminus O_{w}^{\bar{\varepsilon}})\cup Y)$-orbit
of $p_{n-1,\bar{g}}$ is infinite.\parskip0pt

 Thus we replace each occurance of the the neighbourhood $O\subseteq v_{d}\ldots v_{1}(O')$ of the
point $p_{n-1,\bar{g}}$ constructed in the proof of Theorem \ref{ab} by the set $O_{w}
^{\bar{\varepsilon}}$ obtained from Claim $\#$. Consequently, each time we use the condition
$(\ddag )$ instead of condition $(\dag )$ (from the proof of Theorem \ref{ab}) to find a desired tuple
$\bar{g}$. After this modifications the proof of Theorem \ref{ab} works in the case of Theorem
\ref{gab}.\\

\ \ \ \ \ \ \ \ \ \ \ \ \ \ \ \ \ \ \ \ \ \ \ \ \ \ \ \ \ \ \ \ \ \ \ \ \ \ \ \ \ \ \ \ \ \ \ \ \ \ \ \ \ \ \ \ \ \ \ \ \ \ \ \ \ \ \ \ \ \ \ \ \ \ \ \ \ $\square$\\

\subsection{Solving a system of inequalities in the case of almost oscillating words and words with
	non-trivial product of constants}

 Now we would like to apply Theorem \ref{ab} to the broader class of words and larger number of
inequalities. We will introduce the notion of an \emph{almost oscillating} word and the corresponding
set $O_{w}\subseteq X$. Intuitively it describes words, which are oscillating after reduction of certain
subwords. \parskip0pt

Suppose $G$ acts on some perfect Polish space $X$ by homeomorphisms and let $w$ be a
word over $G$ on $t$ variables such that $w$ is reduced in $\mathbb{F} _{t}\ast G$. If $w\notin
\mathbb{F} _{t}$ then we assume that $w$ is in the form
$$w_{X}=u_{X,k_{X}}v_{X,k_{X}}u_{X,k_{X}-1}v_{X,k_{X}-1}\ldots u_{X,1}v_{X,1} ,$$
where $k_{X}\in\mathbb{N}$, $u_{X,i}$ depends only on variables and $v_{X,i}\in G\setminus \{ 1\}$,
$i\leq k_{X}$. Suppose that $w_{X}$ is not oscillating. For an open $A\subseteq X$ let $A^{0}:=
\emph{int}(X\setminus A)$ and $A^{1}:=A$. For any sequence $\bar{\varepsilon} =(\varepsilon_{1},
\ldots , \varepsilon _{k_{X}})\in\{ 0,1\} ^{k}$ we define the set
$$X_{\bar{\varepsilon}}:=\bigcap _{i=1}^{k_{X}} v_{X,1}^{-1}\ldots v_{X,i-1}^{-1}\Big( supp(v_{X,i})
	^{\varepsilon _{i}}\Big) ,$$
$v_{X,0}:=1$. Now we consider the following family of sets:
$$\mathcal{P} ^{1}:=\Big\{ X_{\bar{\varepsilon}}\ \Big|\ \bar{\varepsilon}\in\{ 0, 1\} ^{k}\Big\}\setminus
	\Big\{\emptyset\Big\} .$$
Since $w$ is not oscillating, $X_{(1,\ldots , 1)}=\emptyset$. For any $X_{\bar{\varepsilon}}\in
\mathcal{P}^{1}$ we define a word $w'_{X_{\bar{\varepsilon}}}$ in the following way:
$$w'_{X_{\bar{\varepsilon}}}:=u_{X,k_{X}}v^{\varepsilon _{k_{X}}}_{X,k_{X}}u_{X,k_{X}-1}v
	^{\varepsilon _{k_{X}-1}}_{X,k_{X}-1}\ldots u_{X,1}v^{\varepsilon _{1}}_{X,1},$$
where for any $i\leq k_{X}$, $v^{0}_{X,i}=1$ and $v^{1}_{X,i}=v_{X,i}$. \parskip0pt

 We reduce each $w'_{X_{\bar{\varepsilon}}}$, $X_{\bar{\varepsilon}}\in\mathcal{P} ^{1}$, in
$\mathbb{F}_{t}\ast G$. Suppose that such a reduced word $w'$ is of the form $\bar{w}v'u'$,
where $u'$ contains only variables and $v'\in G\setminus\{ 1\}$. Then we conjugate $w'$ by
$(u')^{-1}$ and denote the obtained word by $w_{X_{\bar{\varepsilon}}}$. If it is not the case,
we simply take $w_{X_{\bar{\varepsilon}}}:=w'$. By $\mathcal{W} ^{1}$ we denote the set of
all non-trivial words $w_{X_{\bar{\varepsilon}}}$, $X_{\bar{\varepsilon}}\in\mathcal{P} ^{1}$. If
there is some oscillating $w_{X_{\bar{\varepsilon}}}\in\mathcal{W} ^{1}$, let
$$\mathcal{P}^{os}:=\Big\{ X_{\bar{\varepsilon}}\in\mathcal{P} ^{1}\ \Big|\ w_{X_{\varepsilon}}\
	\emph{is oscillating}\Big\} .$$

 If there is no oscillating $w_{X_{\bar{\varepsilon}}}\in\mathcal{W} ^{1}$ then for any $X_{\bar{
\varepsilon}}\in\mathcal{P} ^{1}$ we repeat the process described above replacing $X$ by $X
_{\bar{\varepsilon}}$ and $w_{X}$ by $w_{X_{\bar{\varepsilon}}}$. For any $X_{\bar{\varepsilon}}
\in\mathcal{P} ^{1}$ we define the family $\mathcal{P} ^{2} _{X_{\bar{\varepsilon}}}$ and the
corresponding set of words $\mathcal{W} ^{2} _{X_{\bar{\varepsilon}}}$ exactly as $\mathcal{P} ^{1}$
and $\mathcal{W} ^{1}$ were defined above. Now let
$$\mathcal{P} ^{2}:=\bigcup\Big\{ \mathcal{P} ^{2} _{X_{\bar{\varepsilon}}}\ \Big|\ X_{\bar{\varepsilon}}\in
	\mathcal{P} ^{1}\Big\} .$$
Let $\mathcal{W} ^{2}$ be the set of all words $w_{V}$ for $V\in\mathcal{P}^{2}$ defined as $w
_{X_{\bar{\varepsilon}}}$ above. If there is an oscillating word in $\mathcal{W}^{2}$ then we define
$$\mathcal{P} ^{os}:=\Big\{ V\in\mathcal{P} ^{2}\ \Big|\ w_{V}\ \emph{is oscillating}\Big\}$$
and finish the construction. If there is no oscillating word in $\mathcal{W} ^{2}$, then we continue
this procedure. If for some $n\in\mathbb{N}$, $\mathcal{W}^{n}$ contains oscillating words or
$\mathcal{W} ^{n}=\emptyset$, then the procedure terminates. \parskip0pt

\begin{lemm}
	The procedure described above terminates after finitely many steps.
\end{lemm}

\emph{Proof.} \ Suppose that for some $n>1$ we are given a word $w_{V}\in\mathcal{W}^{n-1}$,
which is not oscillating. Notice that if $w\in\mathbb{F} _{t}<\mathbb{F} _{t}\ast G$, then $w$ is
oscillating. Thus $w_{V}$ contains some constants $v_{V,i}$, $1\leq i\leq k_{V}$. \parskip0pt

 Let $\bar{\varepsilon}\in\{ 0, 1\}^{k_{V}}$, $V_{\bar{\varepsilon}}\in\mathcal{P}^{n} _{V}$ and $w
_{V_{\bar{\varepsilon}}}\in\mathcal{W}^{n}$ be obtained from $w_{V}$ as in the construction.
Since $w_{V}$ was not oscillating,
$$U_{1,\ldots ,1}=\bigcap _{i=1}^{k_{V}} v_{V,1}^{-1}\ldots v_{V,i-1}^{-1}\Big( supp(v_{V,i})^{1}
	\Big) =\emptyset.$$
Thus for some $i$, $1\leq i\leq k_{V}$, $v_{V, i}$ becomes $1$ in the word $w_{V
_{\bar{\varepsilon}}}$. Therefore the length of $w_{V_{\bar{\varepsilon}}}$ is strictly smaller than
the length of $w_{V}$. We see that either after finitely many steps we find some $\mathcal{W}
^{n}$ containing an oscillating word or for some $\mathcal{P} ^{n}$ the words $w_{V}$ are equal
to $1$ for all $V\in\mathcal{P} ^{n}$.\\

\ \ \ \ \ \ \ \ \ \ \ \ \ \ \ \ \ \ \ \ \ \ \ \ \ \ \ \ \ \ \ \ \ \ \ \ \ \ \ \ \ \ \ \ \ \ \ \ \ \ \ \ \ \ \ \ \ \ \ \ \ \ \ \ \ \ \ \ \ \ \ \ \ \ \ \ \ $\square$\\

 We say that the initial word $w_{X}$ is \emph{rigid} if the procedure terminates and for some
$n\in\mathbb{N}$ and all $V\in\mathcal{P} ^{n}$, $w_{V}=1$. If it is not the case, then the
procedure teminates producing the set $\mathcal{P} ^{os}$. Then we say that the initial word
$w_{X}$ is \emph{almost oscillating}. \parskip0pt

 Now we rewrite each $w_{V}\notin\mathbb{F} _{t}<\mathbb{F} _{t}\ast G$, $V\in\mathcal{P}
^{os}$, in the form $w_{V}=u_{V,k_{V}}v_{V,k_{V}}\ldots u_{V,1}v_{V,1}$, where $k_{V}\in
\mathbb{N}$, for each $i\leq k_{V}u_{V,i}$ depends only on variables and $v_{V,i}\in G
\setminus\{ 1\}$. Define
$$O_{w}:=\bigcup _{V\in\mathcal{P} ^{os}} \Big( V\cap\bigcap _{i=0} ^{k_{V}-1} v_{V,0}^{-1}v
	_{V,1}^{-1}\ldots v_{V,i}^{-1}\Big( supp(v_{V,i+1})\Big)\Big) ,$$
where for all $V\in\mathcal{P} ^{os}$, $v_{V,0}=1$. In particular, for $w_{V}\in \mathbb{F}_{t}<
\mathbb{F} _{t}\ast G$, $k_{V}=0$ and the contribution of $w_{V}$ to the above sum equals
$V$.\parskip0pt

 Now let $G$ act on a perfect Polish space $X$ by homeomorphisms. Let $w$ be a word
over $G$ on $t$ variables, $y_{1},\ldots , y_{t}$, which is reduced and non-constant (i.e. $w
\notin G$) in $\mathbb{F}_{t}\ast G$. Assume that $w=u_{k}v_{k}\ldots u_{1}v_{1}$, where $u
_{s}$ contains only variables and $v_{s}\in G\setminus\{ 1\}$ for any $s\leq k$. The following
lemma exhibits a relation between the existence of solutions of the inequality $w\neq 1$ and
$w_{U}\neq 1$, where $w_{U}$ is derived from $w$ exactly as in the procedure described
above.

\begin{lemm}\label{ii}
	Suppose that for some $n\geq 1$, $U\in\mathcal{P}^{n}$ and $w_{U}\in\mathcal{W}
	^{n}$, we have $w_{U}(\bar{g})(p)\neq p$ for some point $p\in U$ and some tuple
	$\bar{g}=(g_{1},\ldots , g_{t})\in G^{t}$. If for any $i$, $1\leq i\leq t$, and any $j$,
	$1\leq j\leq k$, $g_{i}$ stabilizes $v_{j}\ldots v_{1}(U)$ setwise, then $w(\bar{g})
	\neq id$.
\end{lemm}

\emph{Proof.}\ Fix some $\bar{g}\in G$ satisfying the conditions of the lemma. Consider a
word $w_{V}\in\mathcal{W}^{r}$, where $r\leq n$ and $V\in\mathcal{P}^{r}_{V'}$ (i.e. $U
\subseteq V\subseteq V'$), which was obtained from some non-oscillating word $w_{V'}\in
\mathcal{W}^{r-1}$ ($\mathcal{W}^{0}:=\{ w\}$). To simplify notation we assume that $w_{V}$
was obtained from $w_{V'}$ by the appropriate reduction, but without the use of conjugation.
Such an assumption can be made wlog, because we have assumed that for any $j\leq k$,
$\bar{g}$ stabilizes $v_{j}\ldots v_{1}(U)$ setwise (see the argument below). We will show
that, if for some $p\in U$, $w_{V}(\bar{g})(p)\neq p$, then we have $w_{V'}(\bar{g})(p)\neq p$.
This proves the lemma by induction starting with the case $r=n$ and $w_{U}\in\mathcal{W}
^{n}$, such that $w_{U}(\bar{g})(p)\neq p$ for some point $p\in U$.\parskip0pt

 Let 
$$w_{V'}:=u_{V',k_{V'}}v_{V',k_{V'}}u_{V',k_{V'}-1}v_{V',k_{V'}-1}\ldots u_{V',1}v_{V',1},$$
where $k_{V'}\in\mathbb{N}\setminus\{ 0\}$, $u_{V',i}$ depends only on variables and $v_{V'
,i}\in G\setminus\{1\}$ for any $i\leq k_{V'}$. Since $w_{V'}$ is non-oscillating and we do not
use conjugation to get $w_{V}$, it follows from the construction above that we may assume
that $w_{V}$ is obtained by reductions from the word
$$w'_{V'_{\bar{\varepsilon}}}=u_{V',k_{V'}}v^{\varepsilon _{k_{V'}}}_{V',k_{V'}}u_{V',k_{V'}-1}
	v^{\varepsilon _{k_{V'}-1}}_{V',k_{V'}-1}\ldots u_{V',1}v^{\varepsilon _{1}}_{V',1},$$
for some $\bar{\varepsilon}\in\{ 0, 1\} ^{k_{V}}$, where for any $i\leq k_{V'}$, $v_{V',i}^{0}:=1$
and $v_{V',i}^{1}:=v_{V',i}$.\parskip0pt

 Now consider the word $w_{V'}(\bar{g})$ and the point $p\in U$. To simplify notation for any
$j\leq k_{V'}$ denote by $p'_{j}$ the point $u_{V',j-1}(\bar{g})v_{V',j-1}\ldots u_{V',1}(\bar{g})
v_{V',1}(p)$ and by $p_{j}$ the point $u_{V',j-1}(\bar{g})v^{\varepsilon _{j-1}}_{V',j-1}\ldots u
_{V',1}(\bar{g})v^{\varepsilon_{1}}_{V',1}(p)$ ($p'_{1}=p_{1}:=p$). In order to show that $w
_{V'}(\bar{g})(p)\neq p$ we prove by induction that $p_{j}=p_{j}'$.\\

 \textbf{Claim.} For any $j\leq k_{V'}$,
$$(\dag )\ u_{V',j}(\bar{g})v_{V',j}(p'_{j})=u_{V',j}(\bar{g})v^{\varepsilon _{j}}_{V',j}(p_{j}).$$

 Note that if the claim holds then for $j=k_{V}'$ we have
$$w_{V'}(\bar{g})(p)=u_{V',k_{V'}}(\bar{g})v_{V',k_{V'}}(p_{k_{V'}}')=w_{V}(\bar{g})(p)\neq p.$$
This will finish the proof of the lemma.\parskip0pt

\emph{Proof of the claim.}\ We apply induction. Fix some $j\leq k_{V'}$ and assume that
$(\dag )$ holds for all $i<j$. It means that $p'_{j}=p_{j}$. If $\varepsilon _{j}=1$, then $u_{V',
j}(\bar{g})v_{V',j}=u_{V',j}(\bar{g})v^{\varepsilon _{j}}_{V',j}$ and we are done.\parskip0pt

 Now assume that $\varepsilon _{j}=0$. We claim that $p'_{j}\in X\setminus supp(v_{V',
j})$. Indeed, since
$$V=V'_{\bar{\varepsilon}}=\bigcap _{i=1}^{k_{V'}} v_{V',1}^{-1}\ldots v_{V',i-1}^{-1}
	\Big( supp(v_{V',i})^{\varepsilon _{i}}\Big) ,$$
for any $i\leq k_{V'}-1$, either $v_{V',i}\ldots v_{V',1}(V)\subseteq supp(v_{V',i+1})$ (when
$\varepsilon _{i+1}=1$) or $v_{V',i}\ldots v_{V',1}(V)\subseteq X\setminus supp(v_{V',i+1})$
(when $\varepsilon _{i+1}=0$). Since for any $i\leq t$ and any $j\leq k$, $g_{i}$ stabilizes
$v_{j}\ldots v_{1}(U)$ setwise and for any $j\leq k_{V'}$, $v_{V',j}$ is some product of the
original constants $v_{1},\ldots , v_{k}$, we see that for any $i\leq t$ and any $j\leq k_{V'}$,
$g_{i}$ stabilizes appropriate $v_{V',j}\ldots v_{V',1}(U)$. Thus
$$u_{V',j-1}(\bar{g})v_{V',j-1}\ldots u_{V',1}(\bar{g})v_{V',1}(U)\subseteq\Big( X\setminus
	supp(v_{V',j})\Big) .$$

 Since $p\in U$, $p'_{j}\in X\setminus supp(v_{V',j})$ and we obtain:
$$u_{V',j}(\bar{g})v_{V',j}(p'_{j})=u_{V',j}(\bar{g})(p'_{j})=u_{V',j}(\bar{g})v^{\varepsilon _{j}}
_{V',j}(p_{j}).$$

\ \ \ \ \ \ \ \ \ \ \ \ \ \ \ \ \ \ \ \ \ \ \ \ \ \ \ \ \ \ \ \ \ \ \ \ \ \ \ \ \ \ \ \ \ \ \ \ \ \ \ \ \ \ \ \ \ \ \ \ \ \ \ \ \ \ \ \ \ \ \ \ \ \ \ \ \ $\square$\\

 We now define one more class of words over $G$.

\begin{defi}
	If $v_{k}\ldots v_{1}\neq 1$ then we say that the word $w$ has \emph{non-trivial
	product of constants} (in the set $V_{w}:=supp(v_{k}\ldots v_{1})\subseteq X$).
\end{defi}

\begin{rema}
	 \emph{In terms of the paper \cite{KhM} the subgroup of $G\ast\mathbb{F}_{t}$
	consisting of words having trivial product of constants is the radical of the system
	of equations $x_{i}=1$, $1\leq i\leq t$.}
\end{rema}

\begin{exam} \label{1}
	\emph{ Consider the Thompson's group $F$ with its standard action on $[0,1]$. We now
	give several illustrations of notions introduced above. We start with two words}\parskip2pt

	$\bullet$ \ \ $w_{2}=x_{[0,\frac{1}{2}],0}^{-1}yx_{[\frac{1}{2},1],1}^{-1}y^{-1}x_{[0,\frac{1}{2}],1}
		yx_{[0,\frac{1}{2}],2}^{-1}$\parskip2pt

	$\bullet$ \ \ $w_{3}=y^{-1}x_{1}yx_{[0,\frac{1}{2}],0}y^{-1}x_{1}^{-1}yx_{[0,\frac{1}{2}],0}
		^{-1}$\emph{.}\parskip2pt

	 \emph{Word $w_{2}$ is not oscillating (because $x_{[0,\frac{1}{2}],2}x_{[0,\frac{1}{2}],1}
	^{-1}((\frac{1}{2},1))\cap (0,\frac{1}{2})=\emptyset$), but it is almost oscillating, whereas $w
	_{3}$ is rigid. To see this we apply the procedure described above.} \parskip0pt

	 \emph{In $w_{2}$ we have four constant segments: $v_{[0,1],1}=x_{[0,\frac{1}{2}],2}^{-1}$,
	$v_{[0,1],2}=x_{[0,\frac{1}{2}],1}$, $v_{[0,1],3}=x_{[\frac{1}{2},1],1}^{-1}$ and $v_{[0,1],4}
	=x_{[0,\frac{1}{2}],0}^{-1}$, which correspond to the sets:}\\

	$\ \ \ supp(v_{[0,1],1})^{1}=\Big(\frac{3}{8},\frac{1}{2}\Big)$ \emph{and} $supp(v_{[0,1],1})^{0}
		=\Big( 0,\frac{3}{8}\Big)\cup\Big(\frac{1}{2},1\Big)$,\\

	$\ \ \ supp(v_{[0,1],2})^{1}=\Big(\frac{1}{4},\frac{1}{2}\Big)$ \emph{and} $supp(v_{[0,1],2})^{0}=
		\Big( 0,\frac{1}{4}\Big)\cup\Big(\frac{1}{2},1\Big)$,\\

	$\ \ \ supp(v_{[0,1],3})^{1}=\Big(\frac{3}{4},1\Big)$ \emph{and} $supp(v_{[0,1],3})^{0}=\Big( 0,
		\frac{3}{4}\Big)$,\\

	$\ \ \ supp(v_{[0,1],4})^{1}=\Big( 0,\frac{1}{2}\Big)$ \emph{and} $supp(v_{[0,1],4})^{0}=\Big(
		\frac{1}{2},1\Big)$.\\

	\emph{Thus the set $\mathcal{P}^{1}$ for $w_{2}$ equals $\{ (0,\frac{1}{4}), (\frac{1}{4},
	\frac{3}{8}), (\frac{3}{8},\frac{1}{2}), (\frac{1}{2},\frac{3}{4}),(\frac{3}{4},1)\}$. Hence we
	obtain five reduced words: $w_{(0,\frac{1}{4})}=x_{[0,\frac{1}{2}],0}^{-1}y$, $w_{(\frac{1}{4},
	\frac{3}{8})}=x_{[0,\frac{1}{2}],0}^{-1}x_{[0,\frac{1}{2}],1}y$, $w_{(\frac{3}{8},\frac{1}{2})}
	=x_{[0,\frac{1}{2}],0}^{-1}x_{[0,\frac{1}{2}],1}yx_{[0,\frac{1}{2}],2}^{-1}$, $w_{(\frac{1}{2},
	\frac{3}{4})}=y$ and $w_{(\frac{3}{4},1)}=yx_{[\frac{1}{2},1],1}^{-1}y^{-1}$. The words $w
	_{(0,\frac{1}{4})}$, $w_{(\frac{1}{4},\frac{3}{8})}$, $w_{(\frac{1}{2},\frac{3}{4})}$ and $w
	_{(\frac{3}{4},1)}$ are oscillating. Note that $w_{(\frac{3}{8},\frac{1}{2})}$ is non-trivial
	and not oscillating, because $supp(x_{[0,\frac{1}{2}],0}^{-1}x_{[0,\frac{1}{2}],1})\cap
	supp(x_{[0,\frac{1}{2}],2}^{-1})=\emptyset$. This finishes the procedure and we see that
	$w_{2}$ is oscillating in the set $(0,\frac{3}{8})\cup (\frac{1}{2},1)$.}

	 \emph{On the other hand $w_{3}$ has also four constant segments: $v_{[0,1],1}=x_{[0,
	\frac{1}{2}],0}^{-1}$, $v_{[0,1],2}=x_{1}^{-1}$, $v_{[0,1],3}=x_{[0,\frac{1}{2}],0}$, $v_{[0,1],4}=x
	_{1}$, which correspond to the sets:}\\

	$\ \ \ supp(v_{[0,1],1})^{1}=\Big( 0,\frac{1}{2}\Big)$ \emph{and} $supp(v_{[0,1],1})^{0}=\Big(
		\frac{1}{2},1\Big)$,\\

	$\ \ \ supp(v_{[0,1],2})^{1}=\Big(\frac{1}{2},1\Big)$ \emph{and} $supp(v_{[0,1],2})^{0}=\Big( 0,
		\frac{1}{2}\Big)$.\\

	$\ \ \ supp(v_{[0,1],3})^{1}=\Big( 0,\frac{1}{2}\Big)$ \emph{and} $supp(v_{[0,1],3})^{0}=\Big(
		\frac{1}{2}, 1\Big)$,\\

	$\ \ \ supp(v_{[0,1],4})^{1}=\Big(\frac{1}{2},1\Big)$ \emph{and} $supp(v_{[0,1],4})^{0}=\Big(0 ,
		\frac{1}{2}\Big)$.\\

	\emph{The set $\mathcal{P}^{1}$ for $w_{3}$ equals $\{ (0,\frac{1}{2}),(\frac{1}{2},1)\}$.
	Thus in this case we obtain two words:}
$$w'_{(0,\frac{1}{2})}=y^{-1}yx_{[0,\frac{1}{2}],0}y^{-1}yx_{[0,\frac{1}{2}],0}^{-1}$$
\emph{and}
$$w'_{(\frac{1}{2},1)}=y^{-1}x_{1}yy^{-1}x_{1}^{-1}y .$$
\emph{In fact both of them are equal $\emph{id}$. Hence $W^{1}$ is empty and therefore $w_{3}$ is
	rigid.}\parskip0pt

	\emph{The fact that $w_{2}$ and $w_{3}$ are not oscillating can be read from the picture.
	For this purpose we look at the rectangles in Picture 8, which correspond to the
	constant subwords of $w_{2}$ (respectively, see Picture 9 for $w_{3}$) - 1st, 3rd, 5th and
	7th from the left. It is easy to see that the constants $x_{[0,\frac{1}{2}],2}^{-1}$, $x_{[0,
	\frac{1}{2}],1}$ and $x_{[0,\frac{1}{2}],0}^{-1}$ are trivial above the middle line (respectively,
	$x_{[0,\frac{1}{2}],0}^{-1}$ and $x_{[0,\frac{1}{2}],0}$ in the case of $w_{3}$). On the other
	hand $x_{[\frac{1}{2},1],1}^{-1}$ is trivial below the middle line (respectively, $x_{1}^{-1}$
	and $x_{1}$ in the case of $w_{3}$). Thus the images of the supports of this elements
	considered in the definition of an oscillating word are disjoint.}\parskip0pt

\setlength{\unitlength}{2mm}
\begin{picture}(45,45)
	\linethickness{0.05mm}
		\put(-5,9){\line(1,0){35}}
		\put(-5,41){\line(1,0){35}}
		\multiput(-5,9)(5,0){8}%
			{\line(0,1){32}}

			\put(-5,21){\line(1,0){5}}
			\put(-5,22){\line(5,1){5}}
			\put(-5,23){\line(5,1){5}}
			\put(-5,25){\line(1,0){5}}
			\put(-4.7,27){\tiny $x_{[0,\frac{1}{2}],2}^{-1}$}

			\put(2,29){\footnotesize $y$}

			\put(5,17){\line(1,0){5}}
			\put(5,21){\line(5,-2){5}}
			\put(5,23){\line(5,-2){5}}
			\put(5,25){\line(1,0){5}}
			\put(5.3,27){\tiny $x_{[0,\frac{1}{2}],1}$}

			\put(11.5,29){\footnotesize $y^{-1}$}

			\put(15,33){\line(1,0){5}}
			\put(15,35){\line(5,2){5}}
			\put(15,37){\line(5,2){5}}
			\put(15.3,27){\tiny $x_{[\frac{1}{2},1],1}^{-1}$}

			\put(22,29){\footnotesize $y$}

			\put(25,13){\line(5,4){5}}
			\put(25,17){\line(5,4){5}}
			\put(25,25){\line(1,0){5}}
			\put(25.3,27){\tiny $x_{[0,\frac{1}{2}],0}^{-1}$}

			\put(-2.2,7){\footnotesize $w_{2}=x_{[0,\frac{1}{2}],0}^{-1}yx_{[\frac{1}{2},1],1}
				^{-1}y^{-1}x_{[0,\frac{1}{2}],1}yx_{[0,\frac{1}{2}],2}^{-1}$}

	\put(-5,9){\circle*{0.5}}
	\put(-5,17){\circle*{0.5}}
	\put(-5,21){\circle*{0.5}}
	\put(-5,25){\circle*{0.5}}
	\put(-5,33){\circle*{0.5}}
	\put(-5,41){\circle*{0.5}}

	\put(0,9){\circle*{0.5}}
	\put(0,17){\circle*{0.5}}
	\put(0,21){\circle*{0.5}}
	\put(0,25){\circle*{0.5}}
	\put(0,33){\circle*{0.5}}
	\put(0,41){\circle*{0.5}}

	\put(5,9){\circle*{0.5}}
	\put(5,17){\circle*{0.5}}
	\put(5,21){\circle*{0.5}}
	\put(5,25){\circle*{0.5}}
	\put(5,33){\circle*{0.5}}
	\put(5,41){\circle*{0.5}}

	\put(10,9){\circle*{0.5}}
	\put(10,17){\circle*{0.5}}
	\put(10,19){\circle*{0.5}}
	\put(10,25){\circle*{0.5}}
	\put(10,33){\circle*{0.5}}
	\put(10,41){\circle*{0.5}}

	\put(15,9){\circle*{0.5}}
	\put(15,17){\circle*{0.5}}
	\put(15,19){\circle*{0.5}}
	\put(15,25){\circle*{0.5}}
	\put(15,33){\circle*{0.5}}
	\put(15,41){\circle*{0.5}}

	\put(20,9){\circle*{0.5}}
	\put(20,17){\circle*{0.5}}
	\put(20,19){\circle*{0.5}}
	\put(20,25){\circle*{0.5}}
	\put(20,33){\circle*{0.5}}
	\put(20,41){\circle*{0.5}}

	\put(25,9){\circle*{0.5}}
	\put(25,17){\circle*{0.5}}
	\put(25,19){\circle*{0.5}}
	\put(25,25){\circle*{0.5}}
	\put(25,33){\circle*{0.5}}
	\put(25,41){\circle*{0.5}}

	\put(30,9){\circle*{0.5}}
	\put(30,21){\circle*{0.5}}
	\put(30,22){\circle*{0.5}}
	\put(30,25){\circle*{0.5}}
	\put(30,33){\circle*{0.5}}
	\put(30,41){\circle*{0.5}}

	\linethickness{0.25mm}

			\put(-5,9){\line(0,1){8}}
			\put(-5,17){\line(0,1){4}}
			\put(-5,25){\line(0,1){8}}
			\put(-5,33){\line(0,1){8}}

			\put(0,9){\line(0,1){8}}
			\put(0,17){\line(0,1){4}}
			\put(0,25){\line(0,1){8}}
			\put(0,33){\line(0,1){8}}

			\put(5,9){\line(0,1){8}}
			\put(5,17){\line(0,1){4}}
			\put(5,25){\line(0,1){8}}
			\put(5,33){\line(0,1){8}}

			\put(10,9){\line(0,1){8}}
			\put(10,17){\line(0,1){2}}
			\put(10,25){\line(0,1){8}}
			\put(10,33){\line(0,1){8}}

			\put(15,9){\line(0,1){8}}
			\put(15,17){\line(0,1){2}}
			\put(15,25){\line(0,1){8}}
			\put(15,33){\line(0,1){8}}

			\put(20,9){\line(0,1){8}}
			\put(20,17){\line(0,1){2}}
			\put(20,25){\line(0,1){8}}
			\put(20,33){\line(0,1){8}}

			\put(25,9){\line(0,1){8}}
			\put(25,17){\line(0,1){2}}
			\put(25,25){\line(0,1){8}}
			\put(25,33){\line(0,1){8}}

			\put(30,9){\line(0,1){12}}
			\put(30,21){\line(0,1){1}}
			\put(30,25){\line(0,1){8}}
			\put(30,33){\line(0,1){8}}

	\linethickness{0.125mm}

			\multiput(-5,13)(2,0){15}%
				{\line(1,0){1}}

			\multiput(-5,17)(2,0){15}%
				{\line(1,0){1}}

			\multiput(-5,21)(2,0){5}%
				{\line(1,0){1}}
			\multiput(11,19)(2,0){7}%
				{\line(1,0){1}}
			\put(25,19){\line(5,3){5}}

			\multiput(1,23)(2,0){2}%
				{\line(1,0){1}}
			\multiput(11,21)(2,0){7}%
				{\line(1,0){1}}
			\put(25,21){\line(5,2){5}}

			\multiput(1,24)(2,0){2}%
				{\line(1,0){1}}
			\put(5,24){\line(5,-1){5}}
			\multiput(11,23)(2,0){7}%
				{\line(1,0){1}}
			\put(25,23){\line(5,1){5}}

			\multiput(-5,25)(2,0){18}%
				{\line(1,0){1}}

			\multiput(-5,33)(2,0){18}%
				{\line(1,0){1}}

			\multiput(-5,35)(2,0){10}%
				{\line(1,0){1}}
			\multiput(20,37)(2,0){5}%
				{\line(1,0){1}}

			\multiput(-5,37)(2,0){10}%
				{\line(1,0){1}}
			\multiput(20,39)(2,0){5}%
				{\line(1,0){1}}

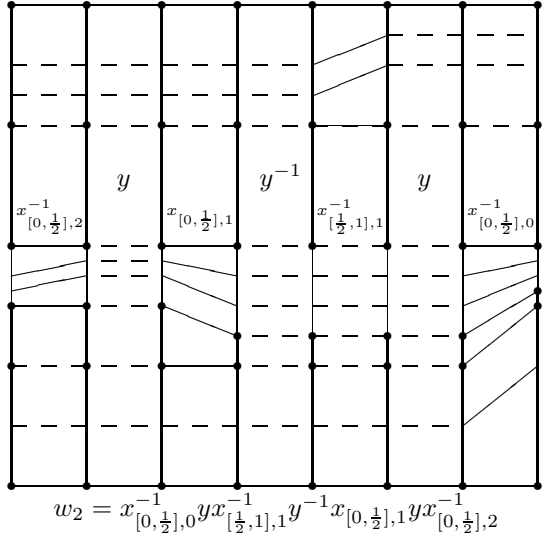
	\put(-5,4){{\footnotesize\emph{Picture 8.\ \ An example of a word, which is not}}}
 	\put(-1,1.5){{\footnotesize oscillating \emph{but is} almost oscillating}}

	\linethickness{0.05mm}
		\put(34,9){\line(1,0){32}}
		\put(34,41){\line(1,0){32}}
		\multiput(34,9)(4,0){9}%
			{\line(0,1){32}}

			\put(34,13){\line(1,1){4}}
			\put(34,17){\line(1,1){4}}
			\put(34,25){\line(1,0){4}}
			\put(34.2,27){\tiny $x_{[0,\frac{1}{2}],0}^{-1}$}

			\put(39.5,23){\footnotesize $y$}

			\put(42,25){\line(1,0){4}}
			\put(42,29){\line(1,1){4}}
			\put(42,33){\line(1,1){4}}
			\put(43.5,27){\tiny $x_{1}^{-1}$}

			\put(47,23){\footnotesize $y^{-1}$}

			\put(50,17){\line(1,-1){4}}
			\put(50,21){\line(1,-1){4}}
			\put(50,25){\line(1,0){4}}
			\put(50.2,27){\tiny $x_{[0,\frac{1}{2}],0}$}

			\put(55.5,23){\footnotesize $y$}

			\put(58,25){\line(1,0){4}}
			\put(58,33){\line(1,-1){4}}
			\put(58,37){\line(1,-1){4}}
			\put(59,27){\tiny $x_{1}$}

			\put(63,23){\footnotesize $y^{-1}$}

			\put(36.5,7){\footnotesize $w_{3}=y^{-1}x_{1}yx_{[0,\frac{1}{2}],0}
				y^{-1}x_{1}^{-1}yx_{[0,\frac{1}{2}],0}^{-1}$}

	\put(34,9){\circle*{0.5}}
	\put(34,25){\circle*{0.5}}
	\put(34,41){\circle*{0.5}}

	\put(38,9){\circle*{0.5}}
	\put(38,25){\circle*{0.5}}
	\put(38,41){\circle*{0.5}}

	\put(42,9){\circle*{0.5}}
	\put(42,25){\circle*{0.5}}
	\put(42,41){\circle*{0.5}}

	\put(46,9){\circle*{0.5}}
	\put(46,25){\circle*{0.5}}
	\put(46,41){\circle*{0.5}}

	\put(50,9){\circle*{0.5}}
	\put(50,25){\circle*{0.5}}
	\put(50,41){\circle*{0.5}}

	\put(54,9){\circle*{0.5}}
	\put(54,25){\circle*{0.5}}
	\put(54,41){\circle*{0.5}}

	\put(58,9){\circle*{0.5}}
	\put(58,25){\circle*{0.5}}
	\put(58,41){\circle*{0.5}}

	\put(62,9){\circle*{0.5}}
	\put(62,25){\circle*{0.5}}
	\put(62,41){\circle*{0.5}}

	\put(66,9){\circle*{0.5}}
	\put(66,25){\circle*{0.5}}
	\put(66,41){\circle*{0.5}}

	\linethickness{0.125mm}

			\multiput(34,25)(2,0){16}%
				{\line(1,0){1}}

			\multiput(38,17)(2,0){6}%
				{\line(1,0){1}}
			\multiput(54,13)(2,0){6}%
				{\line(1,0){1}}

			\multiput(38,21)(2,0){6}%
				{\line(1,0){1}}
			\multiput(54,17)(2,0){6}%
				{\line(1,0){1}}

			\multiput(34,29)(2,0){4}%
				{\line(1,0){1}}
			\multiput(46,33)(2,0){6}%
				{\line(1,0){1}}
			\multiput(62,29)(2,0){2}%
				{\line(1,0){1}}

			\multiput(34,33)(2,0){4}%
				{\line(1,0){1}}
			\multiput(46,37)(2,0){6}%
				{\line(1,0){1}}
			\multiput(62,33)(2,0){2}%
				{\line(1,0){1}}

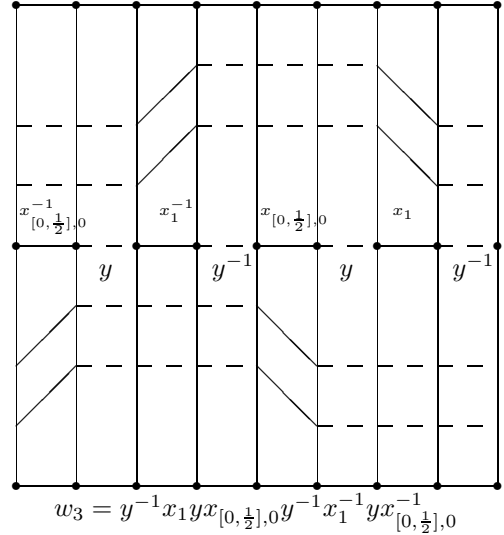
	\put(33.8,4){{\footnotesize\emph{Picture 9.\ \ An example of a word, which}}}
 	\put(37.8,1.5){{\footnotesize \emph{is} rigid}}

\end{picture}

	\emph{Deleting the 5th rectangle in Picture 8 corresponds to the operation of deleting
	all the constants in $w_{2}$ with support contained in $[\frac{3}{4},1]$. We see that after this
	operation we obtain a word, which is oscillating. Thus we may deduce that the original word
	was almost oscillating. On the other hand in the case of word $w_{3}$ no such operation
	produces an oscillating word. This implies that $w_{3}$ is rigid.}\parskip0pt

	\emph{There are words which are oscillating and do not have non-trivial product of
	constants. As an example consider the word}\parskip2pt

	$\bullet$ \ \ $w_{4}=yx_{1}y^{-1}x_{1}^{-1},$\parskip2pt

	\emph{which is ocillating since}
$$O_{w_{4}}=x_{1}\Big(\Big(\frac{1}{2},1\Big)\Big)\cap\Big(\frac{1}{2},1\Big) =\Big(\frac{1}{2},1\Big) .$$
\emph{It does not have non-trivial product of constants because $x_{1} x_{1}^{-1}=id$.}

	\emph{On the other hand the word}\parskip2pt

	$\bullet$ \ \ $w_{5}=yx_{1}y^{-1}x_{[0,\frac{1}{2}],0}y^{2}x_{1}^{-1}$\parskip2pt

	\emph{is not oscillating, because}
	$$O_{w_{5}}=x_{1}x_{[0,\frac{1}{2}],0}^{-1}\Big(\Big(\frac{1}{2},1\Big)\Big)\cap
		x_{1}\Big(\Big( 0,\frac{1}{2}\Big)\Big)\cap\Big(\frac{1}{2},1\Big) =
		\emptyset.$$
	\emph{Since $x_{1}x_{[0,\frac{1}{2}],0}x_{1}^{-1}=x_{[0,\frac{1}{2}],0}$ and
	$supp(x_{[0,\frac{1}{2}],0})=(0,\frac{1}{2})$, $w_{5}$ has non-trivial product of
	constants.}\parskip0pt

	\emph{The property of having non-trivial product of constants can be read from the
	picture as well. For this purpose note that in Picture 10 all the points on the left
	handside of the picture are connected with the points on the right handside situated
	on the same height, whereas in Picture 11 there are points on the left handside,
	which are connected with the points on the right handside situated below them.}

	\emph{Note that the word $w_{2}$ has non-trivial product of constants, while $w_{3}$
	does not have non-trivial product of constants. Combining the examples considered
	above we can also construct examples of words, which are rigid but have	non-trivial
	product of constants.}

\end{exam}

\setlength{\unitlength}{2mm}
\begin{picture}(45,45)
	\linethickness{0.05mm}
		\put(-5,9){\line(1,0){32}}
		\put(-5,41){\line(1,0){32}}
		\multiput(-5,9)(8,0){5}%
			{\line(0,1){32}}

			\put(-5,25){\line(1,0){8}}
			\put(-5,29){\line(2,1){8}}
			\put(-5,33){\line(2,1){8}}
			\put(-1.6,13){\footnotesize $x_{1}^{-1}$}

			\put(6.4,21){\footnotesize $y^{-1}$}

			\put(11,25){\line(1,0){8}}
			\put(11,33){\line(2,-1){8}}
			\put(11,37){\line(2,-1){8}}
			\put(14.5,13){\footnotesize $x_{1}$}

			\put(22.5,21){\footnotesize $y$}

			\put(4,7){\footnotesize $w_{4}=yx_{1}y^{-1}x_{1}^{-1}$}

	\put(-5,9){\circle*{0.5}}
	\put(-5,25){\circle*{0.5}}
	\put(-5,29){\circle*{0.5}}
	\put(-5,33){\circle*{0.5}}
	\put(-5,41){\circle*{0.5}}

	\put(3,9){\circle*{0.5}}
	\put(3,25){\circle*{0.5}}
	\put(3,33){\circle*{0.5}}
	\put(3,37){\circle*{0.5}}
	\put(3,41){\circle*{0.5}}

	\put(11,9){\circle*{0.5}}
	\put(11,25){\circle*{0.5}}
	\put(11,33){\circle*{0.5}}
	\put(11,37){\circle*{0.5}}
	\put(11,41){\circle*{0.5}}

	\put(19,9){\circle*{0.5}}
	\put(19,25){\circle*{0.5}}
	\put(19,29){\circle*{0.5}}
	\put(19,33){\circle*{0.5}}
	\put(19,41){\circle*{0.5}}

	\put(27,9){\circle*{0.5}}
	\put(27,25){\circle*{0.5}}
	\put(27,29){\circle*{0.5}}
	\put(27,33){\circle*{0.5}}
	\put(27,41){\circle*{0.5}}

	\linethickness{0.25mm}

			\put(-5,25){\line(0,1){16}}

			\put(3,25){\line(0,1){16}}

			\put(11,25){\line(0,1){16}}

			\put(19,25){\line(0,1){16}}

			\put(27,25){\line(0,1){16}}

	\linethickness{0.125mm}

			\multiput(-5,25)(2,0){16}%
				{\line(1,0){1}}

			\multiput(3,33)(2,0){4}%
				{\line(1,0){1}}
			\multiput(19,29)(2,0){4}%
				{\line(1,0){1}}

			\multiput(3,37)(2,0){4}%
				{\line(1,0){1}}
			\multiput(19,33)(2,0){4}%
				{\line(1,0){1}}

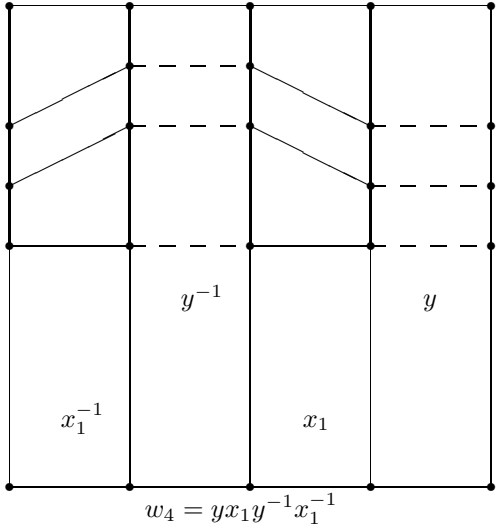
	\put(-5,4){{\footnotesize Picture 10.\ \ An \emph{oscillating} word, which does}}
 	\put(-3,1.5){{\footnotesize not have \emph{non-trivial product of constants}}}

	\linethickness{0.05mm}
		\put(30,9){\line(1,0){36}}
		\put(30,41){\line(1,0){36}}
		\multiput(30,9)(6,0){7}%
			{\line(0,1){32}}

			\put(30,25){\line(1,0){6}}
			\put(30,29){\line(3,2){6}}
			\put(30,33){\line(3,2){6}}
			\put(32.5,27){\tiny $x_{1}^{-1}$}

			\put(38.5,22){\footnotesize $y^{2}$}

			\put(42,17){\line(3,-2){6}}
			\put(42,21){\line(3,-2){6}}
			\put(42,25){\line(1,0){6}}
			\put(42.8,27){\tiny $x_{[0,\frac{1}{2}],0}$}

			\put(50,22){\footnotesize $y^{-1}$}

			\put(54,25){\line(1,0){6}}
			\put(54,33){\line(3,-2){6}}
			\put(54,37){\line(3,-2){6}}
			\put(56,27){\tiny $x_{1}$}

			\put(62,22){\footnotesize $y^{-1}$}

			\put(38.5,7){\footnotesize $w_{5}=yx_{1}y^{-1}x_{[0,
				\frac{1}{2}],0}y^{2}x_{1}^{-1}$}

	\put(30,9){\circle*{0.5}}
	\put(30,17){\circle*{0.5}}
	\put(30,21){\circle*{0.5}}
	\put(30,25){\circle*{0.5}}
	\put(30,29){\circle*{0.5}}
	\put(30,33){\circle*{0.5}}
	\put(30,41){\circle*{0.5}}

	\put(36,9){\circle*{0.5}}
	\put(36,17){\circle*{0.5}}
	\put(36,21){\circle*{0.5}}
	\put(36,25){\circle*{0.5}}
	\put(36,33){\circle*{0.5}}
	\put(36,37){\circle*{0.5}}
	\put(36,41){\circle*{0.5}}

	\put(42,9){\circle*{0.5}}
	\put(42,17){\circle*{0.5}}
	\put(42,21){\circle*{0.5}}
	\put(42,25){\circle*{0.5}}
	\put(42,33){\circle*{0.5}}
	\put(42,37){\circle*{0.5}}
	\put(42,41){\circle*{0.5}}

	\put(48,9){\circle*{0.5}}
	\put(48,13){\circle*{0.5}}
	\put(48,17){\circle*{0.5}}
	\put(48,25){\circle*{0.5}}
	\put(48,33){\circle*{0.5}}
	\put(48,37){\circle*{0.5}}
	\put(48,41){\circle*{0.5}}

	\put(54,9){\circle*{0.5}}
	\put(54,13){\circle*{0.5}}
	\put(54,17){\circle*{0.5}}
	\put(54,25){\circle*{0.5}}
	\put(54,33){\circle*{0.5}}
	\put(54,37){\circle*{0.5}}
	\put(54,41){\circle*{0.5}}

	\put(60,9){\circle*{0.5}}
	\put(60,13){\circle*{0.5}}
	\put(60,17){\circle*{0.5}}
	\put(60,25){\circle*{0.5}}
	\put(60,29){\circle*{0.5}}
	\put(60,33){\circle*{0.5}}
	\put(60,41){\circle*{0.5}}

	\put(66,9){\circle*{0.5}}
	\put(66,13){\circle*{0.5}}
	\put(66,17){\circle*{0.5}}
	\put(66,25){\circle*{0.5}}
	\put(66,29){\circle*{0.5}}
	\put(66,33){\circle*{0.5}}
	\put(66,41){\circle*{0.5}}

	\linethickness{0.125mm}

			\multiput(36,25)(2,0){3}%
				{\line(1,0){1}}
			\multiput(48,25)(2,0){3}%
				{\line(1,0){1}}
			\multiput(60,25)(2,0){3}%
				{\line(1,0){1}}

			\multiput(30,17)(2,0){6}%
				{\line(1,0){1}}
			\multiput(48,13)(2,0){9}%
				{\line(1,0){1}}

			\multiput(30,21)(2,0){6}%
				{\line(1,0){1}}
			\multiput(48,17)(2,0){9}%
				{\line(1,0){1}}

			\multiput(36,33)(2,0){9}%
				{\line(1,0){1}}
			\multiput(60,29)(2,0){3}%
				{\line(1,0){1}}

			\multiput(36,37)(2,0){9}%
				{\line(1,0){1}}
			\multiput(60,33)(2,0){3}%
				{\line(1,0){1}}

	\put(31.5,4){{\footnotesize Picture 11.\ \ A \emph{non-oscillating} word, which}}
 	\put(35.5,1.5){{\footnotesize has \emph{non-trivial product of constants}}}

\end{picture}

 The following theorem gives a sufficient condition for a system of inequalities over $G$ to have
a solution in $G$.

\begin{theo} \label{uab}
	 Let $G$ act on some locally compact perfect Polish space $X$ by
	homeomorphisms. Let $\{ w_{1}, w_{2}, \ldots , w_{m}\}$ be a set of words over $G$
	on $t$ variables, $y_{1},\ldots ,	y_{t}$, which are reduced and non-constant (i. e. $w
	_{i}\notin G$) in $\mathbb{F} _{t}\ast G$. Wlog assume that when $w_{j}\notin
	\mathbb{F} _{t}$, $w_{j}=u_{j,k_{j}}v_{j,k_{j}}\ldots u_{j,1}v_{j,1}$ for $1\leq j\leq m$
	with $v_{i}\in G\setminus\{ 1\}$. If $G$ hereditarily separates $X$ and each $w_{j}$,
	$j\leq m$, is either oscillating or almost oscillating or has non-trivial product of
	constants, then the set of inequalities $w_{1}\neq 1, w_{2}\neq 1,\ldots , w_{m}\neq
	1$ has a solution in $G$. \parskip0pt

	 Moreover, for any collection $\{ O_{j}\}$ such that $O_{j}$ is an open subset of
	the set $O_{w_{j}}$ (or $V_{w_{j}}$ respectively) defined for $w_{j}$ as above, $j
	\leq m$, there is a solution $(g_{1},\ldots , g_{t})$ of this set of inequalities such that
	$supp(g_{i})\subseteq\bigcup _{j=1} ^{m}\mathcal{V}_{w_{j}}(O_{j})$ for $1\leq i\leq t$.
\end{theo}

\emph{Proof.} \ Consider a metric $\rho$ such that $(X,\rho )$ becomes a locally compact
Polish metric space. Fix the collection $\{ O_{j}\}$ from the statement of the theorem. First for
each $j\leq m$ we will choose some non-empty open ball $B_{j}\subseteq O_{j}$ such that the
following conditions are satisfied:\\

$\bullet$ \ \ $B_{j}\cap B_{j'}=\emptyset$ for any $j,j'\in\{ 1,\ldots , m\}$, $j\neq j'$.\\

$\bullet$ \ \ $\mathcal{V} _{w_{j}}(B_{j})\cap \mathcal{V} _{w_{j'}}(B_{j'})=\emptyset$ for any $j,j'
	\in\{ 1,\ldots , m\}$, $j\neq j'$.\\

 We start by choosing a parwise disjoint collection of open balls $B_{1,1},\ldots ,$ $B_{m,1}$
such that $B_{1,j}\subseteq O_{j}$, $j\leq m$. Let
	$$R_{j}(A):=\max _{1\leq i\leq k_{j}} diam\Big( v_{j,i}\ldots v_{j,1}(A)\Big) .$$
 We now perform the following procedure. Choose a ball $B_{1,2,1}\subseteq B_{1,1}$ so that
$R_{1}(B_{1,2,1})< diam(v_{2,1}(B_{2,1}))$. Now (if necessary) replace $B_{2,1}$ with a
smaller ball $B'_{2,1}$ such that $B'_{2,1}\subseteq B_{2,1}$ and $\mathcal{V} _{w_{1}}(B_{1,2,1})
\cap v_{2,1}(B'_{2,1})=\emptyset$. We continue this process for all $s\leq k_{2}$, i. e.: at $s$-th
step we find a ball $B_{1,2,s}\subseteq B_{1,2,s-1}$ so that $R_{1}(B_{1,2,s})< diam(v_{2,s}\ldots
v_{2,1}(B_{2,1}))$ for $B_{2,1}$ updated at $s-1$-st step and (if necessary) we once again replace
$B_{2,1}$ with a smaller ball $B'_{2,1}$ satisfying $B'_{2,1}\subseteq B_{2,1}$ and $\mathcal{V}
_{w_{1}}(B_{1,2,s})\cap v_{2,s}\ldots v_{2,1}(B'_{2,1})=\emptyset$. After $k_{2}$ steps we obtain
a sequence $B_{1,1}\supset B_{1,2,1}\supset\ldots\supset B_{1,2,k_{2}}$ and define $B_{1,2}:=
B_{1,2,k_{2}}$. Note that the updated ball $B_{2,1}$ is a subset of the original one. It is clear that
for any $s\leq k_{2}$, $\mathcal{V} _{w_{1}}(B_{1,2})\cap v_{2,s}\ldots v_{2,1}(B_{2,1})=
\emptyset$ and hence $\mathcal{V} _{w_{1}}(B_{1,2})\cap\mathcal{V}_{w_{2}}(B_{2,1})=
\emptyset$.\parskip0pt

 We repeat this construction starting with $B_{1,2}$ in the obvious way to get the balls $B_{1,3}
\supset\ldots\supset B_{1,m}$. We define $B_{1}:=B_{1,m}$. Since for any $i\leq m$,
$\mathcal{V} _{w_{1}}(B_{1,i})\cap\mathcal{V}_{w_{i}}(B_{i,1})=\emptyset$, we see that for any $i
\leq m$, $\mathcal{V} _{w_{1}}(B_{1})\cap\mathcal{V}_{w_{i}}(B_{i,1})=\emptyset$.

 We now repeat this procedure for all $2\leq j\leq m$, starting with $B_{2,1}$. Since the obtained
collection $B_{2},\ldots , B_{m}$ satisfies $B_{i}\subseteq B_{i,1}$, all the previously satisfied
conditions remain true. Thus for any $j,j'\in\{ 1,\ldots , m\}$, $j\neq j'$, we have $\mathcal{V} _{w
_{j}}(B_{j})\cap \mathcal{V} _{w_{j'}}(B_{j'})=\emptyset$.\parskip0pt

 Now we construct a sequence $(\bar{g} _{0}, \bar{g} _{1},\ldots ,\bar{g} _{m})$, where
$\bar{g} _{i}=(g_{i,1},\ldots , g_{i,t})$, such that for a given $i$, $1\leq i\leq m$ the following
conditions are satisfied:\\

$\bullet\ \ supp(\bar{g _{i}})\subseteq\bigcup _{j=1}^{i}\mathcal{V} _{w_{j}}(B_{j})$,\\

$\bullet\ \ \forall 1\leq j\leq t\ \forall 1\leq i, l\leq m,\ l\neq i\ (g_{i,j}|_{B_{l}}=g_{i-1,j}|_{B_{l}})$,\\

$\bullet\ \ \bar{g}_{i}$ is a solution of the set of inequalities $w_{1}\neq 1,\ldots , w_{i}\neq 1$.\\

  Fix $\bar{g}_{0}:=(1,\ldots , 1)\in G^{t}$. At $n$-th step, $n\leq m$, we will modify the action of
elements from the tuple $\bar{g}_{n-1}$ on the set $B_{n}$. At this step we want $\bar{g}_{n}$
to satisfy $w_{n}(\bar{g}_{n})|_{B_{n}}\neq 1|_{B_{n}}$ (i.e. $\bar{g}_{n}$ is a solution of the
inequality $w_{n}\neq 1$).

 Suppose we have made the first $n-1$ steps and defined the tuple $\bar{g}_{n-1}=(g_{n-1,
1},\ldots , g_{n-1,t})$. For the $n$-th word $w_{n}$ we consider three cases:\\

\textbf{Case 1.} $w_{n}$ has non-trivial product of constants. \parskip4pt

 For each $i\leq t$ we define $g_{n,i}:=g_{n-1,i}$. Thus for each $i\leq t$, we have $g_{n,i}|
_{B_{n}}=g_{n-1,i}|_{B_{n}}=1|_{B_{n}}$. To see that $w_{n}(\bar{g}_{n})\neq 1$ choose any
$p\in B_{n}$. As $B_{n}\subseteq O_{n}\subseteq V_{w_{n}}$, we have $v_{n,k_{n}}v_{n,k_{n}-1}
\ldots v_{n,1}(p)\neq p$. Since for any $j\leq k_{n}$, $u_{n,j}$ depends only on variables and
for all $i\leq t$, $g_{n,i}$ are trivial on the set $\mathcal{V}_{w_{n}}(B_{n})$, we obtain:
$$w_{n}(\bar{g}_{n})(p)=v_{n,k_{n}}v_{n,k_{n}-1}\ldots v_{n,1}(p)\neq p.$$
Thus $\bar{g}_{n}$ is a solution of the inequality $w_{n}\neq 1$.\\

\textbf{Case 2.} $w_{n}$ is oscillating. \parskip4pt

 We apply Theorem \ref{ab} to the word $w_{n}$ and the set $B_{n}\subseteq O_{w_{n}}$ and
obtain some solution of the inequality $w_{n}\neq 1$, $\bar{f}=(f_{1},\ldots , f_{t})\in G^{t}$,
such that $supp(f_{i})\subseteq \mathcal{V}_{w_{n}}(B_{n})$, $1\leq i\leq t$. Now for each $i\leq
t$ we define $g_{n,i}:=f_{i}g_{n-1,i}$. Thus $\bar{g}_{n}$ is also the solution of $w_{n}\neq 1$.
Since
$$\bigcup _{i=1} ^{t} supp(f_{i})\subseteq\mathcal{V}_{w_{n}}(B_{n})\ \land\ \bigcup_{i=1}
	^{t} supp(g_{n-1,i})\subseteq\Big(X\setminus\mathcal{V}_{w_{n}}(B_{n})\Big) ,$$
the tuple $(g_{n,1},\ldots , g_{n,t})$ still is a solution of $w_{j}\neq 1$ for $1\leq j\leq n-1$.\\

\textbf{Case 3.} $w_{n}$ is almost oscillating. \parskip4pt

 If $w_{n}$ is almost oscillating then it follows from the definition that the set $\mathcal{P}^{os}$
corresponding to the word $w_{n}$ is non-empty. Thus there is some word $w_{U}$ derrived
from $w_{n}$ by cancelations of constants and reductions, which is oscillating. By Theorem
\ref{ab} there is some $\bar{f}$ such that for any $i\leq t$, $supp(f_{i})\subseteq \mathcal{V}_{w
_{n}}(B_{n})$, for any $j\leq k_{n}$, $f_{i}$ stabilizes $v_{j}\ldots v_{1}(U)$ setwise and $w_{U}
(\bar{f})(p)\neq p$ for some point $p\in U$. Thus by Lemma \ref{ii} we have $w_{n}(\bar{f})\neq
1$.\parskip0pt

 Now similarly as above for any $i\leq t$ we define $g_{n,i}:=f_{i}g_{n-1,i}$. This gives a
solution of all inequalities $w_{j}\neq 1$ for $j\leq n$.\\

 Thus, after $m$ steps of the algorithm we obtain a tuple $\bar{g} _{m}$, which is the solution of
the set $w_{1}\neq 1,\ldots , w_{m}\neq 1$. Moreover, for any $1\leq i\leq t$, $supp(g
_{i})\subseteq\bigcup _{j=1} ^{m}\mathcal{V}_{w_{j}}(O_{j})$. \\

\ \ \ \ \ \ \ \ \ \ \ \ \ \ \ \ \ \ \ \ \ \ \ \ \ \ \ \ \ \ \ \ \ \ \ \ \ \ \ \ \ \ \ \ \ \ \ \ \ \ \ \ \ \ \ \ \ \ \ \ \ \ \ \ \ \ \ \ \ \ \ \ \ \ \ \ \ $\square$\\

\section{Applications}

 Using Theorem \ref{uab} we can find several interesting limits with respect to groups having
hereditarly separating action.\\

\subsection{Relative limit groups}

 As the first example we have the following proposition, which is a straightforward consequence of
Theorem \ref{uab}.

\begin{prop} \label{pl}
	 Let $t,q\in\mathbb{N}\setminus\{ 0\}$. Suppose that $G=\langle h_{1},\ldots , h_{q}\rangle$
	has a hereditarily separating action on some locally compact perfect Polish space $X$. Let
	$W$ be the set of all non-trivial words from $\mathbb{F} _{t}\ast G$, which are oscillating or
	almost oscillating or have non-trivial products of constants. Then there exist a sequence of
	tuples $(\bar{g}_{n})_{n<\omega}=((g_{n,1},\ldots , g_{n,t}))_{n<\omega}$ from $G$ such
	that $G_{n}=((G, (g_{n,1},\ldots , g_{n,t}, h_{1},\ldots , h_{q})))_{n<\omega}$ is a convergent
	sequence of marked groups such that for any $w\in W$ the inequality $w(\bar{g}_{n})\neq 1$
	is satisfied in almost all $G_{n}$.
\end{prop}

\emph{Proof.} \ Let $W:=\{ w_{0}, w_{1},\ldots\}$. We construct a sequence of marked groups $(G_{n})
_{n<\omega}$, $G_{n}=\langle g_{n,1},\ldots , g_{n,t}, h_{1},\ldots , h_{q}\rangle =G$, where $(g_{n,1},
\ldots , g_{n,t})$ is a solution in $G$ of the system of inequalities $w_{1}\neq 1, w_{2}\neq 1,\ldots , w
_{n}\neq 1$. Since $G$ has hereditarily separating action, it follows from Theorem \ref{uab}, that such
a tuple always exists.\\

\ \ \ \ \ \ \ \ \ \ \ \ \ \ \ \ \ \ \ \ \ \ \ \ \ \ \ \ \ \ \ \ \ \ \ \ \ \ \ \ \ \ \ \ \ \ \ \ \ \ \ \ \ \ \ \ \ \ \ \ \ \ \ \ \ \ \ \ \ \ \ \ \ \ \ \ \ $\square$\\

\begin{exam}
\emph{ Let $G$ act on some locally compact perfect Polish space $X$ by homeomorphisms.
Fix $h_{1},\ldots , h_{k}, h\in G$ and denote by $W=\{ w_{1}, w_{2},\ldots\}$ the set of words over $G$
with $t$ variables $y_{1},\ldots , y_{t}$, which are reduced, non-constant (i.e. $w_{i}\notin G$) in
$\mathbb{F}_{t}\ast G$ and have all constants in $\langle h\rangle$. Thus each $w_{i}$ is
oscillating.}\parskip0pt

\emph{ We now construct a sequence $(\bar{g}_{n})_{n<\omega}$ of tuples $\bar{g}_{n}=(g_{n,1},
\ldots , g_{n,t})_{n<\omega}$ of elements from $G$. For any $n\geq 1$, let $(g_{n,1},\ldots , g_{n,t})$
be the solution of the set of inequalities $w_{1}\neq 1,\ldots , w_{n}\neq 1$ given by Theorem \ref{uab}.
For any $n\in\mathbb{N}$ denote by
$$G_{n}:=\Big(\Big\langle h_{1},\ldots , h_{k}, h, g_{n,1},\ldots , g_{n,t}\Big\rangle , \Big( h_{1},\dots , h
	_{k}, h, g_{n,1},\ldots , g_{n,t}\Big)\Big)$$
a marked group given by the tuple $(h_{1},\ldots , h_{k}, h, g_{n,1},\ldots , g_{n,t})$. Since the space
$\mathcal{G} _{t}$ is compact, there is a convergent subseqence of the sequence of groups $(G_{n})
_{n<\omega}$. It follows from the construction given above that the corresponding limit group contains
$\langle h\rangle\ast\langle g_{1},\ldots , g_{t}\rangle$, where $g_{1},\ldots , g_{t}$ are "limits" of $(g
_{n,1})_{n<\omega},\ldots , (g_{n,1})_{n<\omega}$, respectively.}\parskip0pt

\emph{In the special case, when $G$ is torsion-free and $t=1$, we obtain an embedding of $\mathbb{F}
_{2}$ into $\lim _{n\to\infty} G_{n}$. In particular it is true for the Thompson's group $F$.}

\emph{In the latter case we can develop this idea as follows. Let $W_{0}$ be the set of words over $F$
with $t$ variables, which are reduced, non-constant in $\mathbb{F} _{t}\ast F$ and have all contants in
$\langle x_{0}\rangle$; similarly, let $W_{1}$ be the set of all such words with constants in $\langle x_{1}
\rangle$. Now for any $n\in\mathbb{N}$ let $(g_{n,1},\ldots ,g_{n,t})$ be the solution of the system of
inequalities $w_{1}\neq 1,\dots , w_{n}\neq 1$, where $w_{1},\ldots , w_{n}\in W:=W_{0}\cup W_{1}$.
We consider a sequence $(G_{n})_{n<\omega}$, where
$$G_{n}:=\Big(\Big\langle x_{0}, x_{1}, g_{n,1},\ldots g_{n,t}\Big\rangle , \Big( x_{0}, x_{1}, g_{n,1},\ldots ,
	g_{n,t}\Big)\Big) .$$
 Similarly as above, we see that the corresponding limit group $\lim _{n\to\infty} G_{n}$ marked by a tuple
$(x_{0}, x_{1},g_{1},\ldots , g_{t})$ contains both, $\langle x_{0}\rangle\ast\langle g_{1},\ldots , g_{t}\rangle$
and $\langle x_{1}\rangle\ast\langle g_{1},\ldots , g_{t}\rangle$.}
\end{exam}

 Proposition \ref{pl} leads to the following definition.

\begin{defi}\label{w}
	Fix a group $G$. Let $W$ be any class of words over $G$ with $m$ variables, which are
	reduced and non-trivial in $\mathbb{F}_{t}\ast G$. Let $(G_{n})_{n<\omega}$, $G<
	G_{n}$, be a sequence of marked groups, where $G_{n}:=(\langle g_{n,1},\dots g_{n,m}
	\rangle , (g_{n,1},\dots g_{n,m}))$. We say that $\mathbb{G}$ is a \emph{$W$-limit group}
	of the sequence $(G_{n})_{n<\omega}$ if the following condition is satisfied for any $w\in W$
	$$\forall _{n}^{\infty}\ G_{n}\models\ w(\bar{g}_{n})=1\ \iff\ \mathbb{G}\models\ w(\bar{g})
		=1.$$
\end{defi}

\begin{coro}
	Let $W$ be a class of words over the Thompson's group $F$ with $m$ variables containing
	only words, which are oscillating, almost oscillating or have non-trivial products of constants
	with respect to the natural action of $F$ on the unit interval. Then there exists a sequence of
	groups $(G_{n})_{n<\omega}$, where
	$$G_{n}:=\Big(\Big\langle g_{n,1},\ldots , g_{n,m}, x_{0}, x_{1}\Big\rangle , \Big( g_{n,1},\ldots ,
		g_{n,m}, x	_{0}, x_{1}\Big)\Big)$$
	for $g_{n,i}\in F$, such that $\langle g_{1},\ldots g_{m}\rangle\ast F$ is its $W$-limit.\parskip0pt

	Similarly, let $W$ be a class of words over the Grigorchuk group $G$ with $m$ variables
	containing only words, which are oscillating, almost oscillating or have non-trivial products of
	constants with respect to the natural action of $G$ on the infinite rooted binary tree. Then
	there exists a sequence of groups $(G_{n})_{n<\omega}$, where
	$$G_{n}:=\Big(\Big\langle g_{n,1},\ldots , g_{n,m}, a, b, c, d\Big\rangle , \Big( g_{n,1},\ldots , g_{n,
		m}, a, b, c,d\Big)\Big)\footnote{$a$, $b$, $c$ and $d$ are used to denote the standard
		generators of the first Grigorchuk group (see Introduction for an explicit presentation)}$$
	for $g_{n,i}\in G$, such that $\langle g_{1},\ldots g_{m}\rangle\ast G$ is its $W$-limit.
\end{coro}

\emph{Proof.} \ It follows from Examples \ref{f} and \ref{grig} that both, the action of Thompson's group
on $[0,1]$ and the action of Grigorchuk group on the boundary of infinite rooted binary tree, are
hereditarily separating. Hence we may apply Proposition \ref{pl}. This finishes the proof.\\

\ \ \ \ \ \ \ \ \ \ \ \ \ \ \ \ \ \ \ \ \ \ \ \ \ \ \ \ \ \ \ \ \ \ \ \ \ \ \ \ \ \ \ \ \ \ \ \ \ \ \ \ \ \ \ \ \ \ \ \ \ \ \ \ \ \ \ \ \ \ \ \ \ \ \ \ \ $\square$\\

In the case of Thompson's group we also have the following Proposition.

\begin{prop}\label{iva}
	 Let $F=\langle x_{0}, x_{1}\rangle$ be the Thompson's group and let $H=\langle h_{1},\ldots
	 , h_{t}\rangle$ be a $t$-generated subgroup in $F$. There is an $F$-limit group $G_{H}$,
	which is a homomorphic image of $H\ast F=\langle\hat{\gamma} _{1},\ldots , \hat{\gamma} _{t},
	\hat{x}_{0}, \hat{x}_{1}\rangle$ such that the words, which have non-trivial product of constants,
	are not in the kernel \footnote{i. e. the kernel is a subgroup of the radical of the system
	$x_{i}=1$, $i\leq t$}.
\end{prop}

\emph{Proof.} \ First we will show that any system of inequalities over $F$, $w_{1}\neq\emph{id}, \ldots ,
w_{m}\neq\emph{id}$, where $w_{i}$, $1\leq i\leq m$, are non-constant words on $t$ variables with
non-trivial products of constants, has a solution $(g_{1},\ldots , g_{t})\in F$ such that $\langle g_{1},\ldots
, g_{t}\rangle\cong H$. \parskip0pt

 Fix $w_{1},\ldots w_{m}$. Similarly as earlier wlog we consider only words $w_{j}$, $1\leq j\leq m$, which
can be written in the following reduced and non-degenerated form:
$$\bar{u}_{j,k_{j}}v_{j,k_{j}}\ldots \bar{u}_{j,2}v_{j,2}\bar{u}_{j,1}v_{j,1},$$
where $\bar{u}_{j,i}$ depend only on the letters $y_{1},\ldots ,y_{t}$ and $v_{j,i}\in F\setminus\{\emph{id}\}$,
$1\leq i\leq k_{j}$, where $k_{j}\in\mathbb{N}\setminus\{ 0\}$, $1\leq j\leq m$. \parskip0pt

 Since $w_{1},\ldots w_{m}$ have non-trivial product of constants, we choose points $p_{1},\ldots p
_{m}\in [0,1]$ such that for all $j$, $1\leq j\leq m$, $w_{j}(\emph{id},\ldots , \emph{id})(p_{j})\neq p_{j}$. Next
we choose some dyadic, non-degenerated and closed interval
$$U\subseteq \Big( [0,1]\setminus \bigcup _{j=1}^{m}\Big(\bigcup _{i=0} ^{k_{j}} (w_{j})_{i}(\emph{id},\ldots ,
	\emph{id})(p_{j})\Big)\Big) .$$

 Now we will define the tuple $\bar{g}=(g_{m,1},\ldots , g_{m,t})\in F^{t}$, such that $\bar{g}$ is a
solution of the set of inequalities $w_{1}\neq 1, \ldots , w_{m}\neq 1$ and generates a subgroup of
$F$ isomorphic to $H$. According to the notation from the preliminary section, for each $1\leq i\leq t$
we define $g_{m,i}$ as follows:
$$g_{m,i}=(f_{i})_{U}.$$
Obviously $H_{U}=\langle g_{m,1},\ldots , g_{m,t}\rangle<F_{U}$ is a group isomorphic with the
original $H$. Moreover for each $j$, $1\leq j\leq m$,
$$w_{j}(g_{m,1},\ldots , g_{m,t})(p_{j})=\bar{u}_{j,k_{j}}v_{j,k_{j}}\ldots \bar{u}_{j,2}v_{j,2}\bar{u}_{j,1}v
_{j,1}(g_{m,1},\ldots , g_{m,t})(p_{j})=$$
$$=v_{j,k_{j}}\ldots v_{j,2}v_{j,1}(g_{m,1},\ldots , g_{m,t})(p_{j})\neq p_{j}.$$
 Thus  $(g_{m,1},\ldots , g_{m,t})$ is a solution of the considered set of inequalities. \parskip0pt

 Finally enumerate all non-constant words on $t$ variables and non-trivial product of constants: $w_{1},
\ldots , w_{n},\ldots$. We construct the sequence of groups $(G_{m})_{m<\omega}$ by defining $G
_{m}=(F, (g_{m,1},\ldots , g_{m,t}, x_{0}, x_{1}))$, where $g_{m,1},\ldots , g_{m,t}$ is defined as above
with respect to $w_{1},\ldots , w_{m}$. Since $\langle g_{m,1},\ldots , g_{m,t}\rangle\cong H$ for each
$m$, we see that the limit of sequence $(G_{m})_{m<\omega}$ satisfies the statement of the
proposition.\\

\ \ \ \ \ \ \ \ \ \ \ \ \ \ \ \ \ \ \ \ \ \ \ \ \ \ \ \ \ \ \ \ \ \ \ \ \ \ \ \ \ \ \ \ \ \ \ \ \ \ \ \ \ \ \ \ \ \ \ \ \ \ \ \ \ \ \ \ \ \ \ \ \ \ \ \ \ $\square$\\

We call a weakly branch group $G$ acting on a rooted tree $T$ \emph{self-reproducing} if for any
$n\in\mathbb{N}$, the pointwise stabilizer of the subtree $T_{n}$ consisting of the first $n$ levels
induces $G$ on the subtree $T_{v}$ for every vertex $v$ of the $n$-th level. In fact, the argument of
Proposition \ref{iva} also works in the case of self-reproducing branch groups.

\begin{prop}
	 Let $G=\langle g_{1},\ldots , g_{n}\rangle$ be a self-reproducing weakly branch group with
	respect to the action on a rooted tree $T$. Let $H=\langle h_{1},\ldots , h _{t}\rangle$ be a
	$t$-generated subgroup in $G$. There is an $G$-limit group $G_{H}$, which is a
	homomorphic image of $H\ast G=\langle\hat{\gamma}_{1},\ldots , \hat{\gamma} _{t},\hat{g}
	_{1},\ldots ,\hat{g}_{n}\rangle$ such that the words, which have non-trivial product of
	constants, are not in the kernel.
\end{prop}

\emph{Proof.} \ Similarly as in the proof of Proposition \ref{iva}, for a given set of inequalities $w_{1}\neq
1,\ldots , w_{m}\neq 1$ we find a solution $\bar{f}=(f_{1},\dots , f_{t})\in G$ such that $\langle f_{1},\ldots
, f_{t}\rangle\cong H$.\parskip0pt

 Since $w_{1}, \ldots , w_{m}$ have non-trivial product of constants, we can find the corresponding set of
vertices $v_{1},\ldots ,v_{m}$ in $T$ such that for any $j\leq m$, $v_{j}$ is not stabilized by $w_{j}(1,\ldots
,1)$. By $k$ denote the level of the tree $T$ such that $v_{1},\ldots , v_{m}$ are contained within the first
$k$ levels of $T$. Now choose any $v\in T$ on the $k'$-th level for some $k'>k$. It is easy to see that for
any tuple of elements $\bar{f}=(f_{1},\ldots , f_{t})$ fixing $T_{k}\cup\{ v\}$ and any $j\leq m$ we have
$w_{j}(\bar{f})(v_{j})\neq v_{j}$. Since $G$ is self-reproducing, the stabilizer $G_{v}$ is isomorphic to $G$,
and thus we may choose $\bar{f}$ so that $\langle f_{1},\ldots , f_{t}\rangle\cong H$.\parskip0pt

 We finish the proof exactly as in the case of Proposition \ref{iva}, i. e. by choosing a sequence of tuples
$(\bar{f}_{m})_{m<\omega}=((f_{m,1},\ldots , f_{m,t}))_{m<\omega}$ such that for any $m\in\mathbb{N}$,
$\bar{f} _{m}$ is a solution of the system $w_{1}\neq 1,\ldots , w_{m}\neq 1$ and $\langle f_{m,1},\ldots ,
 f_{m,t}\rangle\cong H$.\\

\ \ \ \ \ \ \ \ \ \ \ \ \ \ \ \ \ \ \ \ \ \ \ \ \ \ \ \ \ \ \ \ \ \ \ \ \ \ \ \ \ \ \ \ \ \ \ \ \ \ \ \ \ \ \ \ \ \ \ \ \ \ \ \ \ \ \ \ \ \ \ \ \ \ \ \ \ $\square$\\

 Since the first Grigorchuk group is an example of a self-reproducing weakly branch group (see \cite{Gri}
for details\footnote{we use the term self-reproducing in a slightly stronger form than in \cite{Gri}, but this
place is correct}), we have the following corollary.

\begin{coro}
	 Let $G=\langle a, b, c, d\rangle$ be the first Grigorchuk group and let $H=\langle h_{1},\ldots
	 , h_{t}\rangle$ be a $t$-generated subgroup in $G$. There is an $G$-limit group $G	_{H}$,
	which is a homomorphic image of $H\ast G=\langle\hat{\gamma} _{1},\ldots , \hat{\gamma}
	_{t},\hat{a}, \hat{b}, \hat{c}, \hat{d}\rangle$ such that the words, which have non-trivial product of
	constants,	are not in the kernel.
\end{coro}

\subsection{HNN-extension of $F$ as an $F$-limit group}

 The most interesting application of Theorem \ref{uab} is as follows. Consider the particular case of
Thompson's group $F$. We will refer to the elements of $F$ as to the homeomorphisms of $[0,1]$.
Consider the following subgroup of $F$
$$H_{1}:=\Big\{ f\in F\ \Big|\ 1\notin\overline{supp}(f)\Big\} .$$
We construct an example of an HNN-extension of $F$ over $H_{1}$, which can be obtained as an
$F$-limit group. In fact $H_{1}$ is not a finitely generated subgroup and, as we will later show (see:
Theorem \ref{mt}), it is a necessary condition.

\begin{theo}\label{h1}
	There is a converegent sequence of groups, $(G_{n})_{n<\omega}$, where $G_{n} = (F,
	(g_{n}, x_{0}, x_{1}))$, $g_{n}\in F$ for $n<\omega$, such that
	$$\lim_{n\to\infty}(G_{n})=\Big( F\ast _{H_{1}}, \Big( g,x_{0}, x_{1}\Big)\Big) :=\Big(\Big
		\langle F, g\ \Big|\ ghg^{-1} =h\ h\in H_{1}\Big\rangle , \Big( g,x_{0}, x_{1}\Big)
		\Big) .$$
\end{theo}

\emph{Proof.}\ Let $W=\{ w_{1}, w_{2},\ldots\}$ be a set of words over $F$ with one variable $y$,
which are reduced, non-constant (i.e. $w_{i}\notin F$) in $F\ast\mathbb{F}$ and have all constants
in $F\setminus H_{1}$. Wlog assume that when $w_{j}\notin\mathbb{F}$, $w_{j}=y^{a_{j,k_{j}}}v_{j,k
_{j}}\ldots y^{a_{j,1}}v_{j,1}$ for $j\geq 1$ with $v_{j,i}\in F\setminus H_{1}$. Now for any $w_{j}\in
W$, for all $i\leq k_{j}$, $1\in\overline{supp}(v_{j,i})$. Thus all the constants $v_{j,1},\ldots , v_{j,k
_{j}}$ are non-trivial in some sufficiently small neighbourhood of $1$. It follows from the definition
that for any $j\geq 1$, $w_{j}$ is oscillating. Let also $Z:=\{z_{1}, z_{2}, \ldots\}$ be an enumeration
of $H_{1}$.\parskip0pt

We now construct a sequence $(g_{n})_{n<\omega}$. Fix $n\in\mathbb{N}$. Since all elements
from $H_{1}$ are trivial in some neighbourhood of $1$ we can find some number $N\in\mathbb{N}$
such that for any $z\in Z_{n}:=\{ z_{1},\ldots , z_{n}\}$, we have $supp(z)\subseteq [0,1-\frac{1}{2^{N}}]$.
On the other hand since all the constants in the words from the set $W_{n}:=\{ w_{1},\ldots , w_{n}\}$
are non-trivial in some neighbourhood of $1$, we know that enlarging $N$ if necessary we can
satisfy $[1-\frac{1}{2^{N}},1]\subseteq\bigcap _{w\in W_{n}}O_{w}$ ($O_{w}$ denotes the set, where
$w$ is oscillating). We now apply Theorem \ref{uab} to the set of oscillating words $W_{n}$ and
obtain a solution $g_{n}$ of the set of inequalities $w_{1}\neq\emph{id},\ldots , w_{m}\neq\emph{id}$.
It also follows from Theorem \ref{uab} that we may choose a solution $g_{n}$ such that $supp(g_{n})
\subseteq (1-\frac{1}{2^{N}},1)$. In this way we construct the whole sequence $(g_{n})_{n<
\omega}$.\parskip0pt

Since for any $h\in H_{1}$ there is $n_{h}\in\mathbb{N}$ such that $supp(h)\subseteq [0,1-\frac{1}{2^{n
_{h}}}]$, we see that for any $n>n_{h}$, the supports of $g_{n}$ and $h$ are disjoint. Thus for any
$h$ and any $n>n_{h}$, $[g_{n},h]=\emph{id}$. Therefore for any $h\in H_{1}$, the relation $[g,h]=
\emph{id}$ is also satisfied in the limit group. Hence $\lim _{n\to\infty}(G_{n})$ is an homomorphic
image of $F\ast _{H_{1}}$.\parskip0pt

Now let $w$ be any word over $F$ with one variable, which is reduced and non-constant in
$F\ast\mathbb{Z}$. Suppose that $F\ast _{H_{1}}\models\ w(g)\neq\emph{id}$. We will show that $\lim
_{n\to\infty}(G_{n})\models\ w(g)\neq\emph{id}$.\parskip0pt

First suppose that $w\in\mathbb{Z}$. Since $w$ is non-trivial, for almost all $n$, $g_{n}\neq\emph{id}$.
Hence we have $\lim _{n\to\infty}(G_{n})\models w(g)\neq\emph{id}$.\parskip0pt

Now assume that $w(g)$ has the following form: $w(g)=g^{a_{k}}v_{k}\ldots g^{a_{1}}v_{1}$, where
$g$ is the letter defining the new generator of $F_{H_{1}}$ and for all $i\leq k$, $a_{i}\in\mathbb{Z}
\setminus\{ 0\}$ and $v_{i}\in F$. Assume that $w(g)$ is reduced in $F\ast_{H_{1}}$. It follows from
Britton's Lemma on irreducible words in HNN-extensions (\cite{LS}, page 181) that for any $v\in H
_{1}$, neither subwords of the form $gvg^{-1}$ nor of the form $g^{-1}vg$ occur in
$w(g)$. Moreover, for any $i>1$, $v_{i}\in F\setminus H_{1}$.\parskip0pt

 Consider $w(g_{n})$. Suppose that for all $i$, $1\leq i\leq k$, $v_{i}\in F\setminus H_{1}$. Then
there is some $N\in\mathbb{N}$ such that for all $n>N$, the word of the form $w=y^{a_{k}}v_{k}
\ldots y^{a_{1}}v_{1}$ belongs to $W_{n}$. Thus it follows from the construction of the sequence
$(g_{n})_{n<\omega}$ that for all $n>N$, $w(g_{n})\neq\emph{id}$.\parskip0pt

 Now suppose that $v_{1}\in H_{1}$. Consider the case where there is at least one constant from
$F\setminus H_{1}$ in $w(g)$, i. e. $v_{2}\neq\emph{id}$. Let $N_{w}$ be large enough to satisfy
$supp(v_{1})\subseteq [0,1-\frac{1}{2^{N_{w}}}]$. Thus for any $n>N_{w}$, $[g_{n},v_{1}]=\emph{id}$.
Using this relation we can switch $g_{n}^{a_{1}}$ and $v_{1}$ with each other and obtain a word $w'(g
_{n})=g_{n}^{a_{k}}v_{k}\ldots g_{n}^{a_{2}}v_{2}v_{1}g_{n}^{a_{1}}$. Since $v_{2}\in F\setminus H_{1}$,
clearly $v_{2}v_{1}\in F\setminus H_{1}$. Thus there is some $N\in\mathbb{N}$ such that for all $n
>N$, the word $w'=y^{a_{k}}v_{k}\ldots y^{a_{2}}v_{2}v_{1}y^{a_{1}}$ belongs to $W_{n}$. Thus it
follows from the construction of the sequence $(g_{n})_{n<\omega}$ that for all $n>N$, $w'(g_{n})
\neq\emph{id}$. Since $w'(g_{n})=w(g_{n})$, we see that $\lim _{n\to\infty}(G_{n})\models w(g)\neq
\emph{id}$.\parskip0pt

 Finally suppose that $v_{1}\in H_{1}$ and $w$ is of the form $w=y^{a}v_{1}$. Since we have
chosen the sequence $(g_{n})_{n<\omega}$ so that $supp(g_{n})\subseteq (1-\frac{1}{2^{N}},1)$,
we see that $supp(g_{n})\cap v_{1}=\emptyset$ for $n$ large enough. Thus $g_{n}^{a}v_{1}$ is
non-trivial for $n$ large enough. The proof is finished.\\

\ \ \ \ \ \ \ \ \ \ \ \ \ \ \ \ \ \ \ \ \ \ \ \ \ \ \ \ \ \ \ \ \ \ \ \ \ \ \ \ \ \ \ \ \ \ \ \ \ \ \ \ \ \ \ \ \ \ \ \ \ \ \ \ \ \ \ \ \ \ \ \ \ \ \ \ \ $\square$\\

\section{Limits of $F$}

 Brin and Squier have shown in $\cite{BS}$ that the Thompson's group $F$ does not satisfy
any group law (see also Corollary 1.2, \cite{A}). In this section we show how to construct words
with constants from $F$, which are equal to identity for any substitution in $F$. Then we apply
this to the main results of the section which find some restrictions on the form of $F$-limits.\\

\subsection{One-variable laws versus multi-variable laws}

\begin{defi}
	 Let $w(y_{1},\ldots , y_{t})$ be a non-trivial word over a group $G$, reduced in the
	group $\mathbb{F} _{t}\ast G$ and containing at least one variable and at least one
	constant. We will call $w$ a \emph{law with constants} in $G$ if for any $\bar{g}=(g
	_{1},\ldots , g_{t})\in G^{t}$, the value of $w(\bar{g})$ is equal to $1$.
\end{defi}

 The following theorem explains why we will be interested only in one-variable laws with
constants. It looks so natural that we would not be surprised if somebody has discovered it
before.

\begin{theo}\label{ng1}
	If a group $G$ satisfies a law with constants in $n$ variables, then $G$ also
	satisfies an one-variable law with constants.
\end{theo}

\emph{Proof.} \ Let $w(x_{1},\ldots , x_{n})$ be a law with constants of $G$. Choose $n$
words $w_{1}(x,y), w_{2}(x,y),\ldots , w_{n}(x,y)$ which freely generate a free subgroup of
the free group $\mathbb{F}(x,y)$. Let $w'(x,y)$ be the reduced word in $\mathbb{F} _{2}\ast
G$ representing $w(w_{1}(x,y),w_{2}(x,y),\ldots , w_{n}(x,y))$. It is clear that $w'(x,y)$ is
a law with constants.\parskip0pt

 Thus we may assume that the initial law with constants $w$ depends only on two
variables $x$ and $y$. Denote by $w(x,y)$ its reduced form in $\mathbb{F}_{2}\ast G$.
Hence we may rewrite $w(x,y)$ (using conjugation if necessary) in the following form
$$w(x,y)=u_{1}x^{\alpha _{1,1}}y^{\beta _{1,1}}\ldots x^{\alpha _{1,l_{1}}}y^{\beta _{1,l
	_{1}}}u_{2}\ldots u_{m}x^{\alpha _{m,1}}y^{\beta _{m,1}}\ldots x^{\alpha _{m,
	l_{m}}}y^{\beta _{m,l_{m}}},$$
where $m, l_{1},\ldots ,l_{m}\in\mathbb{N}\setminus\{ 0\}$, $\alpha _{i,j}$, $\beta _{i,j}\in
\mathbb{Z}$ for all $i\leq m$ and $j\leq l_{i}$, and $u_{1},\ldots u_{m}$ are non-trivial
constants in $G$. Wlog we may also assume that for any $i\leq m$, if $\alpha _{i,1}=0$
then $\beta _{i,1}\neq 0$, if $\beta _{i,l_{i}}=0$, then $\alpha _{i,l_{i}}\neq 0$ and $x
^{\alpha _{i,1}}y^{\beta _{i,1}}\ldots x^{\alpha _{i,l_{i}}}y^{\beta _{i,l_{i}}}$ is reduced in the
free group $\mathbb{F}(x,y)$. If there are no occurances of $x$ or of $y$ in $w(x,y)$,
then the proof is finished. Hence we also assume that there is some $i\leq m$ such that
$\alpha _{i,1}\neq 0$ and there is some $i'\leq m$ such that $\beta _{i',1}\neq
0$.\parskip0pt

 Now let $A=\{ v_{1}(x), v_{2}(x),\ldots ,v_{k}(x)\}$ be the set of all maximal subwords
of $w$, which contain $x$ and do not contain $y$. Suppose that the system of
inequalities $v_{1}(x)\neq 1,\ldots ,v_{k}(x)\neq 1$ has a solution $a\in G$. Then $w(a,
y)$ is an one-variable law with constants in $G$. Indeed,
$$w(a,y)=u_{1}a^{\alpha _{1,1}}y^{\beta _{1,1}}\ldots a^{\alpha _{1,l_{1}}}y^{\beta _{1,l
	_{1}}}u_{2}\ldots u_{n}a^{\alpha _{m,1}}y^{\beta _{m,1}}\ldots a^{\alpha _{m,
	l_{m}}}y^{\beta _{m,l_{m}}},$$
is reduced in the group $\mathbb{Z}\ast G$ and of course for any $b\in G$, $w(a,b)=
1$.\parskip0pt

 Now suppose that there is no solution of the system $v_{1}(x)\neq 1,\ldots ,v_{k}(x)\neq
1$ in $G$. If $k=1$, then $v_{1}(x)$ is an one-variable law with constants and the proof
is finished.\parskip0pt

 Suppose that $k>1$. We choose a sequence of elements $a_{2},\ldots , a_{k}\in G$ in
the following way. Consider the word
$$w_{2}(x,z):=[v_{1}(x),(v_{2}(x))^{z}]=(v_{1}(x))^{-1}z^{-1}(v_{2}(x))^{-1}zv_{1}(x)z^{-1}v
	_{2}(x)z.$$
 If there is any $a\in G$ such that $w_{2}(x,a)\notin G<\mathbb{Z}\ast G$, then we choose
$a_{2}:=a$. If $w_{2}(x,a)\in G$ for any $a\in G$, then presenting $v_{1}$ and $v_{2}$ in
the reduced forms $u'_{1}x^{\alpha _{1}}v'_{1}(x)x^{\beta _{1}}u''_{1}$ and $u'_{2}x
^{\alpha _{2}}v'_{2}(x)x^{\beta _{2}}u''_{2}$, respectively, with $u'_{1}, u''_{1}, u'_{2}, u''
_{2}\in\{ 1, u_{1},\ldots , u_{m}\}$, we rewrite $w_{2}(x,z)$ in the following form
$$w_{2}=(u''_{1})^{-1}\tilde{v}_{1}(x)^{-1}(u'_{1})^{-1}z^{-1}(u''_{2})^{-1}\tilde{v}_{2}(x)
	^{-1}(u'_{2})^{-1}zu'_{1}\tilde{v}_{1}(x)u''_{1}z^{-1}u'_{2}\tilde{v}_{2}(x)u''_{2}
	z,$$
where $\tilde{v}_{1}=x^{\alpha _{1}}v'_{1}(x)x^{\beta _{1}}$ and $\tilde{v}_{2}=x^{\alpha
_{2}}v'_{2}(x)x^{\beta _{2}}$. Since $w_{2}(x,a)\in G$ for any $a$, we see that for any $a
\in G$ at least one of the following equations holds:
$$(u'_{1})^{-1}a^{-1}(u''_{2})^{-1}=1,\ \ \ (u'_{2})^{-1}au'_{1}=1,\ \ \ u''_{1}a^{-1}u'_{2}=1 .$$
But this implies that $G$ has at most three elements and hence obviously satisfies an
one-variable law. Thus we may assume that we can always find $a_{2}$ such that $[v
_{1}(x),(v_{2}(x))^{a_{2}}]\notin G$.\parskip0pt

 Now suppose that we have chosen $a_{2},\ldots , a_{j}$ so that
$$w_{j}(x,a_{j}):=[[\ldots [[v_{1}(x),(v_{2}(x))^{a_{2}}],(v_{3}(x))^{a_{3}}],\ldots ], (v_{j}
	(x))^{a_{j}}]\neq G<\mathbb{Z}\ast G.$$
 Let
$$w_{j+1}(x,z)=[w_{j}(x,a_{j}),(v_{j+1}(x))^{z}] .$$
If there is any $a\in G$ such that $w_{j+1}(x,a)\notin G<\mathbb{Z}\ast G$, then we
choose $a_{j+1}:=a$. If there is no such an $a$, then similarly as above, we present $w
_{j}(x,a_{j})$ and $v_{j+1}(x)$ in the reduced form. Since $w_{j}(x,a_{j})\notin G$, we may
assume that its reduced form is equal to $u'_{j}x^{\alpha _{j}}v'_{j}(x)x^{\beta _{j}}u''_{j}$ for
some $\alpha _{j}, \beta _{j}\in\mathbb{Z}$, $u'_{j},u''_{j}\in G$, $v'_{j}(x)\in\mathbb{Z}\ast
G$, where $x^{\alpha _{j}}v'_{j}(x)x^{\beta _{j}}$ starts and finishes with some occurance
of $x$. Similarly as before, we present $v_{j+1}(x)$ in the reduced form $u'_{j+1}x^{\alpha
_{j+1}}v'_{j+1}(x)x^{\beta _{j+1}}u''_{j+1}$ with $u'_{j+1}, u''_{j+1}\in\{ 1, u_{1},\ldots , u
_{m}\}$. Since $w_{j+1}(x,a)$ belongs to $G$ for any $a$, we can extract the equations
$$(u'_{j})^{-1}a^{-1}(u''_{j+1})^{-1}=1,\ \ \ (u'_{j+1})^{-1}au'_{j}=1,\ \ \ u''_{j}a^{-1}u'_{j+1}=1$$
exactly as above and obtain an one-variable law of $G$.\parskip0pt

 We continue this procedure to obtain a word
$$\tilde{w}(x)=[[\ldots [[v_{1}(x),(v_{2}(x))^{a_{2}}], (v_{3}(x))^{a_{3}}],\ldots ],(v_{k}(x))^{a
	_{k}}]$$
such that $\tilde{w}(x)\notin G$. Since we have assumed that the system of inequalities
$v_{1}(x)\neq 1,\ldots ,v_{k}(x)\neq 1$ does not have a solution in $G$, we see that for
any $a\in G$ at least one of the words $v_{1}(x), v_{2}(x),\ldots ,v_{k}(x)$ is trivial. Thus
for any $a\in G$, $\tilde{w}(a)=1$ and $\tilde{w}(x)$ is an one-variable law with constants
in $G$.\\

\ \ \ \ \ \ \ \ \ \ \ \ \ \ \ \ \ \ \ \ \ \ \ \ \ \ \ \ \ \ \ \ \ \ \ \ \ \ \ \ \ \ \ \ \ \ \ \ \ \ \ \ \ \ \ \ \ \ \ \ \ \ \ \ \ \ \ \ \ \ \ \ \ \ \ \ \ $\square$\\

\subsection{Laws with constants for the Thompson's group $F$}

 The following proposition gives a construction of certain laws with constants in $F$.

\begin{prop} \label{lwc2}
	 Consider the standard action of Thompson's group $F$ on $[0,1]$. Suppose we are given
	two pairwise disjoint open dyadic subintervals $I_{i}=(p_{i}, q_{i})\subseteq [0,1]$, $i=1, 2$,
	and assume that $p_{1}<p_{2}$. Fix any non-trivial $h_{1}\in F_{\bar{I}_{1}}$ and $h_{2}\in
	F_{\bar{I}_{2}}$ and denote:
	$$w^{-}=[h_{1}^{y},h_{2}]=y^{-1}h_{1}^{-1}yh_{2}^{-1}y^{-1}h_{1}yh_{2}$$
	$$w^{+}=[h_{1}^{(y^{-1})},h_{2}]=yh_{1}^{-1}y^{-1}h_{2}^{-1}yh_{1}y^{-1}h_{2}$$
	Then the word $w:=[w^{-},w^{+}]$ is a law with constants in $F$.
\end{prop}

\setlength{\unitlength}{2mm}
\begin{picture}(45,45)
	\linethickness{0.05mm}
		\put(-5,9){\line(1,0){32}}
		\put(-5,41){\line(1,0){32}}
		\multiput(-5,9)(4,0){9}%
			{\line(0,1){32}}

			\put(-5,30){\line(1,0){4}}
			\put(-5,36){\line(1,0){4}}
			\put(-3.8,32.5){\footnotesize $h_{2}$}

			\put(-1.1,30){\line(1,-2){4}}
			\put(-1.05,30){\line(1,-2){4}}
			\put(-1,30){\line(1,-2){4}}
			\put(-0.95,30){\line(1,-2){4}}
			\put(-0.9,30){\line(1,-2){4}}
			\put(0.5,11){\footnotesize $g$}

			\put(3,14){\line(1,0){4}}
			\put(3,20){\line(1,0){4}}
			\multiput(3,21.95)(1,0){4}%
				{\line(1,0){0.5}}
			\multiput(3,22)(1,0){4}%
				{\line(1,0){0.5}}
			\multiput(3,22.05)(1,0){4}%
				{\line(1,0){0.5}}
			\put(4.2,16.5){\footnotesize $h_{1}$}

			\put(6.9,22){\line(1,2){4}}
			\put(6.95,22){\line(1,2){4}}
			\put(7,22){\line(1,2){4}}
			\put(7.05,22){\line(1,2){4}}
			\put(7.1,22){\line(1,2){4}}
			\put(8,11){\footnotesize $g^{-1}$}

			\put(11,30){\line(1,0){4}}
			\put(11,36){\line(1,0){4}}
			\put(12,32.5){\footnotesize $h_{2}^{-1}$}

			\put(14.9,30){\line(1,-2){4}}
			\put(14.95,30){\line(1,-2){4}}
			\put(15,30){\line(1,-2){4}}
			\put(15.05,30){\line(1,-2){4}}
			\put(15.1,30){\line(1,-2){4}}
			\put(16.5,11){\footnotesize $g$}

			\put(19,14){\line(1,0){4}}
			\put(19,20){\line(1,0){4}}
			\multiput(19,21.95)(1,0){4}%
				{\line(1,0){0.5}}
			\multiput(19,22)(1,0){4}%
				{\line(1,0){0.5}}
			\multiput(19,22.05)(1,0){4}%
				{\line(1,0){0.5}}
			\put(20,16.5){\footnotesize $h_{1}^{-1}$}

			\put(22.9,22){\line(1,2){4}}
			\put(22.95,22){\line(1,2){4}}
			\put(23,22){\line(1,2){4}}
			\put(23.05,22){\line(1,2){4}}
			\put(23.1,22){\line(1,2){4}}
			\put(24,11){\footnotesize $g^{-1}$}

	\put(-5,30){\circle*{0.5}}
	\put(-5,36){\circle*{0.5}}

	\put(-1,30){\circle*{0.5}}
	\put(-1,36){\circle*{0.5}}

	\put(3,14){\circle*{0.5}}
	\put(3,20){\circle*{0.5}}
	\put(3,22){\circle*{0.5}}

	\put(7,14){\circle*{0.5}}
	\put(7,20){\circle*{0.5}}
	\put(7,22){\circle*{0.5}}

	\put(11,30){\circle*{0.5}}
	\put(11,36){\circle*{0.5}}

	\put(15,30){\circle*{0.5}}
	\put(15,36){\circle*{0.5}}

	\put(19,14){\circle*{0.5}}
	\put(19,20){\circle*{0.5}}
	\put(19,22){\circle*{0.5}}

	\put(23,14){\circle*{0.5}}
	\put(23,20){\circle*{0.5}}
	\put(23,22){\circle*{0.5}}

	\put(27,30){\circle*{0.5}}

	\linethickness{0.125mm}

			\multiput(-5,14)(1,0){34}%
				{\line(1,0){0.5}}

			\multiput(-5,20)(1,0){34}%
				{\line(1,0){0.5}}

			\multiput(-5,30)(1,0){34}%
				{\line(1,0){0.5}}

			\multiput(-5,36)(1,0){34}%
				{\line(1,0){0.5}}

			\put(30,13.5){\footnotesize $p_{1}$}
			\put(30,19.5){\footnotesize $q_{1}$}
			\put(30,29.5){\footnotesize $p_{2}$}
			\put(30,35.5){\footnotesize $q_{2}$}

	\linethickness{0.75mm}

			\put(-5,30){\line(1,0){4}}
			\put(-5,36){\line(1,0){4}}
			\put(-5,30){\line(0,1){6}}
			\put(-1,30){\line(0,1){6}}

			\put(11,30){\line(1,0){4}}
			\put(11,36){\line(1,0){4}}
			\put(11,30){\line(0,1){6}}
			\put(15,30){\line(0,1){6}}

			\put(3,14){\line(1,0){4}}
			\put(3,20){\line(1,0){4}}
			\put(3,14){\line(0,1){6}}
			\put(7,14){\line(0,1){6}}

			\put(19,14){\line(1,0){4}}
			\put(19,20){\line(1,0){4}}
			\put(19,14){\line(0,1){6}}
			\put(23,14){\line(0,1){6}}

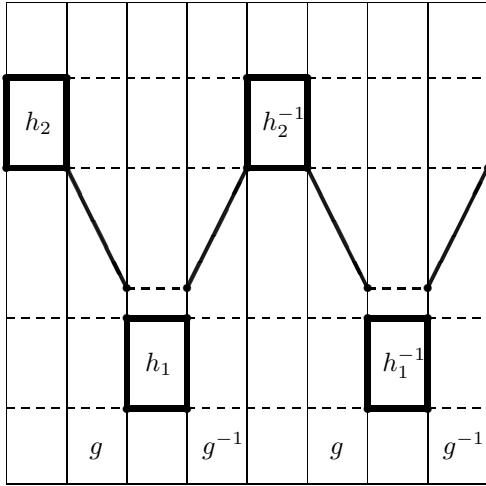
	\put(-5,5.5){{\footnotesize Picture 12.\ \ $w^{-}(g)=\emph{id}$, where $g(p_{2})\geq q_{1}$}}

	\linethickness{0.05mm}
		\put(35,9){\line(1,0){32}}
		\put(35,41){\line(1,0){32}}
		\multiput(35,9)(4,0){9}%
			{\line(0,1){32}}

			\multiput(35,19.95)(1,0){4}%
				{\line(1,0){0.5}}
			\multiput(35,20)(1,0){4}%
				{\line(1,0){0.5}}
			\multiput(35,20.05)(1,0){4}%
				{\line(1,0){0.5}}
			\put(35,30){\line(1,0){4}}
			\put(35,36){\line(1,0){4}}
			\put(36.2,32.5){\footnotesize $h_{2}$}

			\put(38.92,20){\line(1,3){4}}
			\put(38.96,20){\line(1,3){4}}
			\put(39,20){\line(1,3){4}}
			\put(39.04,20){\line(1,3){4}}
			\put(39.08,20){\line(1,3){4}}
			\put(40,11){\footnotesize $g^{-1}$}

			\put(43,14){\line(1,0){4}}
			\put(43,20){\line(1,0){4}}
			\multiput(43,31.95)(1,0){4}%
				{\line(1,0){0.5}}
			\multiput(43,32)(1,0){4}%
				{\line(1,0){0.5}}
			\multiput(43,32.05)(1,0){4}%
				{\line(1,0){0.5}}
			\put(44.2,16.5){\footnotesize $h_{1}$}

			\put(47,31.92){\line(1,-3){4}}
			\put(47,31.94){\line(1,-3){4}}
			\put(47,32){\line(1,-3){4}}
			\put(47,32.04){\line(1,-3){4}}
			\put(47,32.08){\line(1,-3){4}}
			\put(48.5,11){\footnotesize $g$}

			\multiput(51,19.95)(1,0){4}%
				{\line(1,0){0.5}}
			\multiput(51,20)(1,0){4}%
				{\line(1,0){0.5}}
			\multiput(51,20.05)(1,0){4}%
				{\line(1,0){0.5}}
			\put(51,30){\line(1,0){4}}
			\put(51,36){\line(1,0){4}}
			\put(52,32.5){\footnotesize $h_{2}^{-1}$}

			\put(54.92,20){\line(1,3){4}}
			\put(54.96,20){\line(1,3){4}}
			\put(55,20){\line(1,3){4}}
			\put(55.04,20){\line(1,3){4}}
			\put(55.08,20){\line(1,3){4}}
			\put(56,11){\footnotesize $g^{-1}$}

			\put(59,14){\line(1,0){4}}
			\put(59,20){\line(1,0){4}}
			\multiput(59,31.95)(1,0){4}%
				{\line(1,0){0.5}}
			\multiput(59,32)(1,0){4}%
				{\line(1,0){0.5}}
			\multiput(59,32.05)(1,0){4}%
				{\line(1,0){0.5}}
			\put(60,16.5){\footnotesize $h_{1}^{-1}$}

			\put(62.92,32){\line(1,-3){4}}
			\put(62.96,32){\line(1,-3){4}}
			\put(63,32){\line(1,-3){4}}
			\put(63.04,32){\line(1,-3){4}}
			\put(63.08,32){\line(1,-3){4}}
			\put(64.5,11){\footnotesize $g$}

	\put(35,30){\circle*{0.5}}
	\put(35,36){\circle*{0.5}}

	\put(39,20){\circle*{0.5}}
	\put(39,30){\circle*{0.5}}
	\put(39,36){\circle*{0.5}}

	\put(43,14){\circle*{0.5}}
	\put(43,20){\circle*{0.5}}
	\put(43,32){\circle*{0.5}}

	\put(47,14){\circle*{0.5}}
	\put(47,20){\circle*{0.5}}
	\put(47,32){\circle*{0.5}}

	\put(51,20){\circle*{0.5}}
	\put(51,30){\circle*{0.5}}
	\put(51,36){\circle*{0.5}}

	\put(55,20){\circle*{0.5}}
	\put(55,30){\circle*{0.5}}
	\put(55,36){\circle*{0.5}}

	\put(59,14){\circle*{0.5}}
	\put(59,20){\circle*{0.5}}
	\put(59,32){\circle*{0.5}}

	\put(63,14){\circle*{0.5}}
	\put(63,20){\circle*{0.5}}
	\put(63,32){\circle*{0.5}}

	\put(67,20){\circle*{0.5}}

	\linethickness{0.125mm}

			\multiput(33,14)(1,0){34}%
				{\line(1,0){0.5}}

			\multiput(33,20)(1,0){34}%
				{\line(1,0){0.5}}

			\multiput(33,30)(1,0){34}%
				{\line(1,0){0.5}}

			\multiput(33,36)(1,0){34}%
				{\line(1,0){0.5}}

	\linethickness{0.75mm}

			\put(35,30){\line(1,0){4}}
			\put(35,36){\line(1,0){4}}
			\put(35,30){\line(0,1){6}}
			\put(39,30){\line(0,1){6}}

			\put(51,30){\line(1,0){4}}
			\put(51,36){\line(1,0){4}}
			\put(51,30){\line(0,1){6}}
			\put(55,30){\line(0,1){6}}

			\put(43,14){\line(1,0){4}}
			\put(43,20){\line(1,0){4}}
			\put(43,14){\line(0,1){6}}
			\put(47,14){\line(0,1){6}}

			\put(59,14){\line(1,0){4}}
			\put(59,20){\line(1,0){4}}
			\put(59,14){\line(0,1){6}}
			\put(63,14){\line(0,1){6}}

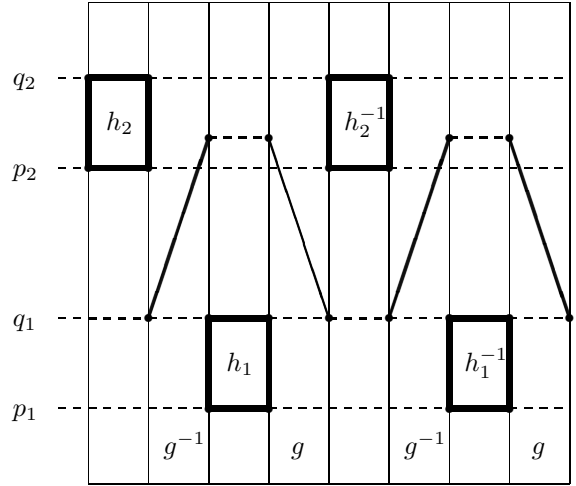
	\put(34.5,5){{\footnotesize Picture 13.\ \ $w^{+}(g)=\emph{id}$, where $g(p_{2})< q_{1}$}}

\end{picture}

\emph{Proof.} \ It is easy to see that $w$ cannot be reduced to a constant. We claim that
\begin{quote} for any $g\in F$ satisfying $w^{-}(g)\neq\emph{id}$ the word $w^{+}(g)$ is equal to
	identity.\end{quote}

 It follows from Remark \ref{GS} that if the supports of any two given elements from $F$ do not intersect,
then the commutator of these two elements is trivial. Hence we see that $w^{-}(g)=[h_{1}^{g},h_{2}]
\neq\emph{id}$ implies $supp(h_{1}^{g})\cap supp(h_{2})\neq\emptyset$. By Remark \ref{rs}, this implies
$g^{-1}(supp(h_{1}))\cap supp(h_{2})\neq\emptyset$. It follows that $g(p_{2})<q_{1}$. Since any $g\in F$
is an increasing function, we obtain the following sequence of implications
$$g(p_{2})<q_{1}\ \Rightarrow\ p_{2}<g^{-1}(q_{1})\ \Rightarrow\ p_{2}<g^{-1}(p_{2})\ \Rightarrow\ q_{1}
	<g^{-1}(p_{2}).$$
This shows that $supp(h_{1})\cap g^{-1}(supp(h_{2}))=\emptyset$. Hence using Remark \ref{rs} and
Remark \ref{GS} similarly as above we see that $w^{+}(g)=\emph{id}$. This proves that for any $g\in F$
either $w^{-}(g)=\emph{id}$ or $w^{+}(g)=\emph{id}$.\\

\ \ \ \ \ \ \ \ \ \ \ \ \ \ \ \ \ \ \ \ \ \ \ \ \ \ \ \ \ \ \ \ \ \ \ \ \ \ \ \ \ \ \ \ \ \ \ \ \ \ \ \ \ \ \ \ \ \ \ \ \ \ \ \ \ \ \ \ \ \ \ \ \ \ \ \ \ $\square$\\

 We have also another version of Proposition \ref{lwc2}, which employs different construction of a law
with constants, which is not based on the conjugation by the inverse of the variable.

\begin{prop} \label{lwc4}
	 Consider the standard action of Thompson's group $F$ on $[0,1]$. Suppose we are given
	four pairwise disjoint closed dyadic subintervals $I_{i}=[p_{i}, q_{i}]\subseteq [0,1]$, $1\leq i
	\leq 4$, and assume that $p_{1}<p_{2}<p_{3}<p_{4}$. Then for any non-trivial $h_{1}\in F
	_{I_{1}}$, $h_{2}\in F_{I_{2}}$, $h_{3}\in F_{I_{3}}$ and $h_{4}\in F_{I_{4}}$, let
	$$w_{14}=[h_{1}^{y},h_{4}]=y^{-1}h_{1}^{-1}yh_{4}^{-1}y^{-1}h_{1}yh_{4}$$
	$$w_{23}=[h_{2}^{y},h_{3}]=y^{-1}h_{2}^{-1}yh_{3}^{-1}y^{-1}h_{2}yh_{3}$$
	Then the word $w$ obtained from $[w_{14},w_{23}]$ by reduction in $\mathbb{Z}\ast F$
	(we treat the variable $y$ as a generator of $\mathbb{Z}$) is a law with constants in $F$.
\end{prop}

\setlength{\unitlength}{2mm}
\begin{picture}(45,45)
	\linethickness{0.05mm}
		\put(-5,9){\line(1,0){32}}
		\put(-5,41){\line(1,0){32}}
		\multiput(-5,9)(4,0){9}%
			{\line(0,1){32}}

			\put(-5,35){\line(1,0){4}}
			\put(-5,39){\line(1,0){4}}
			\put(-3.8,36.5){\footnotesize $h_{4}$}

			\put(-1.08,35){\line(1,-4){4}}
			\put(-1.04,35){\line(1,-4){4}}
			\put(-1,35){\line(1,-4){4}}
			\put(-0.96,35){\line(1,-4){4}}
			\put(-0.92,35){\line(1,-4){4}}
			\put(0.5,16){\footnotesize $g$}

			\put(3,11){\line(1,0){4}}
			\put(3,15){\line(1,0){4}}
			\multiput(3,18.95)(1,0){4}%
				{\line(1,0){0.5}}
			\multiput(3,19)(1,0){4}%
				{\line(1,0){0.5}}
			\multiput(3,19.05)(1,0){4}%
				{\line(1,0){0.5}}
			\put(4.2,12.5){\footnotesize $h_{1}$}

			\put(6.92,19){\line(1,4){4}}
			\put(6.96,19){\line(1,4){4}}
			\put(7,19){\line(1,4){4}}
			\put(7.04,19){\line(1,4){4}}
			\put(7.08,19){\line(1,4){4}}
			\put(8,16){\footnotesize $g^{-1}$}

			\put(11,35){\line(1,0){4}}
			\put(11,39){\line(1,0){4}}
			\put(12,36.5){\footnotesize $h_{4}^{-1}$}

			\put(14.92,35){\line(1,-4){4}}
			\put(14.96,35){\line(1,-4){4}}
			\put(15,35){\line(1,-4){4}}
			\put(15.04,35){\line(1,-4){4}}
			\put(15.08,35){\line(1,-4){4}}
			\put(16.5,16){\footnotesize $g$}

			\put(19,11){\line(1,0){4}}
			\put(19,15){\line(1,0){4}}
			\multiput(19,18.95)(1,0){4}%
				{\line(1,0){0.5}}
			\multiput(19,19)(1,0){4}%
				{\line(1,0){0.5}}
			\multiput(19,19.05)(1,0){4}%
				{\line(1,0){0.5}}
			\put(20,12.5){\footnotesize $h_{1}^{-1}$}

			\put(22.92,19){\line(1,4){4}}
			\put(22.96,19){\line(1,4){4}}
			\put(23,19){\line(1,4){4}}
			\put(23.04,19){\line(1,4){4}}
			\put(23.08,19){\line(1,4){4}}
			\put(24,16){\footnotesize $g^{-1}$}

	\put(-5,35){\circle*{0.5}}
	\put(-5,39){\circle*{0.5}}

	\put(-1,35){\circle*{0.5}}
	\put(-1,39){\circle*{0.5}}

	\put(3,11){\circle*{0.5}}
	\put(3,15){\circle*{0.5}}
	\put(3,19){\circle*{0.5}}

	\put(7,11){\circle*{0.5}}
	\put(7,15){\circle*{0.5}}
	\put(7,19){\circle*{0.5}}

	\put(11,35){\circle*{0.5}}
	\put(11,39){\circle*{0.5}}

	\put(15,35){\circle*{0.5}}
	\put(15,39){\circle*{0.5}}

	\put(19,11){\circle*{0.5}}
	\put(19,15){\circle*{0.5}}
	\put(19,19){\circle*{0.5}}

	\put(23,11){\circle*{0.5}}
	\put(23,15){\circle*{0.5}}
	\put(23,19){\circle*{0.5}}

	\put(27,35){\circle*{0.5}}

	\linethickness{0.125mm}

			\multiput(-5,11)(1,0){34}%
				{\line(1,0){0.5}}

			\multiput(-5,15)(1,0){34}%
				{\line(1,0){0.5}}

			\multiput(-5,20)(1,0){34}%
				{\line(1,0){0.5}}

			\multiput(-5,24)(1,0){34}%
				{\line(1,0){0.5}}

			\multiput(-5,26)(1,0){34}%
				{\line(1,0){0.5}}

			\multiput(-5,30)(1,0){34}%
				{\line(1,0){0.5}}

			\multiput(-5,35)(1,0){34}%
				{\line(1,0){0.5}}

			\multiput(-5,39)(1,0){34}%
				{\line(1,0){0.5}}

			\put(30,10.5){\footnotesize $p_{1}$}
			\put(30,14.5){\footnotesize $q_{1}$}
			\put(30,19.5){\footnotesize $p_{2}$}
			\put(30,23.5){\footnotesize $q_{2}$}
			\put(30,25.5){\footnotesize $p_{3}$}
			\put(30,29.5){\footnotesize $q_{3}$}
			\put(30,34.5){\footnotesize $p_{4}$}
			\put(30,38.5){\footnotesize $q_{4}$}

	\linethickness{0.75mm}

			\put(-5,35){\line(1,0){4}}
			\put(-5,39){\line(1,0){4}}
			\put(-5,35){\line(0,1){4}}
			\put(-1,35){\line(0,1){4}}

			\put(11,35){\line(1,0){4}}
			\put(11,39){\line(1,0){4}}
			\put(11,35){\line(0,1){4}}
			\put(15,35){\line(0,1){4}}

			\put(3,11){\line(1,0){4}}
			\put(3,15){\line(1,0){4}}
			\put(3,11){\line(0,1){4}}
			\put(7,11){\line(0,1){4}}

			\put(19,11){\line(1,0){4}}
			\put(19,15){\line(1,0){4}}
			\put(19,11){\line(0,1){4}}
			\put(23,11){\line(0,1){4}}

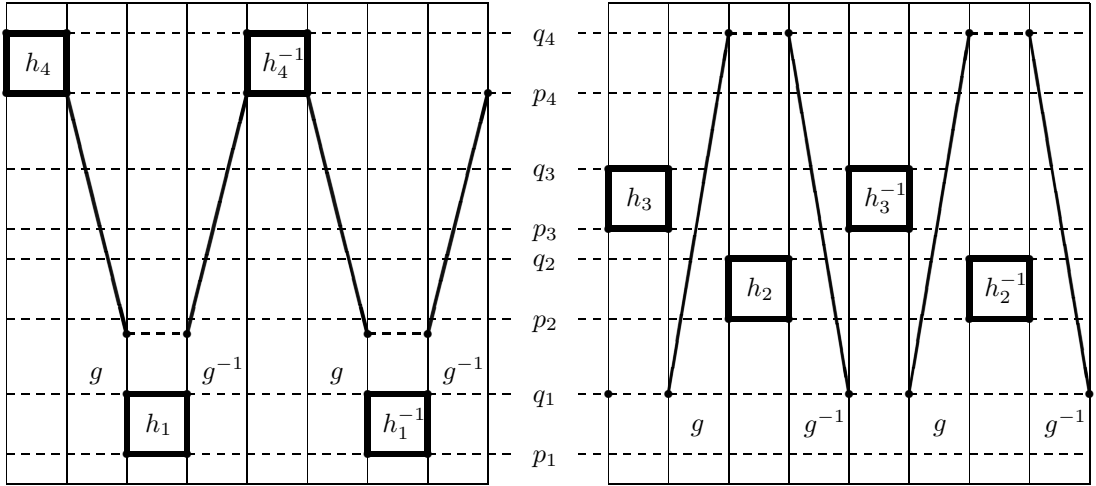
	\put(-5,5.5){{\footnotesize Picture 14.\ \ $w_{14}(g)=\emph{id}$, where $g(p_{4})> q_{1}$}}

	\linethickness{0.05mm}
		\put(35,9){\line(1,0){32}}
		\put(35,41){\line(1,0){32}}
		\multiput(35,9)(4,0){9}%
			{\line(0,1){32}}

			\multiput(35,14.95)(1,0){4}%
				{\line(1,0){0.5}}
			\multiput(35,15)(1,0){4}%
				{\line(1,0){0.5}}
			\multiput(35,15.05)(1,0){4}%
				{\line(1,0){0.5}}
			\put(35,26){\line(1,0){4}}
			\put(35,30){\line(1,0){4}}
			\put(36.2,27.5){\footnotesize $h_{3}$}

			\put(38.94,15){\line(1,6){4}}
			\put(38.97,15){\line(1,6){4}}
			\put(39,15){\line(1,6){4}}
			\put(39.03,15){\line(1,6){4}}
			\put(39.06,15){\line(1,6){4}}
			\put(40.5,12.5){\footnotesize $g$}

			\put(43,20){\line(1,0){4}}
			\put(43,24){\line(1,0){4}}
			\multiput(43,38.95)(1,0){4}%
				{\line(1,0){0.5}}
			\multiput(43,39)(1,0){4}%
				{\line(1,0){0.5}}
			\multiput(43,39.05)(1,0){4}%
				{\line(1,0){0.5}}
			\put(44.2,21.5){\footnotesize $h_{2}$}

			\put(46.94,39){\line(1,-6){4}}
			\put(46.97,39){\line(1,-6){4}}
			\put(47,39){\line(1,-6){4}}
			\put(47.03,39){\line(1,-6){4}}
			\put(47.06,39){\line(1,-6){4}}
			\put(48,12.5){\footnotesize $g^{-1}$}

			\multiput(51,14.95)(1,0){4}%
				{\line(1,0){0.5}}
			\multiput(51,15)(1,0){4}%
				{\line(1,0){0.5}}
			\multiput(51,15.05)(1,0){4}%
				{\line(1,0){0.5}}
			\put(51,26){\line(1,0){4}}
			\put(51,30){\line(1,0){4}}
			\put(52,27.5){\footnotesize $h_{3}^{-1}$}

			\put(54.94,15){\line(1,6){4}}
			\put(54.97,15){\line(1,6){4}}
			\put(55,15){\line(1,6){4}}
			\put(55.03,15){\line(1,6){4}}
			\put(55.06,15){\line(1,6){4}}
			\put(56.6,12.5){\footnotesize $g$}

			\put(59,20){\line(1,0){4}}
			\put(59,24){\line(1,0){4}}
			\multiput(59,38.95)(1,0){4}%
				{\line(1,0){0.5}}
			\multiput(59,39)(1,0){4}%
				{\line(1,0){0.5}}
			\multiput(59,39.05)(1,0){4}%
				{\line(1,0){0.5}}
			\put(60,21.5){\footnotesize $h_{2}^{-1}$}

			\put(62.94,39){\line(1,-6){4}}
			\put(62.97,39){\line(1,-6){4}}
			\put(63,39){\line(1,-6){4}}
			\put(63.03,39){\line(1,-6){4}}
			\put(63.06,39){\line(1,-6){4}}
			\put(64,12.5){\footnotesize $g^{-1}$}

	\put(35,15){\circle*{0.5}}
	\put(35,26){\circle*{0.5}}
	\put(35,30){\circle*{0.5}}

	\put(39,15){\circle*{0.5}}
	\put(39,26){\circle*{0.5}}
	\put(39,30){\circle*{0.5}}

	\put(43,20){\circle*{0.5}}
	\put(43,24){\circle*{0.5}}
	\put(43,39){\circle*{0.5}}

	\put(47,20){\circle*{0.5}}
	\put(47,24){\circle*{0.5}}
	\put(47,39){\circle*{0.5}}

	\put(51,15){\circle*{0.5}}
	\put(51,26){\circle*{0.5}}
	\put(51,30){\circle*{0.5}}

	\put(55,15){\circle*{0.5}}
	\put(55,26){\circle*{0.5}}
	\put(55,30){\circle*{0.5}}

	\put(59,20){\circle*{0.5}}
	\put(59,24){\circle*{0.5}}
	\put(59,39){\circle*{0.5}}

	\put(63,20){\circle*{0.5}}
	\put(63,24){\circle*{0.5}}
	\put(63,39){\circle*{0.5}}

	\put(67,15){\circle*{0.5}}

	\linethickness{0.125mm}

			\multiput(33,11)(1,0){34}%
				{\line(1,0){0.5}}

			\multiput(33,15)(1,0){34}%
				{\line(1,0){0.5}}

			\multiput(33,20)(1,0){34}%
				{\line(1,0){0.5}}

			\multiput(33,24)(1,0){34}%
				{\line(1,0){0.5}}

			\multiput(33,26)(1,0){34}%
				{\line(1,0){0.5}}

			\multiput(33,30)(1,0){34}%
				{\line(1,0){0.5}}

			\multiput(33,35)(1,0){34}%
				{\line(1,0){0.5}}

			\multiput(33,39)(1,0){34}%
				{\line(1,0){0.5}}

	\linethickness{0.75mm}

			\put(35,26){\line(1,0){4}}
			\put(35,30){\line(1,0){4}}
			\put(35,26){\line(0,1){4}}
			\put(39,26){\line(0,1){4}}

			\put(51,26){\line(1,0){4}}
			\put(51,30){\line(1,0){4}}
			\put(51,26){\line(0,1){4}}
			\put(55,26){\line(0,1){4}}

			\put(43,20){\line(1,0){4}}
			\put(43,24){\line(1,0){4}}
			\put(43,20){\line(0,1){4}}
			\put(47,20){\line(0,1){4}}

			\put(59,20){\line(1,0){4}}
			\put(59,24){\line(1,0){4}}
			\put(59,20){\line(0,1){4}}
			\put(63,20){\line(0,1){4}}

	\put(34.5,5){{\footnotesize Picture 15.\ \ $w_{23}(g)=\emph{id}$, where $g(q_{1})\geq p_{4}$}}

\end{picture}

\emph{Proof.} \ It is easy to see that $w$ cannot be reduced to a constant. \parskip0pt

 We claim that
\begin{quote} for any any $g\in F$ satisfying $w_{14}(g)\neq\emph{id}$ the word $w_{23}(g)$ is equal
	to identity \end{quote}
 We repeat the argument of the proof of Lemma \ref{lwc2}. It follows from Remark \ref{GS} that if the
supports of given two elements from $F$ do not intersect, then the commutator of these two
elements is trivial. Hence we see that $w^{14}(g)=[h_{1}^{g},h_{4}]\neq\emph{id}$ implies $supp(h_{1}
^{g})\cap supp(h_{4})\neq\emptyset$. By Remark \ref{rs}, this implies $g^{-1}(supp(h_{1}))\cap supp(h
_{4})\neq\emptyset$. It follows that $g(p_{4})<q_{1}$. Since any $g\in F$ is an increasing function, we
obtain the following sequence of implications
$$g(p_{4})<q_{1}\ \Rightarrow\ p_{4}<g^{-1}(q_{1})\ \Rightarrow\ p_{4}<g^{-1}(p_{2})\ \Rightarrow\ q_{3}
	<g^{-1}(p_{2}).$$

This shows that $supp(h_{3})\cap g^{-1}(supp(h_{2}))=\emptyset$. Hence using Remark \ref{rs} and
Remark \ref{GS} we see similarly as above that $w_{23}(g)=\emph{id}$. This proves that for any $g\in F$
either $w_{14}(g)=\emph{id}$ or $w_{23}(g)=\emph{id}$.\\

\ \ \ \ \ \ \ \ \ \ \ \ \ \ \ \ \ \ \ \ \ \ \ \ \ \ \ \ \ \ \ \ \ \ \ \ \ \ \ \ \ \ \ \ \ \ \ \ \ \ \ \ \ \ \ \ \ \ \ \ \ \ \ \ \ \ \ \ \ \ \ \ \ \ \ \ \ $\square$\\

\subsection{Free products}

 We now apply the construction from Proposition \ref{lwc2} to limits of Thompson's group $F$.

\begin{theo}\label{free}
	 Suppose we are given a convergent sequence of marked groups $((G_{n}, (x_{0}, x_{1},
	g_{n,1},\ldots , g_{n,s})))_{n<\omega}$, where $G_{n}=F$, $(g_{n,1},\ldots , g_{n,s})\in F$,
	$n\in\mathbb{N}$, and denote by $\mathbb{G}$ its limit. Then $\mathbb{G}\neq F\ast G$
	for any non-trivial $G$.
\end{theo}

 Before the proof we formulate a general statement, which exposes the main point of our argument.

\begin{prop}\label{pro}
	 Let $H=\langle h_{1},\ldots , h_{m}\rangle$ be torsion-free group, which satisfies a law with
	or without constants. Let $\mathbb{G}$ be the limit of a convergent sequence of marked
	groups $((G_{n}, (h_{1},\ldots , h_{m}, g_{n,1},\ldots , g_{n,t})))_{n<\omega}$, where $(g_{n,
	1},\ldots , g_{n,t})\in H$, $G_{n}=H$, $n\in\mathbb{N}$. Then $\mathbb{G}\neq H\ast K$ for
	any non-trivial $K<\mathbb{G}$.
\end{prop}

\emph{Proof.} \ It is clear that $\mathbb{G}$ is torsion-free. To obtain a contradiction suppose that
$\mathbb{G}=H\ast K$, $K\neq\{ 1\}$, and $\mathbb{G}$ is marked by a tuple $(h_{1},\ldots , h_{m},
f_{1},\ldots , f_{t})$, where $h_{i}$ are as in the formulation. Suppose that $H$ satisfies a law without
constants and denote it by $w(y_{1},\ldots ,y_{l})$. Similarly as in the proof of Theorem \ref{ng1}, we
may choose $l$ words $w_{1}(x,y),\ldots , w_{l}(x,y)$ which generate a free subgroup of rank $l$ in
the free group $\mathbb{F}(x,y)$. Now $w'(x,y)=w(w_{1}(x,y), w_{2}(x,y),\ldots , w_{n}(x,y))$ is a law
without constants in two variables. Thus we assume that $l=2$ and the initial word $w$ is a two-variable
law without constants. Moreover, since $H$ is torsion-free, we may assume that $w$ is not an one
variable law.\parskip0pt

 Let $g=u(\bar{h},\bar{f})$ be an element of $K\setminus\{ 1\}$ and $h'\in H\setminus\{ 1\}$. Obviously
$w(h',u(\bar{h}, g_{n,1},\ldots , g_{n,t}))=1_{H}$ for all $n<\omega$. It follows from the definition of an
$H$-limit group that $w(h',u(\bar{h},\bar{f}))=1_{\mathbb{G}}$. Since $w$ is not an one-variable law, the
reduced form of $w(h',y)\in H\ast\mathbb{Z}$ has non-trivial occurances of $y$ and elements of $H$.
As $|g|=\infty$, the equality $w(h',g)=1$ contradicts the fact that $\mathbb{G}$ is the free product of
$H$ and $K$.\parskip0pt

 Now assume that $H$ does not satisfy any law without constants and let $w(y_{1},\ldots , y_{l})$ be
a law with constants in $H$. It follows from Theorem \ref{ng1} that $G$ satisfies also an one-variable
law with constants. Hence assume that $w(y)$ is a law with constants in $H$. Similarly as above let
$g=u(\bar{h},\bar{f})$ be an element of $K\setminus\{ 1\}$. Obviously $w(u(\bar{h}, g_{n,1},\ldots , g
_{n,t}))=1_{H}$ for all $n<\omega$. Once again it follows from the definition of an $H$-limit group that
$w(u(\bar{h},\bar{f}))=1_{\mathbb{G}}$. Since $w$ was chosen to be non-trivial, with constants from
$H$ and $|g|=\infty$, we obtain again a contradiction with the fact that $\mathbb{G}$ is the free
product of $H$ and $K$.\\

\ \ \ \ \ \ \ \ \ \ \ \ \ \ \ \ \ \ \ \ \ \ \ \ \ \ \ \ \ \ \ \ \ \ \ \ \ \ \ \ \ \ \ \ \ \ \ \ \ \ \ \ \ \ \ \ \ \ \ \ \ \ \ \ \ \ \ \ \ \ \ \ \ \ \ \ \ $\square$\\

\emph{Proof of Theorem \ref{free}.} \ It follows directly from Theorem \ref{ng1}, that there is some
word $w(y)$, which is a law with constants in $F$, and hence we just apply Proposition \ref{pro} for $H
=F$, $h_{1}=x_{0}$ and $h_{2}=x_{1}$. \\

\ \ \ \ \ \ \ \ \ \ \ \ \ \ \ \ \ \ \ \ \ \ \ \ \ \ \ \ \ \ \ \ \ \ \ \ \ \ \ \ \ \ \ \ \ \ \ \ \ \ \ \ \ \ \ \ \ \ \ \ \ \ \ \ \ \ \ \ \ \ \ \ \ \ \ \ \ $\square$\\

\subsection{HNN-extensions}

 Now we proceed to discuss the case of HNN-extensions. The following theorem is one of the main
results of Section 4.

\begin{theo}\label{mt}
	Let $(G_{n})_{n<\omega}$ be a convergent sequence of groups, where $G_{n} = (F, (g_{n},
	x_{0}, x_{1}))$, $g_{n}\in F$ for $n<\omega$, and let $(G,(g,x_{0}, x_{1}))$, be its limit. Then
	$G$ is not an HNN-extension of F of the form $F\ast _{\alpha}=\langle F, g\ |\ ghg^{-1} =
	\alpha (h)\ h\in H\rangle$ where $H$ is some finitely generated subgroup of $F$ and $\alpha$
	is some embedding of $H$ into $F$.
\end{theo}

 Theorem \ref{h1} shows that Theorem \ref{mt} cannot be generalized to the case of infinitely generated
subgroups. However, it turns out that in the case of centralized HNN-extensions the condition saying that
$H$ is finitely generated, can be replaced by a technical assumption of the following form:
$$(\diamondsuit )\ \exists h',h''\in H\ \bigg( \{ 0,1\}\subseteq\overline{supp}(h')\cup\overline{supp}(h'')\bigg) .$$
 Therefore we prove Theorem \ref{mt} in two steps. First we use $(\diamondsuit )$ and prove the
centralized case.

\begin{theo}\label{dmt}
	Let $(G_{n})_{n<\omega}$ be a convergent sequence of groups, where $G_{n} = (F, (g_{n},
	x_{0}, x_{1}))$, $g_{n}\in F$ for $n<\omega$, and let $(G,(g,x_{0}, x_{1}))$, be its limit. Then
	$G$ is not an HNN-extension of F of the form $F\ast _{H}=\langle F, g\ |\ [g,h]=\emph{id}, h\in H
	\rangle$ where $H$ is some subgroup of $F$ satisfying $(\diamondsuit )$.
\end{theo}

\emph{Proof of Theorem \ref{dmt}.}\ To obtain a contradiction assume that $G=F\ast _{H}$.\parskip0pt

 We start with the case $H=F$. Since $G=F\ast _{H}$, for any $h\in H$, $G\models [g,h]=\emph{id}$. This implies
that for any $h\in H$ there exists $N_{h}\in\mathbb{N}$ such that for all $n>N_{h}$, $G_{n}\models [g_{n},h]
=\emph{id}$. Fix two elements $h_{1}, h_{2}\in F$ such that $\overline{supp}(h_{1})=\overline{supp}(h_{2})=[0,
1]$, both do not have dividing points and they do not have a common root. It follows from Remark \ref{GS}
that $C_{F}(h_{1})=\langle\hat{h}_{1}\rangle$ and $C_{F}(h_{2})=\langle\hat{h}_{2}\rangle$ for some roots
$\hat{h}_{1}$ and $\hat{h}_{2}$ of $h_{1}$ and $h_{2}$ respectively. Since for almost all $n$'s ($n>\max\{ N
_{h_{1}}, N_{h_{2}}\}$), $g_{n}\in\langle\hat{h}_{1}\rangle\cap\langle\hat{h}_{2}\rangle$, we see that for almost
all $n$'s, $g_{n}=\emph{id}$. Thus $G=F\ncong F\ast _{H}$, a contradiction.\parskip0pt

 Let $H\neq F$ and let $h_{1}, h_{2}, \ldots$ be an enumeration of $H$. For any $i\geq 1$ denote by $h
_{i,1},\ldots h_{i,l_{i}}$, $l_{i}\in\mathbb{N}\setminus\{ 0\}$, the elements of the defragmentation of $h_{i}$.\\

 \textbf{Claim.} There is a dyadic non-trivial interval $[p,q]\subseteq [0,1]$ such that $F_{[p,q]}\cap\langle h_{i,j}\
|\ i\geq 1, j\leq l_{i}\rangle$ is trivial or cyclic.\\

\emph{Proof of Claim.}\ Assume the contrary, i.e.:
$$(\dagger )\ \ \forall [p,q]\subseteq [0,1], p,q\in\mathbb{Z}[\frac{1}{2}]\ \Big(\Big( F_{[p,q]}\cap\Big\langle h_{i,j}\
	\Big|\ i\geq 1, j\leq l_{i}\Big\rangle\neq\Big\{ \emph{id}\Big\}\Big)\ \ \land$$
$$\Big( F_{[p,q]}\cap\Big\langle h_{i,j}\ \Big|\ i\geq 1, j\leq l_{i}\Big\rangle\ncong\mathbb{Z}\Big)\Big) .$$

 If there is some non-trivial dyadic interval $[p,q]\subseteq [0,1]\setminus\bigcup _{j=1}^{\infty}\overline{supp}
(h_{j})$, then $F_{[p,q]}\cap H$ is trivial, a contradiction with the assumption. Thus we may assume that
$[0,1]\setminus\bigcup _{j=1}^{\infty}\overline{supp}(h_{j})$ consists of isolated points.\parskip0pt

 Now suppose that there is some number $s$ such that $\bigcup _{i=1} ^{s}\overline{supp}(h_{i}) = [0,1]$
(this happens when for some $s$ the set $[0,1]\setminus\bigcup _{i=1} ^{s}\overline{supp}(h_{i})$ does
not contain an interval). Since for any $i\leq s$ and any $j\leq l_{i}$, the group
$$F_{\overline{supp}(h_{i,j})}\cap\Big\langle h_{i',j'}\ \Big|\ i'\geq 1, j'\leq l_{i'}\Big\rangle$$
is neither trivial nor cyclic, for any pair $(i,j)$ as above there is an element
$$h_{i,j}'\in F_{\overline{supp}(h_{i,j})}\cap\Big\langle h_{i',j'}\ \Big|\ i'\geq 1, j'\leq l_{i'}\Big\rangle ,$$
which does not have a common root with $h_{i,j}$. Enlarging $s$ if necessary we may assume that for all $i
\leq s$ and $j\leq l_{i}$, we have
$$h_{i,j}'\in F_{\overline{supp}(h_{i,j})}\cap\Big\langle h_{i,j}\ \Big|\ i\leq s, j\leq l_{i}\Big\rangle .$$
Since for any $i$ a relation of the form $gh_{i}g^{-1} = h_{i}$ is satisfied in the limit group $G$, there is a
natural number $N$ such that $F$ satisfies $g_{n}h_{i}g_{n}^{-1}=_{F} h_{i}$ when $i\leq s$ and $n>N$.
On the other hand it follows from the definition of defragmentation and Remark \ref{GS} that for any $f,f'\in F$
and any element $f_{0}$ taken from the defragmentation of $f$, if $[f,f']=\emph{id}$, then $[f_{0},f']=\emph{id}$.
Thus we also have that for any $i\leq s$ and any $j\leq l_{i}$, the relations of the form $g_{n}h_{i,j}g_{n}^{-1}=
_{F} h_{i,j}$ and $g_{n}h'_{i,j}g_{n}^{-1}=_{F} h'_{i,j}$ are satisfied for all $n>N$. We claim that this can happen
only if $g_{n}=\emph{id}$, $n>N$ (contradicting the assumption $g_{n}\to g\neq\emph{id}$).\parskip0pt

 Indeed, consider the defragmentation of $g_{n}$ for some fixed $n>N$. Suppose that it is non-empty
and denote by $g_{n,0}$ a non-trivial element of this defragmentation.  Since $g_{n}$ commutes with
each $h_{i,j}$, for all $i\leq s$ and $j\leq l_{i}$, the element $g_{n,0}$ stablizes the ends of any interval
$\overline{supp}(h_{i,j})$. This fact together with the equality $\bigcup _{j=1} ^{s}\overline{supp}(h_{j})= [0,
1]$, imply that there is a pair of elements $h_{i_{0}, j_{0}}$ and $h_{i_{0},j_{0}}'$ as above in $F
_{\overline{supp}(g_{n,0})}$, which do not have a common root and such that
$$[g_{n,0},h_{i_{0}, j_{0}}]=[g_{n,0},h_{i_{0},j_{0}}']=\emph{id} .$$
Once again it follows from Remark \ref{GS} that
$$C_{F_{\overline{supp}(h_{i_{0},j_{0}})}}(h_{i_{0},j_{0}})=\Big\langle\hat{h}_{i_{0},j_{0}}\Big\rangle\ \ \land\ \ C_{F
	_{\overline{supp}(h'_{i_{0},j_{0}})}}(h_{i_{0},j_{0}})=\Big\langle\hat{h}'_{i_{0},j_{0}}\Big\rangle$$
for some roots $\hat{h}_{i_{0},j_{0}}$ and $\hat{h}'_{i_{0}, j_{0}}$ of $h_{i_{0},j_{0}}$ and $h'_{i_{0},j_{0}}$
respectively. Since $g_{n,0}\in\langle\hat{h}_{i_{0},j_{0}}\rangle\cap\langle\hat{h}'_{i_{0},j_{0}}\rangle$, we see
that $g_{n,0}=\emph{id}$, a contradiction.\parskip0pt

 Now consider the case when for any $s\in\mathbb{N}$, $\bigcup _{i=1} ^{s}\overline{supp}(h_{i})\neq [0,
1]$. For any $i\geq 1$, $j\leq l_{i}$, denote by $[p_{i,j},q_{i,j}]$ the support of $h_{i,j}$. It follows from the
assumptions of the theorem that some initial and some final subintervals of $[0,1]$ are covered by the
supports of two elements from $H$. Thus there is some dyadic non-trivial interval $[u,v]$, $u\neq 0$, $v
\neq 1$, such that for any $s\in\mathbb{N}$,
$$\bigcup _{i=1}^{s}\Big(\bigcup _{j=1}^{l_{i}} [p_{i,j},q_{i,j}]\Big)\neq [u,v] .$$

Since $[0,1]\setminus\bigcup _{j=1}^{\infty}\overline{supp}(h_{j})$ consists of isolated points we may
assume that there are some two elements $h_{r_{1}}$ and $h_{r_{2}}$ such that $u\in\overline{supp}(h_{r
_{1}})$ and $v\in\overline{supp}(h_{r_{2}})$. Let $h_{r_{1},j_{1}}$ and $h_{r_{2},j_{2}}$, where $j_{1}\leq l
_{r_{1}}$, $j_{2}\leq l_{r_{2}}$, be such two elements of the decompositions of $h_{r_{1}}$ and $h_{r
_{2}}$ respectively, that $u\in\overline{supp}(h_{r_{1},j_{1}})$ and $v\in\overline{supp}(h_{r_{2},j_{2}})$.
Now
$$[u,v]\nsubseteq\overline{supp}(h_{r_{1}})\cup\overline{supp}(h_{r_{2}})\ \Longrightarrow\ \overline{supp}
	(h_{r_{1},j_{1}})\subseteq\Big( [0,1]\setminus\overline{supp}(h_{r_{2},j_{2}})\Big) .$$

 Since the groups
$$\Big\langle h_{i,j}\ \Big|\ i\geq 1, j\leq l_{i}\Big\rangle\cap F_{[p_{r_{1},j_{1}},q_{r_{1},j_{1}}]}$$
and
$$\Big\langle h_{i,j}\ \Big|\ i\geq 1, j\leq l_{i}\Big\rangle\cap F_{[p_{r_{2},j_{2}},q_{r_{2},j_{2}}]}$$
are neither trivial nor cyclic, our standard argument shows that there is some $N$ such that for all $n>N$
$$g_{n}|_{[p_{r_{1},j_{1}},q_{r_{1},j_{1}}]\cup [p_{r_{2},j_{2}},q_{r_{2},j_{2}}]}=\emph{id}|_{[p_{r_{1,j_{1}}},q_{r_{1},
	j_{1}}]\cup [p_{r_{2,j_{2}}},q_{r_{2},j_{2}}]}.$$
Since $g_{n}$ does not move any points within $[p_{r_{1},j_{1}},q_{r_{1},j_{1}}]\cup [p_{r_{2},j_{2}},q_{r
_{2},j_{2}}]$ and is increasing, we also have that $g_{n}([0,p_{r_{1},j_{1}}])=[0,p_{r_{1},j_{1}}]$, $g_{n}([q
_{r_{1},j_{1}},p_{r_{2},j_{2}}])=[q_{r_{1},j_{1}},p_{r_{2},j_{2}}]$ and $g_{n}([q_{r_{2},j_{2}},1])=[q_{r_{2},j
_{2}},1]$ for all $n>N$.\parskip0pt

 Let $f_{0}\in F_{[p_{r_{1},j_{1}},q_{r_{2},j_{2}}]}$ be any element satisfying $f_{0}(q_{r_{1},j_{1}})=
p_{r_{2},j_{2}}$. Consider the word $w_{1}=[y,y^{f_{0}}]$, where $y$ denotes a variable. Fix some $n
>N$. For any $t\in [0,p_{r_{1},j_{1}}]$ we have:
$$w_{1}(g_{n})(t)=g_{n}^{-1}f_{0}^{-1}g_{n}^{-1}f_{0}g_{n}f_{0}^{-1}g_{n}f_{0}(t)=g_{n}^{-1}f_{0}^{-1}g_{n}
	^{-1}f_{0}g_{n}f_{0}^{-1}g_{n}(t).$$
But $g_{n}(t)\in [0,p_{r_{1},j_{1}}]$ and thus we obtain the further reductions of $w_{1}$:
$$w_{1}(g_{n})(t)=g_{n}^{-1}f_{0}^{-1}g_{n}^{-1}f_{0}g_{n}f_{0}^{-1}(g_{n}(t))=g_{n}^{-1}f_{0}^{-1}g_{n}^{-1}
	(g_{n}^{2}(t))=t.$$
Let $t\in [p_{r_{1},j_{1}},q_{r_{1},j_{1}}]$. Since $p_{r_{2},j_{2}}=f_{0}(q_{r_{1},j_{1}})$, we have that $g_{n}
(f_{0}(t))\leq f_{0}(q_{r_{1},j_{1}})$. Hence $f_{0}^{-1}g_{n}f_{0}(t)\leq q_{r_{1},j_{1}}$. On the other hand
since
$$f_{0}([0,p_{r_{1},j_{1}}])=g_{n}([0,p_{r_{1},j_{1}}])=[0,p_{r_{1},j_{1}}] ,$$
we also have that $f_{0}^{-1}g_{n}f_{0}(t)\geq p_{r_{1},j_{1}}$. Thus it follows that $f_{0}^{-1}g_{n}f_{0}(t)\in
[p_{r_{1},j_{1}},q_{r_{1},j_{1}}]$ and we have:
$$w_{1}(g_{n})(t)=g_{n}^{-1}f_{0}^{-1}g_{n}^{-1}f_{0}g_{n}(f_{0}^{-1}g_{n}f_{0}(t))=g_{n}^{-1}f_{0}^{-1}g_{n}
	^{-1}f_{0}(f_{0}^{-1}g_{n}f_{0}(t))=g_{n}^{-1}(t).$$
But $t$ is in $[p_{r_{1},j_{1}},q_{r_{1},j_{1}}]$ and hence $g_{n}^{-1}(t)=t$.\parskip0pt

Let $t\in [q_{r_{1},j_{1}},q_{r_{2},j_{2}}]$. Since $f_{0}(t)\geq p_{r_{2},j_{2}}$, we have:
$$w_{1}(g_{n})(t)=g_{n}^{-1}f_{0}^{-1}g_{n}^{-1}f_{0}g_{n}f_{0}^{-1}g_{n}f_{0}(t)=g_{n}^{-1}f_{0}^{-1}g_{n}
	^{-1}f_{0}g_{n}f_{0}^{-1}f_{0}(t)=$$
$$=g_{n}^{-1}f_{0}^{-1}g_{n}^{-1}(f_{0}g_{n}(t))=g_{n}^{-1}f_{0}^{-1}(f_{0}g_{n}(t))=t.$$
Finally for $t\in [q_{r_{2},j_{2}},1]$ we apply the same argument as to $z\in [0,p_{r_{1},j_{1}}]$.\parskip0pt

 We see that for any $n>N$, $w_{1}(g_{n})=\emph{id}$. Thus $G\models w_{1}(g)=\emph{id}$. On the other
hand since
$$[u,v]\subseteq [p_{r_{1},j_{1}},q_{r_{2},j_{2}}]\subseteq\Big( [p_{r_{1},j_{1}},q_{r_{1},j_{1}}]\cup\overline{supp}
	(f_{0})\cup [p_{r_{2},j_{1}},q_{r_{2},j_{2}}]\Big)$$
and
$$[u,v]\neq\bigcup _{i=1}^{s}\Big(\bigcup _{j=1}^{l_{i}} [p_{i,j},q_{i,j}]\Big) ,$$
the element $f_{0}$ does not belong to $H$. Thus by Britton's Lemma $w_{1}$ cannot be reduced in the
HNN-extension $F\ast _{H}$. We obtain a contradiction with $w_{1}(g)=\emph{id}$. This finishes the proof
of the claim.\\

\ \ \ \ \ \ \ \ \ \ \ \ \ \ \ \ \ \ \ \ \ \ \ \ \ \ \ \ \ \ \ \ \ \ \ \ \ \ \ \ \ \ \ \ \ \ \ \ \ \ \ \ \ \ \ \ \ \ \ \ \ \ \ \ \ \ \ \ \ \ \ \ \ \ \ \ \ $\square$\\

 By claim there is some dyadic non-trivial interval $[p,q]\subseteq [0,1]$ such that $F_{[p,q]}\cap\langle h_{i,
j}\ |\ i\geq 1, j\leq l_{i}\rangle$ is trivial or cyclic.\parskip0pt

Denote by $h$ the generator of $F_{[p,q]}\cap\langle h_{i,j}\ |\ i\geq 1, j\leq l_{i}\rangle$ if it is cyclic (if
$F_{[p,q]}\cap\langle h_{i,j}\ |\ i\geq 1, j\leq l_{i}\rangle=\{ \emph{id}\}$ we consider $h:=\emph{id}$). If $F
_{[p,q]}\cap\langle h_{i,j}\ |\ i\geq 1, j\leq l_{i}\rangle$ is non-trivial, then by shortening the interval $[p,q]$
if necessary, we may assume that $\overline{supp}(h)=[p,q]$. We now find two pairwise disjoint closed
and dyadic subintervals $I_{i}=[a_{i},b_{i}]\subseteq [p,q]$, $i=1, 2$, $a_{1}<a_{2}$. Next we choose two
non-trivial elements $f_{1}\in F_{I_{1}}$ and $f_{2}\in F_{I_{2}}$. It follows from Remark \ref{GS} that $f
_{1}$ and $f_{2}$ do not have a common root with $h$ (except the trivial case when $h=\emph{id}$). Then
it follows from Proposition \ref{lwc2} that the word
$$w_{2}:=[y^{-1}f_{1}^{-1}yf_{2}^{-1}y^{-1}f_{1}yf_{2},yf_{1}^{-1}y^{-1}f_{2}^{-1}yf_{1}y^{-1}f_{2}]$$
satisfies $w_{2}(f)=\emph{id}$ for all $f\in F$. This implies that for all $g_{n}$, $w_{2}(g_{n})=\emph{id}$
and thus $G\models w_{2}(g)=\emph{id}$.\parskip0pt

 On the other hand $w_{2}$ is non-trivial, non-constant and contains two types of constants $f_{1}^{\pm
1}$ and $f_{2}^{\pm 1}$. Since $H<\langle h_{i,j}\ |\ i\geq 1, j\leq l_{i}\rangle$, $f_{1},f_{2}\notin H$. Thus
it once again follows from Britton's Lemma on irreducible words in HNN-extensions that the word:
$$w_{2}:=[g^{-1}f_{1}^{-1}gf_{2}^{-1}g^{-1}f_{1}gf_{2},gf_{1}^{-1}g^{-1}f_{2}^{-1}gf_{1}g^{-1}f_{2}]$$
is non-trivial in $G=F\ast_{H}$. This contradiction finishes the proof.\\

\ \ \ \ \ \ \ \ \ \ \ \ \ \ \ \ \ \ \ \ \ \ \ \ \ \ \ \ \ \ \ \ \ \ \ \ \ \ \ \ \ \ \ \ \ \ \ \ \ \ \ \ \ \ \ \ \ \ \ \ \ \ \ \ \ \ \ \ \ \ \ \ \ \ \ \ \ $\square$\\

Now we can proceed with the second step of the proof of Theorem \ref{mt}.\\

\emph{Proof of Theorem \ref{mt}.}\ To obtain a contradiction assume that $G=F\ast _{\alpha}$. Let $\{ h
_{1},\ldots , h_{s}\}$ be the generating set of $H$. First consider the centralized case, i.e. $\alpha = id$.
If $\bigcup _{i=1} ^{s}\overline{supp}(h_{i}) = [0,1]$, then condition $(\diamondsuit )$ holds and we
apply Theorem \ref{dmt}. If $\bigcup _{i=1} ^{s}\overline{supp}(h_{i})\neq [0,1]$, then there exists some
dyadic non-trivial subinterval $[p,q]\subseteq [0,1]$ such that $F_{[p,q]}\cap\langle h_{i,j}\ |\ i\leq s, j\leq l_{i}
\rangle =\{\emph{id}\}$, where $h_{i,j}$ are elements of defragmentations of generators of $H$. Hence we refer
to the proof of Theorem \ref{dmt} and construct a suitable law with constants. The rest of the the proof in
this case is exactly as in the final part of the proof above.\parskip0pt

 Generally, let us consider the situation, where $G=F\ast _{\alpha}$ and $\alpha$ is an arbitrary
embedding. Consider a sequence $(g_{n})_{n<\omega}$ in $F$ such that for every $i\leq s$, $g_{n}h
_{i}g_{n}^{-1}=\alpha (h_{i})$ holds for almost all $n$. It follows that there is at least one element $f\in F$
such that for all $i\leq s$, $fh_{i}f^{-1}=\alpha (h_{i})$. Fix such $f$. We now apply the argument of
Theorem \ref{dmt} to the sequence $(f^{-1}g_{n})_{n<\omega}$. This sequence is convergent to some
element $g'$ (in fact $g'=f^{-1}g$). Each $f^{-1}g_{n}$ commutes with any $h_{i}$, $i\leq s$, and hence
with any $h$ from $H$. We can now repeat step by step the proof above to get a contradiction.\\

\ \ \ \ \ \ \ \ \ \ \ \ \ \ \ \ \ \ \ \ \ \ \ \ \ \ \ \ \ \ \ \ \ \ \ \ \ \ \ \ \ \ \ \ \ \ \ \ \ \ \ \ \ \ \ \ \ \ \ \ \ \ \ \ \ \ \ \ \ \ \ \ \ \ \ \ \ $\square$\\

\vspace*{10mm}

\begin{flushleft}
\begin{footnotesize}
Roland Zarzycki\\
University of Wroclaw\\
pl. Grunwaldzki 2/4\\
50-384 Wroc{\l}aw\\
Poland\\
email: zarzycki@math.uni.wroc.pl
\end{footnotesize}
\end{flushleft}


\begin{thebibliography}{99}

	\bibitem{A} M. Abert, Group laws and free subgroups in topological groups, \emph{Bull.
		London Math. Soc.} 37 (2005), 525-534.
	\bibitem{Ast} A. Akhmedov, M. Stein, J. Taback, Free limits of Thompson's group $F$,
		announcement, conference on Geometric and Asymptotic Group Theory with
		Applications, March 9-12, 2009, Stevens Institute, Hoboken, NJ.
	\bibitem{Gr} L. Bartholdi, R.I. Grigorchuk, Z. Sunik, Branch groups, [in:] \emph{Handbook
		of Algebra}, vol. 3, North-Holland, Amsterdam, 2003, 989-1112.
	\bibitem{BaB} B. Baumslag, Residually free groups, \emph{Proc. London Math. Soc.} (3),
		17:402-418, 1967.
	\bibitem{BaG} G. Baumslag, On generalised free products, \emph{Math. Z.}, 78:423-438,
		1962.
	\bibitem{Bri} M. G. Brin, The free group of rank 2 is a limit of Thompson's group F, eprint
		arXiv:0904.2626v2.
	\bibitem{BS} M. G. Brin, C. C. Squier, Groups of piecewise linear homeomorphisms of the
		real line, \emph{Invent. Math.} 79 (1985), 485-498.
	\bibitem{CFP} J.W. Cannon, W.J. Floyd, W.R. Parry, Introductory notes on Richard
		Thompson's groups, \emph{Enseign. Math.} (2) 42 (1996), 215-256.
	\bibitem{GC} C. Champetier, V. Guirardel, Limit groups as limits of free groups:
		compactifying the set of free groups, \emph{Israel J. Math.} 146 (2005), 1-76.
	\bibitem{EE} E. S. Esyp, On identities in Thompson's group, eprint arXiv:0902.0199v1.
	\bibitem{FGM} B. Fine, A. M. Gaglione, A. Myasnikov, G. Rosenberger, D. Spellman, A
		classification of fully residually free groups of rank three or less, \emph{J.
		Algebra}, 200(2):571-605, 1998.
	\bibitem{Gri} R. Grigorchuk, Just Infinite branch groups, [in:] \emph{New Horizons in pro-$p$
		Groups}, M. Du Santoy, D. Segal, A. Shalev (eds), Progress in Mathematics, Vol.
		184, Birkh\"{a}user, Boston (2000), 121-180.
	\bibitem{GS} V. Guba, M. Sapir, Diagram groups, \emph{Memoirs of the American
		Mathematical Society}, Volume 130, Number 620, November 1997.
	\bibitem{G} L. Guyot, Limits of dihedral groups, eprint arXiv:math/0710.1495.
	\bibitem{SG} L. Guyot and Y. Stalder, Limits of Baumslag-Solitar groups, \emph{Groups
		Geom. Dyn.} 2 (2008), 353-381.
	\bibitem{H} M. Hall, The Theory of Groups, The McMillan Co., New York, 1965.
	\bibitem{KM} M. Kassabov, F. Matucci, The simultaneous conjugacy problem in
		Thompson's group F, eprint arXiv:math/0607167.
	\bibitem{KhM} O. Kharlampovich, A. Myasnikov, Irreducible affine varieties over a free
		group. II: Systems in row-echelon form and description of residually free groups,
		\emph{J. Algebra}, V. 200, 517--570 (1998).
	\bibitem{LS} R. Lyndon, P. Schupp, Combinatorial group theory, Springer, Berlin 1977.
	\bibitem{Sel} Z. Sela, Diophantine geometry over groups I: Makanin-Razborov diagrams,
		\emph{Publications Mathematiques de l'IHES} 93(2001), 31-105.
	\bibitem{Zar} R. Zarzycki, Limits of Thompson's group $F$, to appear in \emph{Trends
		in  Mathematics}, Birkhauser Verlag.

\end{thebibliography}
\end{document}